\documentclass[11pt,epsfig,amsfonts]{amsart}
\usepackage{stmaryrd}
\usepackage{xcolor}
\usepackage{epsfig}
\usepackage{enumerate}
\usepackage{amsmath}
\usepackage{amssymb}
\usepackage{amscd}
\usepackage[normalem]{ulem}
\usepackage{graphicx}
\usepackage{lineno}
\usepackage{pstricks}
\usepackage{tocvsec2}
\usepackage{mathtools}
\allowdisplaybreaks[4]
\usepackage{setspace}

\usepackage [backref, colorlinks, citecolor=red, anchorcolor=black, linkcolor=blue]{hyperref}

\usepackage{amsfonts}
\usepackage[top=35mm, bottom=35mm, left=30mm, right=30mm]{geometry}
\usepackage{verbatim}
\usepackage{graphicx, amssymb, amsmath,amsthm}
\allowdisplaybreaks[4]
\usepackage{appendix}
 \usepackage{mathrsfs}
 \usepackage{cite}

\usepackage{setspace}
 \makeatletter
 \@namedef{subjclassname@2020}{%
 	\textup{2020} Mathematics Subject Classification}
 \makeatother
\usepackage{xcolor}

\usepackage{xcolor}

\newtheorem{definition}{Definition}[section]
\newtheorem{lemma}{Lemma}[section]
\newtheorem{theorem}{Theorem}%[section]
\newtheorem{corollary}{Corollary}[section]
\newtheorem{remark}{Remark}[section]

\newtheorem{proposition}{Proposition}[section]

\numberwithin{equation}{section}
\usepackage{epsfig}
\newtheorem{claim}{Claim}[section]

\numberwithin{equation}{section}
\renewcommand*{\backref}[1]{}
\renewcommand*{\backrefalt}[4]{\quad \tiny
  \ifcase #1 (\textbf{NOT CITED.})%
  \or    (Cited on Section~#2.)%
  \else   (Cited on Section~#2.)%
  \fi}

\def \P{\mathbb{P}}

\def \A{\mathbb{A}}

    \def \E{\mathcal{E}}

\makeatletter

\def\MRbibitem{\@ifnextchar[\my@lbibitem\my@bibitem}

\def\mybiblabel#1#2{\@biblabel{{\hyperref{http://www.ams.org/mathscinet-getitem?mr=#1}{}{}{#2}}}}

\def\myhyperanchor#1{\Hy@raisedlink{\hyper@anchorstart{cite.#1}\hyper@anchorend}}

\def\my@lbibitem[#1]#2#3#4\par{%
  \item[\mybiblabel{#2}{#1}\myhyperanchor{#3}\hfill]#4%
  \@ifundefined{ifbackrefparscan}{}{\BR@backref{#3}}%
  \if@filesw{\let\protect\noexpand\immediate% write to aux-file
    \write\@auxout{\string\bibcite{#3}{#1}}}\fi\ignorespaces%
}

\def\my@bibitem#1#2#3\par{%
  \refstepcounter\@listctr% standard tex item counter for the generic item number
  \item[\mybiblabel{#1}{\the\value\@listctr}\myhyperanchor{#2}\hfill]#3%
  \@ifundefined{ifbackrefparscan}{}{\BR@backref{#2}}%
  \if@filesw\immediate\write\@auxout% write to aux-file
    {\string\bibcite{#2}{\the\value\@listctr}}\fi\ignorespaces%
}

\makeatother


\author{Xue Liu}
\address[Xue Liu]
{School of Mathematics\\ Southeast University\\
 Nanjing 211189, PR China}
\email[X.~Liu]{xueliuseu@seu.edu.cn}

\author{Xiaomin Zhou}
\address[Xiaomin Zhou]
{School of Mathematics and Statistics, Huazhong University of Science and Technology, Wuhan, Hubei 430074, China}
\email[X.~Zhou]{zxm12@mail.ustc.edu.cn}

\author{Xiaoyao Zhou}
\address[Xiaoyao Zhou]
{School of Mathematical Sciences, Ministry of Education Key Laboratory of NSLSCS, Nanjing Normal University, Nanjing, 210023, Jiangsu, P.R.China}
\email[X.~Zhou]{zhouxiaoyaodeyouxian@126.com}

\begin{document}

\begin{abstract}
We establish quenched level-2 large deviation estimates for empirical measures of measurable skew products with compact metric fibers over an invertible ergodic measure-preserving transformation on a Polish space. Our hypotheses are formulated through lower and upper energy functions, which describe the exponential decay of conditional measures of fiber Bowen balls. Under almost-sure fiber specification, we prove lower bounds for open sets. Upper bounds for closed convex sets require only the upper energy condition, and upper semicontinuity of fiber entropy removes the convexity restriction. When these hypotheses hold with a common energy function \(T\), we obtain a quenched level-2 large deviation principle with a good rate function \(I_T(\rho)=\int T\,d\rho-h_\rho(F)\) on invariant measures with the prescribed base marginal, and \(+\infty\) otherwise. We also derive quenched level-1 large deviation principles for bounded continuous observables and a quenched Laplace principle. As an application, we consider skew products that are Anosov and topologically mixing on fibers. For Gibbs equilibrium states associated with random potentials $f$ having integrable two-sided Bowen variation, we construct bounded strictly positive common energy functions by coboundary corrections. This construction identifies the large deviation rate on invariant measures as \(\pi_F(f)-h_\rho(F)-\int f\,d\rho\), where \(\pi_F(f)\) denotes the fiber topological pressure.
\end{abstract}

\title[Quenched large deviations for skew products]{Quenched large deviation estimates for skew product systems with the almost-sure fiber specification property}
\maketitle

\parskip 10pt
\begingroup
  \setcounter{tocdepth}{2} % Include sections and subsections
  \setstretch{0.9}
  \small
  \setlength{\parskip}{0pt}
  \tableofcontents
\endgroup
\section{Introduction}\label{section 1}

\subsection{Quenched level--2 large deviations for random dynamical systems}
Ergodic theorems describe the typical asymptotic behavior of trajectories, leaving open the quantitative question of how likely substantial statistical deviations are over long  observation times.
Large deviation theory addresses this question by quantifying the asymptotic decay of these probabilities on an exponential scale. In this paper, we study this problem for randomly driven systems, where the evolution of the environment determines the sequence of maps acting on the phase space.
Let $M$ be a compact metric space and let $(\Omega,\mathcal{F},\P,\theta)$ be an invertible ergodic measure-preserving dynamical system on a Polish space $\Omega$ describing a general driving system. Let $\{F_\omega\}_{\omega\in\Omega}$  be a family of homeomorphisms on $M$ so that the  skew product system $\Theta(\omega,x)=(\theta\omega,F_\omega x)$ is measurable. This skew product defines a two-sided continuous random dynamical system (RDS)
\begin{equation*}
  F_\omega^n=\begin{cases}
               F_{\theta^{n-1}\omega}\circ \cdots \circ F_\omega, & \mbox{if } n\geq 1  \\
               \operatorname{id}, & \mbox{if } n=0 \\
               (F_{\theta^n\omega})^{-1} \circ \cdots \circ(F_{\theta^{-1}\omega})^{-1}, & \mbox{if }n<0.
             \end{cases}
\end{equation*}
For standard references for RDSs, see Arnold's book \cite{ArnoldBook1998} and the survey by Kifer and Liu \cite{KL06}.

Following the quenched viewpoint in \cite{Kifer1996RandomEnvironments}, we use the following formulation of the level--2 large deviation principle.
Let \(\nu\) be a reference probability measure on $\Omega\times M$ with marginal $\P$ on $\Omega$, and we denote its disintegration by $\{\nu_\omega\}_{\omega\in\Omega}$. Let \(\operatorname{Pr}(\Omega\times M)\) denote the space of probability measures on $\Omega\times M$ endowed with the narrow topology. The random dynamical system $F$ (or the skew product $\Theta$) is said to satisfy a quenched level--2 large deviation principle with a reference measure $\nu$ and a rate function $I:\operatorname{Pr}(\Omega\times M)\to[0,+\infty]$ if  there exists a $\P$-full-measure set $\Omega_0$ such that for every closed subset $K\subset \operatorname{Pr}(\Omega\times M)$,
\begin{equation*}
   \limsup_{n\to+\infty}\frac{1}{n}\log\nu_\omega\{x\in M:\ \E_n(\omega,x)\in K\}\leq -\inf_{\rho\in K}I(\rho),\ \forall \omega\in\Omega_0
\end{equation*}
and for every open subset $G\subset \operatorname{Pr}(\Omega\times M)$,
\begin{equation*}
  \liminf_{n\to +\infty}\frac{1}{n}\log \nu_\omega\{x\in M:\ \mathcal{E}_n(\omega,x)\in G\}\geq -\inf_{\rho\in G}I(\rho),\ \forall \omega\in\Omega_0,
\end{equation*}
with \(\inf\varnothing=+\infty\) and \(\log0=-\infty\), where the empirical measure is defined by
\begin{equation*}
  \E_n(\omega,x)
:=\frac1n\sum_{j=0}^{n-1}\delta_{\Theta^j(\omega,x)}=\frac1n\sum_{j=0}^{n-1}\delta_{(\theta^j\omega,F_\omega^j x)}.
\end{equation*}
 Empirical measures provide a natural level of description because they record the statistical distribution of an orbit and simultaneously determine the time averages of all bounded continuous observables. Establishing a quenched level--2 large deviation principle therefore quantifies, on the exponential scale, the likelihood of atypical statistical behavior and the reliability of asymptotic ergodic predictions. It also provides a common framework for studying fluctuations of different observables, since the level--1 large deviation principles for their time averages follow from the corresponding principle for empirical measures.

  When the external driving is absent, represented by taking $\Omega$ to be a
singleton, the quenched
level--2 large deviation principle reduces to its deterministic counterpart.
Several complementary methods have been developed for this problem and for
the related deviations of time averages.

Thermodynamic formalism relates large deviations to pressure and
equilibrium states. Kifer \cite{Kifer1990} developed a general criterion
based on limiting logarithmic moment generating functionals and suitable
uniqueness conditions for equilibrium states. This approach also yields
level--2 principles for rational maps satisfying the topological
Collet--Eckmann condition \cite{CommanRiveraLetelier2011}. Spectral
methods study perturbations of transfer operators, as in Keller and
Nowicki's work on Collet--Eckmann unimodal maps \cite{KellerNowicki1992}.
Inducing and tower constructions provide another route for nonuniformly
hyperbolic systems, with exponential return-time tails yielding large
deviation principles for ergodic averages \cite{ReyBelletYoung2008}.

These approaches have random counterparts. Kifer
\cite{Kifer1991RandomExpanding} and Liu and Zhao \cite{LiuZhao2003}
obtained quenched large deviation results for random expanding maps and
small random perturbations of Axiom~A systems, respectively. Spectral
methods have been extended to random expanding and hyperbolic systems
through transfer-operator cocycles and their leading Lyapunov exponents
\cite{DragicevicFroylandGonzalezTokmanVaienti2018,
DragicevicFroylandGonzalezTokmanVaienti2020}. For i.i.d. random systems
admitting Young towers with exponential tails, Hafouta \cite{Hafouta2022}
established quenched local large deviation principles for observables.

The method most directly related to the present work combines entropy
estimates with control of the masses of Bowen balls and an orbit-gluing
property. Young \cite{Young1990} established general large deviation bounds
of this type, with applications to Axiom~A attractors and Gibbs measures on
shift spaces. The role of specification is particularly clear in the lower
bound: it allows orbit segments typical for different ergodic measures to be
concatenated, so that nonergodic invariant measures can also contribute to the
variational estimate. Entropy controls the number of distinguishable orbit
segments, while the Bowen-ball estimates control their mass under the
reference measure.

Classical specification can be weakened substantially in this argument.
Pfister and Sullivan \cite{PfisterSullivan2005} introduced the approximate
product property and proved that it implies entropy density of ergodic
measures. Together with suitable energy-function hypotheses, their framework
gives level--2 large deviation bounds and, in particular, a full principle for
every $\beta$-shift with respect to its measure of maximal entropy.
Climenhaga, Thompson and Yamamoto \cite{ClimenhagaThompsonYamamoto2017}
established level--2 principles using specification on a collection of words,
edit approachability of the language, and a corresponding weak Gibbs
property. Their applications include $\beta$-shifts, $S$-gap shifts, and their
factors, with suitable equilibrium states as reference measures. Pfister and
Sullivan \cite{PfisterSullivan2018} further clarified the weak Gibbs
formulation: entropy density supplies the lower bound for all invariant
measures, and upper semicontinuity of entropy supplies the upper bound for
arbitrary closed sets.

Specification also has a history in random thermodynamic formalism.
Gundlach and Kifer \cite{Kifer2000specification} introduced specification and
expansiveness for random bundle transformations, proved uniqueness of
equilibrium states under their hypotheses, and discussed large deviations
for random sofic shifts. In their study of Anosov systems driven by a quasi-periodic force, Huang, Lian, and Lu \cite{HLLScichina} introduced a random specification
property. Its fiber realization is the fiber specification property
formulated by Nian Liu and Xue Liu
\cite{Liu2025variatonalprinciple,LiuLiu2026Irregular} in their work on
the multifractal analysis of skew product transformations.

In the present paper, we develop an energy-function approach to quenched
large deviations for measurable skew products. Proposition~\ref{proposition lower energy function}
characterizes the existence of lower and nondegenerate upper energy functions in terms of
conditional fiber Bowen-ball masses. Theorem~\ref{thm 2} gives upper bounds for closed
convex sets without specification, while Theorem~\ref{thm main 1} gives lower bounds
under $\P$-almost-sure fiber specification. If fiber entropy is upper semicontinuous
and one function satisfies both energy conditions, Theorem~\ref{proposition normalized good rate function}
gives a quenched level--2 large deviation principle with a good rate function equal
to energy minus fiber entropy. We also obtain level--1 large deviation principles
and a quenched Laplace principle.

We apply this framework to skew products that are Anosov and
topologically mixing on fibers in Section \ref{section example}. For Gibbs equilibrium states associated
with random potentials with integrable two-sided Bowen variation, we
construct bounded positive common energy functions and identify the
resulting rate functions in terms of fiber pressure, entropy, and the
potential.

\subsection{Setting and main results}
Throughout this paper, let $\Omega$ be a Polish space, let $\P$ be a Borel probability measure on $\Omega$, and let $\mathcal{F}$ be the Borel $\sigma$-algebra.
Let $\theta$ be an invertible bimeasurable map that preserves $\P$ and is ergodic. Let $M$ be a compact metric space with metric $d_M$, and let $F:\Omega\to\operatorname{Homeo}(M,M)$ be a Borel measurable map, where $\operatorname{Homeo}(M,M)$ is the space of homeomorphisms from $M$ to itself endowed with the compact-open topology. Then the generator $\{F_\omega\}_{\omega\in\Omega}$ defines a skew-product system $\Theta: \Omega\times M \to \Omega\times M$ by $\Theta(\omega,x)=(\theta\omega, F_{\omega}(x))$. We note that $\Theta$ may not be continuous.

 Let $\operatorname{Pr}(\Omega\times M)$ denote the space of Borel probability measures on the product space $(\Omega\times M,\mathcal{F}\otimes \mathcal{B}(M))$. The narrow topology $\tau$ on $\operatorname{Pr}(\Omega\times M)$ is the coarsest topology such that $\mu\mapsto \int_{\Omega\times M}g\,d\mu$ is continuous for every bounded continuous function $g\in C_b(\Omega\times M)$.  Let $\operatorname{Pr}_{\P}(\Omega\times M)$ denote the space of probability measures on the product space $(\Omega\times M,\mathcal{F}\otimes\mathcal{B}(M))$ whose marginal on $\Omega$ equals $\P$.  By Rokhlin's disintegration theorem, any $\mu\in \operatorname{Pr}_\P(\Omega\times M)$  disintegrates $d\mu(\omega,x)=d\mu_\omega(x)d\P(\omega)$, where for $\P$-a.e. $\omega\in\Omega$, $\mu_\omega$ is a Borel probability measure on $M$.  Let $I_\P(\Omega\times M)$ be the set of $\Theta$-invariant Borel probability measures with marginal $\P$ on $\Omega$, and let $I_\P^e(\Omega\times M)$ be its subset of ergodic measures. We write $h_\rho(F)$ for the fiber entropy, equivalently the entropy of $\Theta$ relative to $\mathcal B_\Omega:=\pi_\Omega^{-1}\mathcal F$; its definition is recalled in Section~\ref{section 3.4}.

 \begin{definition}
   An $\mathcal{F}\otimes\mathcal{B}(M)$-measurable function $T:\Omega\times M\to\mathbb{R}$ is called an upper normal integrand if there exists a sequence of random continuous functions $\phi_k\in L^1(\Omega,C(M))$ such that $\phi_k(\omega,x)\searrow T(\omega,x)$ for $\P$-a.s. $\omega\in\Omega$ and all $x\in M$.  Here $L^1(\Omega,C(M))$ is the linear space of all random continuous functions $\phi:\Omega\times M\to \mathbb{R}$ satisfying:
    \begin{enumerate}
      \item for all $x\in M$, the mapping $\omega\mapsto \phi(\omega,x)$ is measurable;
      \item for all $\omega\in \Omega$, the mapping $x\mapsto \phi(\omega,x)$ is continuous and bounded;
      \item $\sup\{|\phi(\cdot,x)|:\ x\in M\}\in L^1(\P)$.
    \end{enumerate}
 \end{definition}
A sufficient criterion and a semicontinuity property for upper normal integrands
are given in Section~\ref{subsection energy function}. In particular, for a
nonnegative upper normal integrand $T$, the map
$\rho\mapsto\int_{\Omega\times M}T\,d\rho$ is upper semicontinuous on
$\operatorname{Pr}_\P(\Omega\times M)$ in the random narrow topology $\tau_{rn}$.
This topology is generated by $\mu\mapsto\mu(f)$ for $f\in L^1(\Omega,C(M))$;
see \cite[Definition~3.15]{Hans02}.

\begin{definition}\label{definition lower energy}
  Let $\nu\in \operatorname{Pr}_\P(\Omega\times M)$. A function $T:\Omega\times M\to\mathbb R_{\geq0}$ is called a lower energy function for $\nu$ if $T$ is an upper normal integrand and, for $\P$-a.s. $\omega\in\Omega$,
  \begin{equation}\label{eq def lower energy}
    \lim_{\epsilon \searrow 0}\liminf_{N\to +\infty}\inf_{x\in M}\left(\frac{1}{N}\log\nu_\omega(B_\omega(x,N,\epsilon))+\int Td(\E_N(\omega,x))\right)\geq 0,
  \end{equation}
where the fiber empirical measure is defined by $\E_n(\omega,x):=\frac{1}{n}\sum_{i=0}^{n-1}\delta_{\Theta^i(\omega,x)}$, and the fiber Bowen ball is defined by $B_\omega(x,n,\epsilon):=\{y\in M:\ d_{\omega}^n(x,y)=\max\limits_{0\leq i<n}d_M(F_\omega^ix,F_\omega^iy)<\epsilon\}.$

    The function $T:\Omega\times M\to \mathbb{R}_{\geq 0}$ is called an upper energy function for $\nu$ if $(-T)$ is an upper normal integrand, $\operatorname*{ess\,sup}\limits_{\omega}\sup_{x\in M} T(\omega,x)<\infty$, and for $\P$-a.s. $\omega\in\Omega$,
  \begin{equation}\label{eq def upper energy}
     \lim_{\epsilon \searrow 0}\limsup_{N\to +\infty}\sup_{x\in M}\left(\frac{1}{N}\log\nu_\omega(B_\omega(x,N,\epsilon))+\int Td(\E_N(\omega,x))\right)\leq 0.
  \end{equation}
\end{definition}
It is clear that $T\equiv 0$ is an upper energy function for any $\nu\in \operatorname{Pr}_\P(\Omega\times M)$.   An upper energy function $T$ for $\nu$ is called nondegenerate if $
    \int_\Omega \min_{x\in M}T(\omega,x)d\P(\omega)>0.$ The following proposition characterizes the existence of lower and nondegenerate upper energy functions in terms of the asymptotic decay of the measures of fiber Bowen balls.
\begin{proposition}\label{proposition lower energy function}
Let $\Theta$ be the skew-product system on the compact metric space $M$ driven by $(\Omega,\mathcal{F},\P,\theta)$, and let $\nu\in \operatorname{Pr}_\P(\Omega\times M)$. Then $\nu$ admits a lower energy function if and only if
  \begin{equation}\label{eq proposition exists energy}
    \operatorname*{ess\ inf}_{\omega\in\Omega}\lim_{\epsilon\searrow0}\liminf_{n\to+\infty}\inf_{x\in M}\frac{1}{n}\log \nu_\omega (B_\omega(x,n,\epsilon))>-\infty.
  \end{equation}
The measure
   $\nu\in \operatorname{Pr}_\P(\Omega\times M)$ admits a nondegenerate upper energy function if and only if
  \begin{equation}\label{eq proposition exists upper energy}
    \operatorname*{ess\ sup}_{\omega\in\Omega}\lim_{\epsilon\searrow0}\limsup_{n\to+\infty}\sup_{x\in M}\frac{1}{n}\log \nu_\omega (B_\omega(x,n,\epsilon))<0,
  \end{equation}
  where $ \operatorname{ess\ inf}$ and $ \operatorname{ess\ sup}$ are taken with respect to $\P$.
\end{proposition}
We first establish the following general quenched level--2 large deviation upper bounds for fiber empirical measures. For an energy function $T$, define the extended-real-valued entropy--energy functional $I_T$ on $\operatorname{Pr}(\Omega\times M)$ by
\begin{equation}\label{eq rate function I T}
  I_T(\rho):=
  \begin{cases}
    \displaystyle
    \int_{\Omega\times M}T\,d\rho-h_\rho(F),
       &\rho\in I_\P(\Omega\times M),\\[1.2ex]
    +\infty,
       &\rho\notin I_\P(\Omega\times M).
  \end{cases}
\end{equation}
Under a single energy condition, $I_T$ need not be nonnegative and may take
the value $-\infty$. For example, $I_0(\rho)=-h_\rho(F)$ on invariant measures.
Theorem~\ref{proposition normalized good rate function} gives conditions under
which $I_T$ is a good rate function: it takes values in $[0,+\infty]$, is lower
semicontinuous, and has compact sublevel sets
$\{\rho:I_T(\rho)\leq\alpha\}$ for every finite $\alpha\geq0$.

\begin{theorem}\label{thm 2}
  Let $\nu\in \operatorname{Pr}_\P(\Omega\times M)$ be the reference measure, and let $T:\Omega\times M\to\mathbb{R}_{\geq 0}$ be an upper energy function for $\nu$. There exists a $\P$-full-measure set $\Omega_{\rm up}$ such that, for every $\omega\in\Omega_{\rm up}$ and every closed convex set $K\subset \operatorname{Pr}(\Omega\times M)$,
  \begin{equation}\label{eq thm 2}
    \limsup_{n\to+\infty}\frac{1}{n}\log\nu_\omega\{x\in M:\ \E_n(\omega,x)\in K\}\leq -\inf_{\rho\in K}I_T(\rho).
  \end{equation}
If, furthermore, the fiber entropy map $\mu\mapsto h_\mu(F)$ is upper semicontinuous on $I_\P(\Omega\times M)$ with respect to $\tau_{rn}$ (equivalently, the restriction of $\tau$, by Lemma~\ref{lemma lemma diffrent converges}), then the conclusion also holds when $K\subset \operatorname{Pr}(\Omega\times M)$ is assumed only to be closed.
\end{theorem}
In order to establish the quenched level--2 large deviation lower bounds for fiber empirical measures, we need the following $\P$-a.s. fiber specification property, which is weaker than the fiber specification property considered in \cite{Liu2025variatonalprinciple}.
\begin{definition}\label{def fiber specification}
The skew product $\Theta$ has the $\P$-a.s. fiber specification property
if, for every $\epsilon>0$, there exists $m=m(\epsilon)\in\mathbb N$
such that, for $\P$-a.e. $\omega\in\Omega$, the following holds:
for any $k\geq1$, points $x_1,\ldots,x_k\in M$, and integers
$0\leq a_1\leq b_1<\cdots<a_k\leq b_k$ satisfying
$a_{i+1}>b_i+m$ for $1\leq i<k$, there exists $x\in M$ such that
\begin{equation*}
  \max_{0\leq j\leq b_i-a_i}
  d_M\bigl(F_\omega^{a_i+j}x,F_{\theta^{a_i}\omega}^{\,j}x_i\bigr)
  <\epsilon,\qquad i=1,\ldots,k.
\end{equation*}
If the same condition holds for every $\omega\in\Omega$, we call it the
fiber specification property \cite{Liu2025variatonalprinciple}.
\end{definition}
\begin{theorem}\label{thm main 1}
Assume that the skew-product system $\Theta$, driven by $(\Omega,\mathcal{F},\P,\theta)$ on the compact metric space $M$, has the $\P$-a.s. fiber specification property. Assume the reference measure $\nu \in \operatorname{Pr}_\P(\Omega\times M)$ admits a lower energy function, and let $T:\Omega\times M\to \mathbb{R}_{\geq 0}$ be any lower energy function for $\nu$. Let $\{\nu_\omega\}_{\omega\in\Omega}$ be its disintegration. There exists a $\P$-full-measure set $\Omega_{\rm low}$ such that, for every $\omega\in\Omega_{\rm low}$ and every $\tau$-open set $G\subset \operatorname{Pr}(\Omega\times M)$,
  \begin{equation*}
    \liminf_{n\to +\infty}\frac{1}{n}\log \nu_\omega\{x\in M:\ \mathcal{E}_n(\omega,x)\in G\}\geq -\inf_{\rho\in G}I_T(\rho).
  \end{equation*}
\end{theorem}
The hypothesis of Theorem~\ref{thm main 1} is satisfied by several natural classes of random hyperbolic systems. First, every skew product that is Anosov and topologically mixing on fibers in Section~\ref{section example} has the fiber specification property by \cite[Theorem~A]{Liu2025variatonalprinciple}, and hence has the $\P$-a.s.\ fiber specification property for every ergodic invariant measure $\P$ on the base. In particular, we can apply our results to the systems introduced in \cite[Section 7]{HLLScichina}. More generally, suppose that on a common invariant full-measure set there are
homeomorphisms $h_\omega$ such that
$h_{\theta\omega}\circ F_\omega=\widetilde F_\omega\circ h_\omega$,
and that both $\{h_\omega\}$ and $\{h_\omega^{-1}\}$ are equicontinuous.
The fiber specification property is then preserved: choose a shadowing tolerance
using the common modulus of continuity, transfer the prescribed segments by
$h_\omega$, and pull the shadowing point back by $h_\omega^{-1}$.
See also \cite[Lemma 2.3]{LiuMa2026} for specification and induced random dynamics.  The equicontinuous random conjugacy argument covers, in particular, sufficiently small random perturbations of a topologically mixing Anosov diffeomorphism, since random structural stability gives an equicontinuous random conjugacy between the perturbed system and the direct product; see \cite[Lemma~2.2]{LiuMa2026}.

\begin{theorem}\label{proposition normalized good rate function}
Assume that the skew-product system $\Theta$ has the $\P$-a.s. fiber
specification property and that the fiber entropy map
$\rho\mapsto h_\rho(F)$ is upper semicontinuous on
$I_\P(\Omega\times M)$ with respect to $\tau_{rn}$.  Let $\nu\in \operatorname{Pr}_\P(\Omega\times M)$ be the reference measure and assume
that a single function
$T:\Omega\times M\to\mathbb R_{\geq0}$ is simultaneously a lower
energy function and an upper energy function for $\nu$.
The functional $I_T$ is a good rate function on
$(\operatorname{Pr}(\Omega\times M),\tau)$, and
\begin{equation}\label{eq normalization I T}
  \inf_{\rho\in \operatorname{Pr}(\Omega\times M)}I_T(\rho)=0.
\end{equation}
There exists a $\P$-full-measure set $\Omega_{\mathrm{LD}}$ such that, for every $\omega\in\Omega_{\mathrm{LD}}$, the following bounds hold simultaneously for all open sets $G\subset \operatorname{Pr}(\Omega\times M)$ and all closed sets $K\subset \operatorname{Pr}(\Omega\times M)$:
\begin{equation}\label{eq LDP lower I T}
  \liminf_{n\to+\infty}\frac1n
  \log \nu_\omega\{x\in M:\ \E_n(\omega,x)\in G\}
  \geq-\inf_{\rho\in G}I_T(\rho),
\end{equation}
and, for every closed set $K\subset \operatorname{Pr}(\Omega\times M)$,
\begin{equation}\label{eq LDP upper I T}
  \limsup_{n\to+\infty}\frac1n
  \log \nu_\omega\{x\in M:\ \E_n(\omega,x)\in K\}
  \leq-\inf_{\rho\in K}I_T(\rho).
\end{equation}
\end{theorem}
\begin{corollary}[Quenched level--1 large deviations]
\label{corollary quenched level one LDP}
Assume the hypotheses of
Theorem~\ref{proposition normalized good rate function}, and let
$g\in C_b(\Omega\times M)$. On the same set $\Omega_{\mathrm{LD}}$ as in that theorem, the following quenched level--1 large deviation principle holds. For every $\omega\in\Omega_{\mathrm{LD}}$ and every open set $O\subset\mathbb R$,
\begin{equation}\label{eq level one lower bound}
  \liminf_{n\to+\infty}\frac{1}{n}
  \log\nu_\omega\{x\in M:\ \int_{\Omega\times M} g\,d(\E_n(\omega,x))\in O \}
  \geq-\inf_{a\in O}J_g(a),
\end{equation}
and, for every closed set $K\subset\mathbb R$,
\begin{equation}\label{eq level one upper bound}
  \limsup_{n\to+\infty}\frac{1}{n}
  \log\nu_\omega\{x\in M:\int_{\Omega\times M} g\,d(\E_n(\omega,x))\in K\}
  \leq-\inf_{a\in K}J_g(a),
\end{equation}
where $J_g:\mathbb{R}\to [0,+\infty]$ is a good rate function defined by
\begin{equation}\label{eq level one rate function}
  J_g(a):=
  \inf\left\{
  \int_{\Omega\times M}T\,d\rho-h_\rho(F):
  \rho\in I_\P(\Omega\times M),\
  \int_{\Omega\times M}g\,d\rho=a
  \right\},
\end{equation}
and the infimum over the empty set is understood to be $+\infty$. The set $\Omega_{\mathrm{LD}}$ is independent of $g$, $O$ and $K$.
\end{corollary}

\begin{corollary}[Quenched Laplace principle]
\label{corollary quenched Laplace principle}
Assume the hypotheses of
Theorem~\ref{proposition normalized good rate function}.  For every
$\Phi\in C_b(\operatorname{Pr}(\Omega\times M))$ and every $\omega\in\Omega_{\mathrm{LD}}$,
\begin{equation}\label{eq quenched Laplace principle}
  \lim_{n\to+\infty}\frac{1}{n}
  \log\int_M
  \exp\bigl(n\Phi(\E_n(\omega,x))\bigr)d\nu_\omega(x)=
  \sup_{\rho\in \operatorname{Pr}(\Omega\times M)}
  \{\Phi(\rho)-I_T(\rho)\}.
\end{equation}
In particular, for every $g\in C_b(\Omega\times M)$ and every $\omega\in\Omega_{\mathrm{LD}}$,
\begin{equation*}
  \lim_{n\to+\infty}\frac1n
  \log\int_M
  \exp\left(\sum_{j=0}^{n-1}g(\Theta^j(\omega,x))\right)d\nu_\omega(x)=
  \sup_{\rho\in I_\P(\Omega\times M)}
  \left\{
  h_\rho(F)+\int_{\Omega\times M}g\,d\rho
  -\int_{\Omega\times M}T\,d\rho
  \right\}.
\end{equation*}
\end{corollary}
In Section~\ref{section example}, we make our large deviation estimates explicit for skew products that are Anosov and topologically mixing on fibers. Following Gundlach and Kifer \cite{Kifer2000specification},
we denote by $V_\P^\pm(F)$ the class of random potentials with
integrable two-sided Bowen variation. For every $f\in V_\P^\pm(F)$, the Gundlach--Kifer theory provides a Gibbs equilibrium state $\nu_f$, and we construct a common lower and upper energy function $T_f$ for $\nu_f$, such that
\begin{equation*}
   \int T_f\,d\rho
  =
  \pi_F(f)-\int f\,d\rho,\ \forall \rho\in I_\P(\Omega\times M),
\end{equation*} where $\pi_F(f)$ is the fiber topological pressure of $F$ with respect to $f$; see \cite[Definition~1.2.3]{KL06}.

\subsection{Contributions and relation to previous work}
\label{subsection contributions}

The main contribution of this paper is a quenched formulation of the
entropy--energy approach to level--2 large deviations for measurable
skew products with compact fibers. Earlier results for random expanding
maps, random environments, and random perturbations of Axiom~A systems
\cite{Kifer1991RandomExpanding,Kifer1996RandomEnvironments,LiuZhao2003}
provide important precedents. Our abstract results treat a reference
measure through one-sided estimates for its conditional fiber Bowen-ball
masses. The reference measure need not be invariant or arise from
thermodynamic formalism, and the assumptions do not require a spectral
gap for transfer operators or uniqueness of equilibrium states for
perturbed potentials.

The principal geometric step is the lower bound for arbitrary invariant
random measures, which may not be ergodic. A key step is Lemma \ref{Lemma lower estimates}, in which we combine approximation by
finite convex combinations of ergodic measures with random Katok entropy
estimates, and concatenate the resulting orbit segments along one fixed
typical base orbit. This construction bypasses establishing
entropy density of ergodic measures, or its relative version.  This extends the orbit-gluing mechanism of
\cite{Young1990,PfisterSullivan2005} to conditional probabilities in the
fibers.

The energy assumptions also distinguish the roles of the reference
measure and the dynamics. Proposition~\ref{proposition lower energy function}
characterizes the existence of lower and nondegenerate upper energy
functions. Theorems~\ref{thm 2} and~\ref{thm main 1} use the respective
one-sided conditions separately: the upper bound does not require
specification, while the lower bound uses orbit gluing. A common energy
function and upper semicontinuity of fiber entropy then give the good
rate function and the level--2 conclusions of
Theorem~\ref{proposition normalized good rate function}. In the Gibbs
applications, Theorem~\ref{thm simultaneously lower and upper} constructs
a bounded positive common energy function by two coboundary corrections
for potentials with integrable two-sided Bowen variation. Controlling
the resulting endpoint terms on typical fibers allows us to recover the
usual pressure--entropy--potential rate within the abstract framework.

Our results complement the literature on statistical limit laws for
random dynamical systems. Spectral methods yield quenched large deviations
and limit theorems for expanding and hyperbolic systems
\cite{DragicevicFroylandGonzalezTokmanVaienti2018,
DragicevicFroylandGonzalezTokmanVaienti2020}, with extensions to
vector-valued observables \cite{DragicevicHafouta2020} and random
$U(1)$ extensions \cite{ChungNakanoWittsten2023}. Related results include
fiberwise central limit theorems and laws of the iterated logarithm
\cite{Kifer1998LimitTheorems}, and almost sure invariance principles
\cite{DragicevicFroylandGonzalezTokmanVaienti2018ASIP,
DragicevicHafouta2021ASIP}.

Finally, our base transformation is assumed to be invertible and
ergodic, without an additional mixing hypothesis.  A different perspective is provided
by Hafouta \cite{Hafouta2024MixingBase}, who obtained limit theorems,
moderate deviations, and concentration estimates for skew products
with mixing bases and expanding-on-average fibers. Those results
concern orbit sums under the invariant measure of the whole skew
product, so fluctuations of the environment also enter. Our probabilities
are conditioned on a typical environment; the two settings therefore
address distinct, complementary statistical questions.

The rest of this paper is organized as follows. Section~\ref{section example} presents
applications to Gibbs states for skew products that are Anosov and
topologically mixing on fibers. Section~\ref{sec:preliminaries} collects
preliminary results on narrow topologies, energy functions, measure
approximation, and random Katok entropy estimates.
Sections~\ref{sec:upper-bound} and~\ref{sec:lower-bound} establish the
large deviation upper and lower bounds, respectively, with the auxiliary
lemmas for the lower bound proved in Section~\ref{sec proof of lemma}.
Finally, Section~\ref{subsection consequences of main results} proves
that the rate function is good and establishes the level--1 large
deviation and Laplace principles.

\section{Skew products that are Anosov and topologically mixing on fibers}\label{section example}
In this section, we apply the main results to systems that are Anosov and topologically mixing on fibers, which were studied systematically by Huang, Lian, and Lu in \cite{HLLScichina}. We consider the special case in which $M$ is a connected closed smooth Riemannian manifold and $\theta:\Omega\rightarrow \Omega$ is a homeomorphism of a compact metric space $(\Omega,d_\Omega)$. Let $\P$ be any $\theta$-invariant ergodic Borel probability measure on $\Omega$. Let $\operatorname{Diff}^2(M)$ denote the space of $C^2$ diffeomorphisms of $M$ with the $C^2$ topology, and let $F:\Omega\to\operatorname{Diff}^2(M)$ be continuous.
Since the $C^2$ topology is stronger than the compact-open topology, the natural inclusion $\operatorname{Diff}^2(M)\hookrightarrow\operatorname{Homeo}(M,M)$ is continuous. Hence $F$ is also measurable as a map into $\operatorname{Homeo}(M,M)$.
The skew product $\Theta:\Omega\times M\rightarrow \Omega\times M$ induced by $F$ and $\theta$ is defined by
$
  \Theta(\omega,x)=(\theta\omega,F_\omega x),\ \forall\omega\in\Omega,\ x\in M,
$
where we write $F(\omega)$ as $F_\omega$. The skew-product system $\Theta$ is Anosov on fibers if, for every $(\omega,x)\in \Omega\times M$, there is a splitting $
  T_xM=E_{(\omega,x)}^s\oplus E_{(\omega,x)}^u$ that depends continuously on $(\omega,x)\in\Omega\times M$, has $\dim E^s_{(\omega,x)},\ \dim E_{(\omega,x)}^u>0$, and satisfies
\begin{equation*}
  D_xF_\omega E_{(\omega,x)}^u=E_{\Theta(\omega,x)}^u,\ \ D_xF_\omega E_{(\omega,x)}^s=E_{\Theta(\omega,x)}^s
\end{equation*}
and
\begin{equation*}
  \begin{cases}
    |D_xF_\omega\xi|\geq e^{\lambda_0}|\xi|, & \ \forall \xi\in E_{(\omega,x)}^u, \\
    |D_xF_\omega\eta|\leq e^{-\lambda_0}|\eta|, &\ \forall \eta\in E_{(\omega,x)}^s,
  \end{cases}
\end{equation*}
where $\lambda_0>0$ is a constant; and
$\Theta$ is topologically mixing on fibers if, for any nonempty open sets $U,V\subset M$, there exists $N>0$ such that, for every $n\geq N$ and $\omega\in\Omega$, $
  \Theta^n(\{\omega\}\times U)\cap (\{\theta^n\omega\}\times V)\not=\emptyset.
$

For $f\in L^1(\Omega,C(M))$, we write $S_{m,n}f(\omega,x)=\sum_{j=m}^{n-1}f(\Theta^j(\omega,x))$. When $m=0$, we simply denote $S_nf(\omega,x)=\sum_{j=0}^{n-1}f(\Theta^j(\omega,x))$. Define the two-sided Bowen ball by
\begin{equation*}
  B_{m,n}^\omega(x,r):=\{y\in M: d_M(F_\omega^j x,F_\omega^j y)< r\mbox{ for every }m\leq j<n\}.
\end{equation*}
For $f\in L^1(\Omega,C(M))$, we define
\begin{equation*}
  K_f^+(\omega):=\sup_{n\geq 1, x,y\in M,\ y\in B_{0,n}^\omega(x,\epsilon_0)}|S_{0,n}f(\omega,x)-S_{0,n}f(\omega,y)|
  \end{equation*}
and
\begin{equation*}
  K_f^-(\omega):=\sup_{m<0, x,y\in M,\ y\in B_{m,0}^\omega(x,\epsilon_0)}|S_{m,0}f(\omega,x)-S_{m,0}f(\omega,y)|,
\end{equation*}
where the constant $\epsilon_0>0$ is the fiber expansive constant for the Anosov-on-fibers system $F$ (see \cite{HLLScichina}), i.e.,
\begin{equation*}
  \mbox{ if $d_M(F_\omega^n x,F_\omega^n y)\leq \epsilon_0$ for $n\in\mathbb{Z}$, then $x=y$.}
\end{equation*}
We choose $\epsilon_0$ smaller, if necessary, so that the uniform local
stable/unstable manifold estimates and local product structure apply at this scale.
 The class of potential functions considered by Gundlach and Kifer in \cite{Kifer2000specification} is
\begin{equation*}
  V_{\P}^\pm(F):=\{f\in L^1(\Omega,C(M)):\ \int_\Omega\bigl(K_f^+(\omega)+K_f^-(\omega)\bigr)\,d\P(\omega)<\infty\}
\end{equation*}
A convenient uniform subclass is
\begin{equation*}
  B^\pm(F):=\{f\in C(\Omega\times M):\ \sup_{\omega\in\Omega} K_f^\pm(\omega)<\infty\}\subset V_\P^\pm(F).
\end{equation*}

\begin{lemma}[A sufficient regularity condition]\label{lemma sufficent regularity}
Assume that $(\Omega\times M,\Theta)$ is Anosov on fibers. Suppose that $f\in C(\Omega\times M)$ is equi-H\"older continuous in the fiber variable; that is, there exist constants $H_f>0$ and $\alpha\in(0,1]$ such that
  \begin{equation*}
    |f(\omega,x)-f(\omega,y)|\leq H_fd_M(x,y)^\alpha,\ \forall \omega\in\Omega,\ \forall x,y\in M.
  \end{equation*}
Then $f\in B^\pm(F)$.
\end{lemma}
\begin{proof}[Proof of Lemma \ref{lemma sufficent regularity}]
  Uniform hyperbolicity and the local product structure on the fibers (see \cite{HLLScichina}) give, for each fixed $\lambda\in(0,\lambda_0)$, a constant $L_F\geq1$ independent of $\omega,x,y,n$ such that, if $y\in B_{0,n}^\omega(x,\epsilon_0)$,
  \begin{equation*}
    d_M(F_\omega^i x,F_\omega^i y)\leq L_F\cdot \epsilon_0\cdot(e^{-\lambda i}+e^{-\lambda(n-i)}),\ 0\leq i<n.
  \end{equation*}
Indeed, use the local product point of $x$ and $y$. The stable component
is contracted from the initial endpoint, and the unstable component is
estimated backwards from the final endpoint. The uniform local product
bounds control both endpoint distances by a constant times $\epsilon_0$;
changing $n-1-i$ to $n-i$ only changes that constant.
By direct summation and the inequality $(a+b)^\alpha\leq a^\alpha+b^\alpha$ for $\alpha\in(0,1]$, we have
  \begin{equation*}
    K_f^+(\omega)\leq \frac{2 H_f\cdot (L_F\cdot \epsilon_0)^\alpha}{1-e^{-\alpha\lambda}},\ \forall \omega\in\Omega,
  \end{equation*}
for any $f\in C(\Omega\times M)$ equi-H\"older continuous in the fiber variable. The same finite-segment estimate, applied to the segment starting at $\theta^m\omega$ for $m<0$, gives the identical bound for $K_f^-(\omega)$.
\end{proof}
Two points $x,y\in M$ are called $(\omega,n,\epsilon)$-separated if $d_\omega^n(x,y)\geq \epsilon$. For any $f\in L^1(\Omega,C(M))$, $n\geq 1$, and constant $\epsilon>0$, set
\begin{equation}\label{eq pi F varphi}
  \pi_F^*(f)(\omega,\epsilon,n)=\sup\{Z_n(\omega,f,Q):\ Q\text{ is a maximal $(\omega,n,\epsilon)$-separated subset of }M\},
\end{equation}
where
\begin{equation*}
  Z_n(\omega,f,Q)=\sum_{x\in Q}\exp(S_{0,n}f(\omega,x))=\sum_{x\in Q}\exp\left(\sum_{i=0}^{n-1}f(\Theta^i(\omega,x))\right).
\end{equation*}
Here maximal means maximal by inclusion. Every separated set can be extended to such a set, and the positive weights ensure that its weighted sum cannot decrease. The fiber topological pressure of $F$ (or the relative topological pressure of $\Theta$) with respect to $f\in L^1(\Omega,C(M))$ is
\begin{equation}\label{eq fiber top pressure}
\begin{split}
  \pi_F(f)    =& \lim_{\epsilon\to0}\limsup_{n\to \infty}\frac{1}{n}\int\log\pi_F^*(f)(\omega,\epsilon,n)d\mathbb{P}(\omega) \\
    = & \lim_{\epsilon\to 0}\limsup_{n\to \infty}\frac{1}{n}\log\pi_F^*(f)(\omega,\epsilon,n)
    =\lim_{\epsilon\to 0}\liminf_{n\to \infty}\frac{1}{n}\log\pi_F^*(f)(\omega,\epsilon,n)
\end{split}
\end{equation}
for $\P$-a.s. $\omega\in\Omega$, where the last two equalities follow from
\cite[Proposition 1.2.6]{KL06}, since $\mathbb{P}$ is ergodic with respect to $\theta$. The quantity $h_{\mathrm{top}}(F)=\pi_F(0)$ is called the fiber topological entropy of $\Theta$.
The pressure is finite in the present setting. Indeed, compactness of the base
and continuity into $\operatorname{Diff}^2(M)$ give
\[
 d=\dim M,\qquad
 L=\max\left\{1,\sup_{\omega,x}\|D_xF_\omega\|\right\}<\infty.
\]
An ordinary ball of radius $\epsilon L^{-(n-1)}/2$ has $n$-step Bowen
diameter less than $\epsilon$. The covering estimate on a compact
$d$-dimensional manifold therefore gives $h_{\mathrm{top}}(F)\leq d\log L$.
Writing $b_f(\omega)=\|f(\omega,\cdot)\|_\infty$, we obtain
\[
 h_{\mathrm{top}}(F)-\int b_f\,d\P
 \leq\pi_F(f)\leq h_{\mathrm{top}}(F)+\int b_f\,d\P<\infty.
\]
The first assertion of the following lemma is
\cite[Theorem~A]{Liu2025variatonalprinciple}; that result does not require
unique ergodicity of the base.
\begin{lemma}
  Assume that $(\Omega\times M,\Theta)$ is a system that is Anosov and topologically mixing on fibers and is driven by a homeomorphism $(\Omega,\theta)$. Then $\Theta$ has the fiber specification property. Moreover, the RDS $F$ corresponding to $\Theta$ has the Gundlach--Kifer specification property defined in \cite[Def. 2.5]{Kifer2000specification}.
\end{lemma}
\begin{remark}\label{remark GK specification}
  The Gundlach--Kifer specification property is defined as follows; see \cite[Definition~2.5]{Kifer2000specification}.
  Let $\epsilon_0:\Omega\to\mathbb{R}^+$ be the expansive characteristic for RDS $F$, i.e., for any $\omega\in\Omega$
   \begin{equation*}
     \mbox{if $d_M(F_\omega^n(x),F_\omega^n(y))\leq \epsilon_0(\theta^n\omega)$ for $n\in\mathbb{Z}$, then $x=y$.}
   \end{equation*}
 The RDS $F$ has $k$-specification if for each constant $c>0$ and $\mathbb P$-a.s. $\omega\in\Omega$, there exists an $\mathbb{N}$-valued random variable $L_c=L_c(\omega)\geq 1$ such that for any points $x_0,\ldots,x_k\in M$ and integers $-\infty\leq b_0<a_1<b_1<\cdots<a_k<b_k<a_{k+1}\leq \infty$ satisfying $a_{i+1}\geq b_i+L_c(\theta^{b_i}\omega),$ $i=0,\ldots,k$ whenever $b_i$ is finite, one can find $z\in M$ such that
  \begin{align}
    d_M(F_\omega^j(x_i),F_\omega^j(z))&\leq c\epsilon_0(\theta^j\omega),\quad a_i\leq j\leq b_i,\ 1\leq i\leq k\label{kifer s1}\\
    d_M(F_\omega^j(x_0),F_\omega^j(z))&\leq c\epsilon_0(\theta^j\omega)\ \ \forall j\geq a_{k+1} \mbox{ and }\forall j\leq b_0.\label{kifer s2}
  \end{align}
If $F$ satisfies $k$-specification for any $k\in \mathbb{N}_+$ with $L_c$ independent of $k$, then $F$ is said to have Gundlach--Kifer specification property.

For systems that are Anosov and topologically mixing on fibers, the
expansivity characteristic is constant and the specification gap is independent
of $k$ and $\omega$.

The compact fiber and the uniform hyperbolic estimates give uniform bounds
for the fixed-scale covering numbers and strong expansivity times used in
\cite{Kifer2000specification}. At each fixed scale the specification gap is
deterministic. The remaining potential-dependent quantities are bounded by
finite sums of translates of $K_f^+$, $K_f^-$ and
$\|f(\omega,\cdot)\|_\infty$, which are integrable. Thus the covering and
integrability conditions in \cite[Lemma~3.8 and Theorem~3.9]{Kifer2000specification}
are satisfied in this uniform setting.
% AUTHOR-CHECK R03: verify the exact constants in the original GK (3.15)
% and the same-scale open-ball comparison; see the revision report.
\end{remark}

The following is extracted from the Gibbs-state results of Gundlach and Kifer; see Propositions 3.2 and 3.6, Lemma 3.8 and Theorem 3.9 in \cite{Kifer2000specification}.
\begin{proposition}\label{proposition 2.1}
 Assume that $(\Omega\times M,\Theta)$ is a system that is Anosov and topologically mixing on fibers and is driven by a homeomorphism $(\Omega,\theta)$.
  For every $f\in V_\P^\pm(F)$, there is a unique equilibrium state $\nu_f\in \operatorname{Pr}_\P(\Omega\times M)$ such that
  \begin{equation*}
    \pi_F(f)=h_{\nu_f}(F)+\int_{\Omega\times M} f d\nu_f.
  \end{equation*}
Moreover, $\nu_f$ is $\Theta$-invariant, ergodic, and has the Gibbs property, i.e., for every $0<c\leq 1/4$, its fiber Bowen ball masses satisfy
  \begin{equation}\label{eq error Bowenball}
    \log (\nu_f)_\omega(B_\omega(x,n,c\epsilon_0))=S_nf(\omega,x)-\log \pi_F^*(f)(\omega,c\epsilon_0,n)+E_{n,c}(\omega,x),
  \end{equation}
where the error satisfies
  \begin{equation}\label{eq estimate error}
    \lim_{n\to\infty}\frac{\sup_{x\in M}|E_{n,c}(\omega,x)|}{n}=0,\ \mbox{ for $\P$-a.s. $\omega\in\Omega$.}
  \end{equation}
  Furthermore, if $f\in B^\pm(F)$, then for every $0<c\leq 1/4$, there is a deterministic constant $C_c\geq 1$ such that for $\P$-a.s. $\omega\in\Omega$,
  \begin{equation*}
    (C_c)^{-1}\frac{e^{S_nf(\omega,x)}}{\pi_F^*(f)(\omega,c\epsilon_0,n)}\leq (\nu_f)_\omega(B_\omega(x,n,c\epsilon_0))\leq  C_c\frac{e^{S_nf(\omega,x)}}{\pi_F^*(f)(\omega,c\epsilon_0,n)}.
  \end{equation*}
\end{proposition}
\begin{proof}[Proof of Proposition \ref{proposition 2.1}]
The proof is exactly a computation following \cite{Kifer2000specification}.
We apply Gundlach--Kifer with $\epsilon_0/2$ as the expansivity
characteristic. Since closed Bowen balls at this scale are contained
in the corresponding open Bowen balls of radius $\epsilon_0$,
their forward and backward variations are bounded by
$K_f^+$ and $K_f^-$, respectively. By remark \ref{remark GK specification}, the specification gaps and the
auxiliary expansivity times can be chosen uniformly in $\omega$.
Consequently, the covering and integrability hypotheses in
\cite[Lemma~3.8 and Theorem~3.9]{Kifer2000specification} hold: the
potential-dependent terms in their condition~(3.15) are bounded by
finitely many translates of $K_f^++K_f^-$ and
$\|f(\omega,\cdot)\|_\infty$.
Proposition~3.6, Lemma~3.8 and Theorem~3.9 of that paper therefore give
an invariant ergodic Gibbs state $\nu_f$ and identify it as the unique
equilibrium state. Their limiting pressure agrees with $\pi_F(f)$,
since strict and non-strict separation give partition sums lying
between one another at radii $r$ and $r/2$.

Fix $0<c\leq1/4$ and put $
 r=c\epsilon_0$ and $s=\frac{c\epsilon_0}{2}.
$ With $\epsilon_0/2$ as the expansivity characteristic,
\cite[Proposition~3.2, (3.3)]{Kifer2000specification}
gives estimates for closed Bowen balls of radius $s$.
We now compare these estimates with the open Bowen balls of radius
$r$ appearing in our statement.
 To pass to our open balls and separation convention, observe
that uniform expansivity bounds the cardinality of any
$s$-separated subset of $B_\omega(x,n,r)$ by a constant independent
of $\omega,x,n$: separation can occur only within a fixed number of
iterates of the two endpoints, where compactness gives a uniform
packing bound. Covering by the corresponding closed $s$-Bowen balls
and using $K_f^+(\omega)$ to compare their weights therefore changes
the Gibbs bounds by at most a deterministic factor times
$e^{K_f^+(\omega)}$. The endpoint factors in (3.3), and the finitely
many short lengths excluded from its upper bound, can be absorbed
into a nonnegative function $A_c\in L^1(\P)$, giving
\[
 \left|\log
 \frac{(\nu_f)_\omega(B_\omega(x,n,r))\,
       \pi_F^*(f)(\omega,r,n)}{e^{S_nf(\omega,x)}}\right|
 \leq A_c(\omega)+A_c(\theta^n\omega)
\]
for all $n\geq1$ and $x\in M$, on a set of full measure.
Here $A_c$ is bounded by a constant plus finitely many translates of
$K_f^++K_f^-$ and $\|f(\omega,\cdot)\|_\infty$.

For every $\eta>0$, invariance of $\P$ and integrability give
\[
 \sum_{n\geq1}\P\{A_c\circ\theta^n>n\eta\}
 =\sum_{n\geq1}\P\{A_c>n\eta\}\leq \frac{1}{\eta}\int_\Omega A_cd\P<\infty.
\]
Thus $A_c(\theta^n\omega)/n\to0$ almost surely by Borel--Cantelli lemma,
which proves \eqref{eq estimate error} for the error in
\eqref{eq error Bowenball}, uniformly in $x$.
If $f\in B^\pm(F)$, the same construction makes $A_c$ bounded, so
$C_c=\exp(2\|A_c\|_\infty)$ gives the asserted uniform Gibbs bounds.
\end{proof}

\begin{theorem}\label{thm gibbsstate energy function}
  Assume that $(\Omega\times M,\Theta)$ is a system that is Anosov and topologically mixing on fibers and is driven by a homeomorphism $(\Omega,\theta)$. For every $f\in V_\P^\pm(F)$, let $\nu=\nu_f$ be the corresponding random Gibbs state. Then $\nu$ admits a lower energy function and a nondegenerate upper energy function.
\end{theorem}
\begin{proof}[Proof of Theorem \ref{thm gibbsstate energy function}]
Step 1. We prove that $\nu$ admits a lower energy function.
 For $f\in V_\P^\pm(F)\subset L^1(\Omega,C(M)),$ define
  \begin{equation*}
    \mbox{$a_n(\omega):=\min\limits_{x\in M}S_nf(\omega,x)$ and $b_n(\omega):=\max\limits_{x\in M}S_nf(\omega,x)$.}
  \end{equation*} Then the cocycle identity implies that
$
   a_{n+m}(\omega)\geq a_n(\omega)+a_m(\theta^n\omega)\mbox{ and }b_{n+m}(\omega)\leq b_n(\omega)+b_m(\theta^n\omega).
$ Moreover, $a_n,b_n\in L^1(\P)$ because
 \begin{equation}\label{eq bound an bn}
   |a_n(\omega)|+|b_n(\omega)|\leq 2\sum_{j=0}^{n-1}\|f(\theta^j\omega,\cdot)\|_{\infty}\in L^1(\P)
 \end{equation}
since $f\in V_\P^\pm(F)\subset L^1(\Omega,C(M))$. Applying Kingman's subadditive ergodic theorem to $-a_n$ and $b_n$ and using the ergodicity of $\P$, we obtain the following two limits:
 \begin{equation}\label{eq def beta}
   \beta_-(f):=\lim_{n\to\infty}\frac{a_n(\omega)}{n},\ \beta_+(f):=\lim_{n\to\infty}\frac{b_n(\omega)}{n},\ \mbox{ for $\P$-a.s.} \omega\in\Omega.
 \end{equation}
Moreover, \eqref{eq bound an bn} and Birkhoff's ergodic theorem imply that $|\beta_\pm(f)|\leq \int_\Omega\|f(\omega,\cdot)\|_\infty d\P(\omega)<+\infty.$

 For fixed $c\in(0,1/4]$, set $\epsilon=c\epsilon_0$ and take the infimum or supremum over $x\in M$ in \eqref{eq error Bowenball}. Combining \eqref{eq estimate error}, \eqref{eq fiber top pressure}, and \eqref{eq def beta}, and then letting $n\to\infty$ and $c\searrow 0$, yields, for $\P$-a.s. $\omega\in\Omega$,
 \begin{equation}\label{eq admits lower energy}
   \lim_{\epsilon\searrow0}\liminf_{n\to+\infty}\inf_{x\in M}\frac{1}{n}\log \nu_\omega (B_\omega(x,n,\epsilon))=\beta_-(f)-\pi_F(f)>-\infty,
 \end{equation}
and
 \begin{equation}\label{eq admit uuper energy}
   \lim_{\epsilon\searrow0}\limsup_{n\to+\infty}\sup_{x\in M}\frac{1}{n}\log \nu_\omega (B_\omega(x,n,\epsilon))=\beta_+(f)-\pi_F(f).
 \end{equation}
Therefore, \eqref{eq proposition exists energy} holds, and $\nu$ admits a lower energy function.

 Step 2. To prove that $\nu$ admits a nondegenerate upper energy function, it suffices to show that $\beta_+(f)-\pi_F(f)<0$. For $\omega\in\Omega$ in a $\P$-full-measure set on which \eqref{eq def beta} and \eqref{eq fiber top pressure} hold, let $y_n=y_n(\omega)\in M$ satisfy $b_n(\omega)=S_nf(\omega,y_n)$. Extend the separated singleton $\{y_n\}$ to a maximal separated set. Since all weights are positive, \eqref{eq pi F varphi} gives
 \begin{equation*}
 	\pi_F^*(f)(\omega,\epsilon,n)\geq e^{b_n(\omega)},
 \end{equation*}
and therefore $\pi_F(f)\geq \beta_+(f)$. To prove strict inequality, we construct exponentially many separated orbits whose potential sums remain sufficiently large.

Choose $\delta>0$ sufficiently small that $\delta<\min\{\operatorname{diam}(M)/8,\epsilon_0\}$. Let $L=L(\delta)$ be the gap length in the fiber specification property corresponding to $\delta$. Define
\begin{equation}\label{eq def K delta}
	E_\delta(\omega):=\sup_{n\geq 1, x,y\in M, y\in B_{0,n}^\omega(x,\delta)}|S_{0,n}f(\omega,x)-S_{0,n}f(\omega,y)|\leq K_f^+(\omega).
\end{equation}  Define
\begin{equation*}
	\Gamma(\omega):=E_\delta(\theta^{L+1}\omega)+2\sum_{i=-L}^{i=L}\|f(\theta^i\omega,\cdot)\|_\infty\in L^1(\P).
\end{equation*}
Choose $D>0$ such that $\P\{\Gamma\leq D\}>3/4$. By Birkhoff's ergodic theorem, for $\P$-a.s. $\omega\in\Omega$, we have
\begin{equation}\label{eq positive recurrent}
	\lim_{n\to \infty}\frac{1}{n}\sum_{i=0}^{n-1}1_{\{\Gamma\leq D\}}(\theta^i\omega)>3/4.
\end{equation}
From now on, fix $\omega\in\Omega$ in a $\P$-full-measure set on which \eqref{eq fiber top pressure}, \eqref{eq def beta}, and \eqref{eq positive recurrent} hold. For sufficiently large $n\in\mathbb{N}$, choose $I_n(\omega)\subset \{L+1,...,n-L-2\}$ such that
\begin{enumerate}
	\item $|i-j|\geq 2L+3, \forall i\not=j\in I_n(\omega)$;
	\item  for every $i\in I_n(\omega)$, $\Gamma(\theta^i\omega)\leq D$,
	\item and
	\begin{equation}\label{eq lower density}
		c_0:=\liminf\limits_{n\to\infty}\frac{\#I_n(\omega)}{n}\geq \frac{\P\{\Gamma\leq D\}}{2L+3}>0.
	\end{equation}
\end{enumerate}
This can be done inductively as follows. Let $i_1\geq L+1$ be the first integer such that $\Gamma(\theta^{i_1}\omega)\leq D$. Once $i_k$ is defined, let $i_{k+1}$ be the first integer greater than or equal to $i_{k}+2L+3$ and less than or equal to $n-L-2$ such that $\Gamma(\theta^{i_{k+1}}\omega)\leq D$. Estimate \eqref{eq lower density} is a consequence of \eqref{eq positive recurrent}.

For each $i\in I_n(\omega)$, choose $e_i\in M$ such that $d_M(e_i,F_\omega^i y_n)\geq \frac{\operatorname{diam}(M)}{2}$. Recall that $y_n=y_n(\omega)$ satisfies $b_n(\omega)=S_nf(\omega,y_n)$. For each subset $J\subset I_n(\omega)$, define a specification as follows. Begin with the orbit of $y_n$ on $[0,n-1]$. For every $i\in J$, remove the window $[i-L,i+L]$ and insert the singleton orbit segment $e_i$ at time $i$; the $L$ iterates on either side serve as gaps in the fiber specification property. The windows are disjoint, and successive prescribed intervals are more than $L$ apart, so the fiber specification property produces a point $x_J$.

The shadowing points $x_J$ corresponding to distinct choices of $J$ are $(\omega,n,\frac{\operatorname{diam}(M)}{2}-2\delta)$-separated. Let $r_0=\frac{\operatorname{diam}(M)}{2}-2\delta.$ Thus, we have constructed $2^{\#I_n(\omega)}$ distinct $(\omega,n,r_0)$-separated orbits.

We now estimate $S_nf(\omega,x_J)$. If $J=\varnothing$, take $x_J=y_n$.
If $J=\{i_1<\cdots<i_q\}$ is nonempty, the retained intervals are
\[
 [0,i_1-L-1],\quad
 [i_r+L+1,i_{r+1}-L-1]\ (1\leq r<q),\quad
 [i_q+L+1,n-1].
\]
The variation on the initial retained interval is at most $E_\delta(\omega)$;
on each subsequent retained interval it is at most
$E_\delta(\theta^{i_r+L+1}\omega)$. On a removed window we use
\[
 \left|f(\theta^j\omega,F_\omega^j y_n)
       -f(\theta^j\omega,F_\omega^j x_J)\right|
 \leq2\|f(\theta^j\omega,\cdot)\|_\infty,
 \qquad i-L\leq j\leq i+L.
\]
Note that for $i\in I_n(\omega)$, $\Gamma(\theta^i\omega)\leq D.$ Therefore, we have
\begin{equation*}
\begin{split}
 S_nf(\omega,x_J)\geq   & b_n(\omega)-E_\delta(\omega)-\sum_{i\in J}(E_\delta(\theta^{i+L+1}\omega)+2\sum_{p=-L}^{L}\|f(\theta^{i+p}\omega,\cdot)\|_\infty) \\
    \geq  & 	b_n(\omega)-E_\delta(\omega)-D\cdot \#J.
\end{split}
\end{equation*}

 Extending the constructed separated set to a maximal one, we obtain
 \begin{equation*}
 	\pi_F^*(f)(\omega,r_0,n)\geq \sum_{J\subset I_n(\omega)} e^{S_nf(\omega,x_J)}\geq e^{b_n(\omega)-E_\delta(\omega)}\sum_{J\subset I_n(\omega)}e^{-D\#J}.
 \end{equation*}
Note that
 \begin{equation*}
 	\sum_{J\subset I_n(\omega)}e^{-D\#J}=\sum_{k=0}^{\#I_n(\omega)}\begin{pmatrix}
 		\#I_n(\omega)\\ k
 	\end{pmatrix}e^{-Dk}=(1+e^{-D})^{\#I_n(\omega)}.
 \end{equation*}
Therefore, we have
 \begin{equation*}
 \pi_F(f)\geq 	\liminf_{n\to \infty}\frac{1}{n}\log \pi_F^*(f)(\omega,r_0,n)\geq \beta_+(f)+c_0\log(1+e^{-D})>\beta_+(f).
 \end{equation*}
 The two equalities \eqref{eq admit uuper energy} and \eqref{eq admits lower energy} imply that, for the Gibbs state $\nu_f$ corresponding to $f\in V_\P^\pm(F)$, the constant function $\pi_F(f)-\beta_-(f)$ is a lower energy function for $\nu_f$, while the constant function $\pi_F(f)-\beta_+(f)$ is an upper energy function for $\nu_f$.
This completes the proof.
\end{proof}
\subsection{A common energy function}
The Gibbs estimates \eqref{eq error Bowenball} and the definition of energy function motivate the following construction of a common
lower and upper energy function. We seek an energy
function $T_f$ such that
\begin{equation*}
  \lim_{n\to\infty}\frac1n\sup_{x\in M}
  \left|S_nT_f(\omega,x)+S_nf(\omega,x)-n\pi_F(f)\right|=0
\end{equation*}
for $\P$-a.s. $\omega\in\Omega$, so that its orbit sums cancel the
exponential contribution of the Gibbs potential uniformly on each
fiber. The natural candidate is $\pi_F(f)-f$. However, it need not be nonnegative
or essentially bounded for a general $f\in V_\P^\pm(F)$.
We therefore seek a coboundary correction
\[
  T_f=\pi_F(f)-f+H-H\circ\Theta.
\]
Such a correction changes orbit sums only by the endpoint term
$H-H\circ\Theta^n$; provided this term is uniformly $o(n)$ on each
typical fiber, it preserves the required exponential estimates.
The construction proceeds in two steps.
 First, choosing
$\beta_+(f)<a<\pi_F(f)$, we use the maximal accumulated excess
$u=\sup_{n\geq0}(S_nf-na)$, with $S_0f=0$, to obtain
$g=\pi_F(f)-f+u-u\circ\Theta\geq \pi_F(f)-a>0$.
The strict inequality $a>\beta_+(f)$ ensures that this supremum is
finite on typical fibers. Second, choosing
$b>\pi_F(f)-\beta_-(f)$, we apply the same construction to $g$:
\[
  v=\sup_{n\geq0}(S_ng-nb),
  \qquad
  T_f=g+v\circ\Theta-v.
\]
Since $v\geq0$, this second correction imposes the upper bound
without destroying the positive lower bound, giving
$\pi_F(f)-a\leq T_f\leq b$. In the following Theorem \ref{thm simultaneously lower and upper}, we also prove that
$\int T_f\,d\rho=\pi_F(f)-\int f\,d\rho$ for every
$\rho\in I_\P(\Omega\times M)$.
Consequently, the entropy--energy expression in \eqref{eq rate function I T} becomes
\[
  \int T_f\,d\rho-h_\rho(F)
  =
  \pi_F(f)-h_\rho(F)-\int f\,d\rho,
\]
which is the standard Gibbs large-deviation rate on invariant
measures. Thus the auxiliary parameters change the energy
representative while leaving the resulting rate function unchanged.
\begin{theorem}\label{thm simultaneously lower and upper}
Assume the hypotheses of Theorem~\ref{thm gibbsstate energy function}. Let $f\in V_\P^\pm(F)\subset L^1(\Omega,C(M))$, and let $\nu_f$ be its Gibbs equilibrium state. Recall that $\beta_-(f)\leq \beta_+(f)<\pi_F(f)$ are defined in \eqref{eq def beta}. Choose any constants $a,b$ such that
\begin{equation*}
  \beta_+(f)<a<\pi_F(f),\mbox{ and }b>\pi_F(f)-\beta_-(f).
\end{equation*}
On a measurable $\theta$-invariant full-measure set $\Omega_*\subset\Omega$ where the limits in \eqref{eq def beta}, the two applications of Lemma~\ref{lemma upper bounded representative} below, and the required endpoint estimates all hold, define
\begin{equation}\label{eq def u,v}
\begin{split}
   u(\omega,x)=\sup_{n\geq 0}\{S_nf(\omega,x)-na\}, &\ \ g(\omega,x)=\pi_F(f)-(f(\omega,x)+u(\Theta(\omega,x))-u(\omega,x)) \\
    v(\omega,x)=\sup_{n\geq 0}\{S_ng(\omega,x)-nb\}, & \ \ T_f(\omega,x)=g(\omega,x)+v(\Theta(\omega,x))-v(\omega,x).
\end{split}
\end{equation}
with the convention that $S_0f=S_0g=0$.
The set $\Omega_*$ is constructed successively in the proof; after each step we take its $\theta$-invariant full-measure core. Set $T_f(\omega,x)=\pi_F(f)-a$ on the exceptional fibers.  Then the function $T_f$ is measurable, continuous on each fiber and satisfies
  \begin{equation}\label{eq bound of Tf}
    0<\pi_F(f)-a\leq T_f\leq b<+\infty.
  \end{equation}
Moreover, $T_f$ is simultaneously a lower and an upper energy function for $\nu_f$ satisfying
  \begin{equation}\label{eq limit Sntf}
    \lim_{n\to+\infty}\frac{1}{n}\sup_{x\in M}|S_nT_f(\omega,x)+S_nf(\omega,x)-n\pi_F(f)|=0\mbox{ for $\P$-a.s. $\omega\in\Omega$,}
  \end{equation}
and
  \begin{equation}\label{eq int Tf}
    \int_{\Omega\times M} T_fd\rho=\pi_F(f)-\int_{\Omega\times M} fd\rho,\mbox{ for every }\rho\in I_\P(\Omega\times M).
  \end{equation}
\end{theorem}
To prove Theorem \ref{thm simultaneously lower and upper}, we need the following two lemmas.
\begin{lemma}\label{lemma converges to 0}
  Let $R:\Omega\to [0,+\infty)$ be measurable and finite on a $\theta$-invariant $\P$-full-measure set. If $R\circ \theta-R\in L^1(\P)$, then $R(\theta^n\omega)/n\to 0$ for $\P$-a.s. $\omega\in\Omega$ as $n\to+\infty$.
\end{lemma}
\begin{proof}[Proof of Lemma \ref{lemma converges to 0}]
  Note that we do not assume $R\in L^1(\P)$. Set $d=R\circ \theta- R\in L^1(\P)$. By Birkhoff ergodic theorem, the following
  \begin{equation*}
    \lim_{n\to+\infty}\frac{1}{n}\sum_{j=0}^{n-1}d(\theta^j\omega)=\lim_{n\to+\infty}\frac{R(\theta^n\omega)-R(\omega)}{n}=\lim_{n\to+\infty}\frac{R(\theta^n\omega)}{n}
  \end{equation*}
has a $\P$-a.s. limit. Moreover, for any $\epsilon>0$, the $\theta$-invariance gives
  \begin{equation*}
    \lim_{n\to+\infty}\P\{\omega\in\Omega:\ R(\theta^n\omega)>n\epsilon\}=\lim_{n\to+\infty}\P\{\omega\in\Omega:\ R(\omega)>n\epsilon\}=0
  \end{equation*}
Thus the sequence $\frac{R(\theta^n\omega)}{n}$ converges to $0$ in probability. Therefore, its almost-sure limit is $0$.
\end{proof}

\begin{lemma}\label{lemma upper bounded representative}
  Let $k\in L^1(\Omega,C(M))$. Suppose
    \begin{equation}\label{eq limit exists gamma}
      \lim\limits_{n\to+\infty}\frac{1}{n}\max\limits_{x\in M}S_nk(\omega,x)=\gamma \mbox{ for some
$\gamma\in\mathbb{R}$ and  $\P$-a.s. $\omega\in\Omega$.}
    \end{equation}
Then for $c>\gamma$, $h(\omega,x):=\sup_{n\geq0}\{S_nk(\omega,x)-nc\}$
is a finite, nonnegative, measurable function, continuous on $\P$-a.s. fiber.
 The following potential $
 \hat k=k+h\circ\Theta-h$
satisfies
\begin{equation}\label{eq:k-bounds}
\hat{k}=\min\{k+h\circ \Theta, c\},\qquad \min\{k,c\}\leq\hat k\leq c,
 \qquad |h\circ\Theta-h|\leq|k-c|,
\end{equation}
on $\P$-a.s. fibers. In particular $\hat k\in L^1(\Omega,C(M))$. If
$H(\omega)=\max_{x\in M} h(\omega,x)$, then
\begin{equation}\label{eq:w-sublinear}
 \frac{H(\theta^n\omega)}n\to 0,\mbox{ and }
 \frac{1}{n}\sup_{x\in M}|S_n\widehat k (\omega,x)-S_nk(\omega,x)|\to 0,\mbox{ for $\P$-a.s. $\omega\in\Omega$}
\end{equation}
\end{lemma}
\begin{proof}[Proof of Lemma \ref{lemma upper bounded representative}]
Fix $\omega$ from  a measurable $\theta$-invariant full-measure set on which \eqref{eq limit exists gamma} holds and all the fiber
norms of $k$ along the base orbit are finite. Since $\gamma<c$, for all sufficiently large $n\geq 0$, $S_nk(\omega,x)<nc$ for every $x$. Then $h(\omega,x)=\sup_{n\geq0}\{S_nk(\omega,x)-nc\}$ is therefore a maximum over finitely many functions, including the zero function when $n=0$. This proves the finiteness and continuity. Measurability follows from the countable supremum. Shifting the summation index gives
\begin{equation}\label{eq expression of h}
  h(\omega,x)=\sup_{n\geq 0}\{S_nk(\omega,x)-nc\}=\max\{0,k(\omega,x)-c+h\circ \Theta(\omega,x)\}\geq 0.
\end{equation}
Then
\begin{equation*}
  \hat{k}=k+h\circ \Theta-h=k+h\circ \Theta-\max\{0,k-c+h\circ \Theta\}=\min\{k+h\circ \Theta,c\}\geq \min\{k,c\}
\end{equation*}
and
\begin{equation*}
  |h\circ\Theta-h|=|\hat k-k|=|\min\{k+h\circ\Theta,c\}-k|\leq |k-c|.
\end{equation*}
Since $\hat{k}=\min\{k+h\circ \Theta,c\}\leq c$ on $\P$-a.s. fibers, its fiber norm is integrable. Therefore, $\hat{k}\in L^1(\Omega,C(M)).$

Since $F_\omega$ is onto for all $\omega\in\Omega$, we have
\begin{equation*}
  H(\theta\omega)=\max_{x\in M}h(\theta\omega,x)=\max_{x\in M} h(\theta\omega,F_\omega x)=\max_{x\in M}h(\Theta(\omega,x)).
\end{equation*}
Consequently, we have
\begin{equation*}
  |H\circ \theta-H|\leq \max_{x\in M}|(h(\Theta(\omega,x))-h(\omega,x))|\leq \|k(\omega,\cdot)-c\|_\infty\in L^1(\P).
\end{equation*}
Lemma \ref{lemma converges to 0} implies that $\frac{H(\theta^n\omega)}{n}\to 0$ for $\P$-a.s. $\omega\in\Omega$ as $n\to +\infty$. Telescoping gives for $\P$-a.s. $\omega\in\Omega$,
\begin{equation*}
  \frac{1}{n}\sup_{x\in M}|S_n\hat{k}(\omega,x)-S_nk(\omega,x)|\leq\frac{1}{n}( H(\omega)+H(\theta^n\omega))\to 0,
\end{equation*}
which completes the proof of Lemma \ref{lemma upper bounded representative}.
\end{proof}
Now we are ready to prove Theorem \ref{thm simultaneously lower and upper}.
\begin{proof}[Proof of Theorem \ref{thm simultaneously lower and upper}]
We construct the common full-measure set successively. Start with the
invariant core on which the two Kingman limits in \eqref{eq def beta}
hold and all fiber norms of $f$ along the base orbit are finite.
Each application of Lemma~\ref{lemma upper bounded representative}
uses its measurable extensions on exceptional fibers. After imposing
its endpoint limit, pass again to the invariant full-measure core.

 Let $f\in V_\P^\pm(F)\subset L^1(\Omega,C(M))$. We first apply Lemma \ref{lemma upper bounded representative} with $k=f$ and $c=a>\gamma=\beta_+(f)$. Set $\hat{f}=f+u\circ \Theta-u$, where $u$ is defined in \eqref{eq def u,v}. It follows that the function $u$ is finite, nonnegative, and continuous on $\P$-a.s. fibers. By the second estimate in \eqref{eq:k-bounds}, we have
 \begin{equation}\label{eq bound of g}
   \pi_F(f)-a\leq \pi_F(f)-\hat{f}\leq \pi_F(f)-\min\{f,a\}.
 \end{equation}
Thus the function $g=\pi_F(f)-\hat{f}$ defined in \eqref{eq def u,v} is strictly positive and has integrable fiber norm. As a consequence, $g\in L^1(\Omega,C(M))$. Set $U(\omega)=\max_{x\in M}u(\omega,x)$. Then \eqref{eq:w-sublinear} implies that
 \begin{equation}\label{eq estimate U}
   \lim_{n\to+\infty}\frac{U(\theta^n\omega)}{n}=0,\mbox{ for $\P$-a.s. $\omega\in\Omega$}.
 \end{equation}
Note that $S_ng=n\pi_F(f)-S_nf-u\circ\Theta^n+u.$ Therefore, for $\P$-a.s. $\omega\in\Omega$,
 \begin{equation*}
   \limsup_{n\to+\infty}\frac{1}{n}\sup_{x\in M}|S_ng(\omega,x)-n\pi_F(f)+S_nf(\omega,x)|\leq \limsup_{n\to+\infty}\frac{1}{n}(U(\omega)+U(\theta^n\omega))=0.
 \end{equation*}
Hence
 \begin{equation*}
   \lim_{n\to+\infty}\frac{1}{n}\max_{x\in M}S_ng(\omega,x)=\pi_F(f)-\beta_-(f)<b,\mbox{ for $\P$-a.s. $\omega\in\Omega$.}
   \end{equation*}
Now, we can apply Lemma \ref{lemma upper bounded representative} with $k=g\in L^1(\Omega,C(M))$, $c=b$ and $\gamma=\pi_F(f)-\beta_-(f)$. It follows that the function $v$ defined in \eqref{eq def u,v} is finite, nonnegative, and continuous on $\P$-a.s. fibers, and
   \begin{equation*}
     T_f=g+v\circ \Theta-v\overset{\eqref{eq:k-bounds}}=\min\{b,g+v\circ \Theta\}\in L^1(\Omega,C(M)).
   \end{equation*}
Note that $g\geq \pi_F(f)-a$ by \eqref{eq bound of g}, it follows that
   \begin{equation*}
     \pi_F(f)-a\leq T_f\leq b,
   \end{equation*}
which proves \eqref{eq bound of Tf}. For $V(\omega)=\max_{x\in M}v(\omega,x)$, it follows from \eqref{eq:w-sublinear} that
   \begin{equation}\label{eq estimate V}
     \lim_{n\to+\infty}V(\theta^n\omega)/n=0,\mbox{ for $\P$-a.s. $\omega\in\Omega$}.
   \end{equation}
Telescoping gives $S_nT_f=n\pi_F(f)-S_nf+v\circ\Theta^n-v+u-u\circ\Theta^n$. Therefore, for $\P$-a.s. $\omega\in\Omega$, we have
   \begin{equation*}
      \lim_{n\to+\infty}\frac{1}{n}\sup_{x\in M}|S_nT_f(\omega,x)+S_nf(\omega,x)-n\pi_F(f)|\leq \lim_{n\to+\infty}\frac{U(\omega)+U(\theta^n\omega)+V(\omega)+V(\theta^n\omega)}{n}=0,
        \end{equation*}
which proves \eqref{eq limit Sntf}.

Intersect these invariant cores with those of the Gibbs estimates at a
countable sequence of radii tending to zero, and call the resulting set
$\Omega_*$. The functions above agree there with \eqref{eq def u,v}.
On its complement define $T_f=\pi_F(f)-a$, as in the statement.
This extension is measurable and fiberwise continuous and preserves
\eqref{eq bound of Tf} everywhere.

Next, we show that $T_f$ is simultaneously a lower and an upper energy function for $\nu_f$. Both $T_f$ and $-T_f$ belong to $L^1(\Omega,C(M))$; thus both are upper normal integrands, and $T_f\leq b$ is essentially bounded. Taking $r=c\epsilon_0$ in \eqref{eq error Bowenball}, we obtain
 \begin{equation*}
   \begin{split}
       & \frac{1}{n}\log(\nu_f)_\omega(B_\omega(x,n,r))+\frac{1}{n}S_nT_f(\omega,x)\\
        =&\pi_F(f)-\frac{1}{n}\log\pi_F^*(f)(\omega,r,n)+\frac{1}{n}\left(E_{n,c}(\omega,x)+v(\Theta^n(\omega,x))-v(\omega,x)+u(\omega,x)-u(\Theta^n(\omega,x))\right).
   \end{split}
 \end{equation*}
The terms $\frac{1}{n}\left(E_{n,c}(\omega,x)+v(\Theta^n(\omega,x))-v(\omega,x)+u(\omega,x)-u(\Theta^n(\omega,x))\right)$ tend to $0$ uniformly in $x$ by \eqref{eq estimate error}, \eqref{eq estimate U} and \eqref{eq estimate V}. The identity \eqref{eq fiber top pressure} gives
 \begin{equation*}
 \begin{split}
     &  \lim_{r\downarrow 0}\liminf_{n\to+\infty}\inf_{x}\left\{ \frac{1}{n}\log(\nu_f)_\omega(B_\omega(x,n,r))+\frac{1}{n}S_nT_f(\omega,x)\right\} \\
      =& \pi_F(f)-\lim_{r\downarrow 0}\limsup_{n\to+\infty}\frac{1}{n}\log \pi_F^*(f)(\omega,r,n)=0,
 \end{split}
 \end{equation*}
and
  \begin{equation*}
 \begin{split}
     &  \lim_{r\downarrow 0}\limsup_{n\to+\infty}\sup_{x}\left\{ \frac{1}{n}\log(\nu_f)_\omega(B_\omega(x,n,r))+\frac{1}{n}S_nT_f(\omega,x)\right\} \\
      =& \pi_F(f)-\lim_{r\downarrow 0}\liminf_{n\to+\infty}\frac{1}{n}\log \pi_F^*(f)(\omega,r,n)=0,
 \end{split}
 \end{equation*} Thus $T_f$ satisfies both energy conditions.

 Finally, we prove \eqref{eq int Tf}. Fix any $\rho\in I_\P(\Omega\times M)$. The transfer functions $u$ and $v$ need not be integrable with respect to $\rho$. We therefore use the uniform sublinear estimate instead of integrating the individual coboundaries. Since $T_f\leq b$ and $f\in L^1(\Omega,C(M))\subset L^1(\rho)$, the function $T_f+f-\pi_F(f)\in L^1(\rho)$. The equation \eqref{eq limit Sntf} implies that the Birkhoff averages of the function $T_f+f-\pi_F(f)$ tend to zero $\rho$-almost everywhere. By Birkhoff's ergodic theorem, we have $\int_{\Omega\times M}\bigl(T_f+f-\pi_F(f)\bigr)\,d\rho=0$, which implies \eqref{eq int Tf}. This completes the proof of Theorem \ref{thm simultaneously lower and upper}.
 \end{proof}
\subsection{Some consequences}
A skew product that is Anosov on fibers is fiber expansive with a uniform expansivity constant; see \cite{HLLScichina}. Then by \cite[Theorem 1.3.5]{KL06}, the fiber entropy map $\mu\mapsto h_\mu(F)$ is upper semicontinuous in the narrow topology on $I_\P(\Omega\times M)$.
\begin{corollary}
  Assume the hypotheses of Theorem~\ref{thm gibbsstate energy function}. For every $f\in V_\P^\pm(F)$, let $\nu_f$ be the corresponding random Gibbs state, and let $T_f$ be the energy function for $\nu_f$ constructed in Theorem \ref{thm simultaneously lower and upper}. Then the following conclusions hold on the common full-measure set supplied by Theorem~\ref{proposition normalized good rate function} for this reference measure.
   \begin{itemize}
     \item Quenched level--2 large deviation principle for the reference measure $\nu_f$, as in Theorem \ref{proposition normalized good rate function} with the good rate function $\rho\mapsto I_{T_f}(\rho)$ on $(\operatorname{Pr}(\Omega\times M),\tau)$ defined by
\begin{equation*}
    I_{T_f}(\rho):=
  \begin{cases}
    \displaystyle\pi_F(f)-
    \int_{\Omega\times M}f\,d\rho-h_\rho(F),
       &\rho\in I_\P(\Omega\times M),\\[1.2ex]
    +\infty,
       &\rho\notin I_\P(\Omega\times M).
  \end{cases}
\end{equation*}
     \item Quenched level--1 large deviation principles as in Corollary \ref{corollary quenched level one LDP} for the reference measure $\nu_f$ and any observable $g\in C_b(\Omega\times M)$ with the good rate function $J_g:\mathbb{R}\to [0,+\infty]$ defined by
         \begin{equation*}
            J_g(a):=
  \inf\left\{
 \pi_F(f)- \int_{\Omega\times M}f\,d\rho-h_\rho(F):
  \rho\in I_\P(\Omega\times M),\
  \int_{\Omega\times M}g\,d\rho=a
  \right\},
         \end{equation*}
and the infimum over the empty set is understood to be $+\infty$.
     \item Quenched Laplace principle:  for every
$\Phi\in C_b(\operatorname{Pr}(\Omega\times M))$ and every $\omega\in\Omega_{\mathrm{LD}}$,
\begin{equation*}
  \lim_{n\to+\infty}\frac{1}{n}
  \log\int_M
  \exp\bigl(n\Phi(\E_n(\omega,x))\bigr)d(\nu_f)_\omega(x)=
  \sup_{\rho\in \operatorname{Pr}(\Omega\times M)}
  \{\Phi(\rho)-I_{T_f}(\rho)\}.
\end{equation*}
   \end{itemize}
 \end{corollary}

\section{Preliminaries}\label{sec:preliminaries}
\subsection{Narrow topologies and random probability measures}

By the disintegration theorem, for every $\mu\in \operatorname{Pr}_{\P}(\Omega\times M)$, there exists a measurable disintegration $\omega\mapsto \mu_\omega$, unique up to a $\P$-null set, such that
\begin{equation*}
  \int_{\Omega\times M}h(\omega, x)d\mu(\omega, x)=\int_\Omega\int_Mh(\omega, x)d\mu_\omega(x)d\P(\omega)
\end{equation*}
for every bounded measurable $h: \Omega \times M\rightarrow\mathbb{R}$ (see \cite[Proposition 3.6]{Hans02}).
Conversely, each random probability measure $\omega\mapsto \mu_{\omega}$ defines, through the above equation, a unique Borel probability measure with marginal $\mathbb{P}$ on $(\Omega\times M,\mathcal{F}\otimes \mathcal{B}(M))$.
Owing to the uniqueness of the disintegration, we shall identify $\mu\in \operatorname{Pr}_{\mathbb{P}}(\Omega\times M)$ with its disintegration $\omega\mapsto \mu_\omega$ without further notice; see also Proposition 3.3(ii) and Proposition 3.6 in \cite{Hans02}.

The following lemma from Prohorov theory for random probability measures can be found in \cite[Cor. 4.10]{Hans02} and is critical to the proof of our main theorem.
\begin{lemma}\label{lemma basis of narrow topology}
  The narrow topology $\tau_{rn}$ on $\operatorname{Pr}_\P(\Omega\times M)$ is already generated by the maps $\mu\mapsto \mu(g)$, where $g\in L^1(\Omega,C(M))$ is a random continuous function satisfying
$0\leq g\leq 1$ and
\begin{equation*}
  [g(\omega,\cdot)]_L:=\sup_{x\not=y}\frac{|g(\omega,x)-g(\omega,y)|}{d_M(x,y)}\leq 1,\ \mbox{ for $\P$-a.s. $\omega\in\Omega$}.
\end{equation*}
Such a function $g$ is called a random Lipschitz function with bound $1$.
\end{lemma}

Because the Borel $\sigma$-algebra $\mathcal{F}$ on Polish space $\Omega$ is countably generated modulo $\P$, there is a countable algebra $\{G_m\}_{m\in \mathbb{N}}$ such that, for any $F\in \mathcal{F}$, there exists $G\in \sigma(\{G_m\}_{m\in \mathbb{N}})$ satisfying $\mathbb{P}(F\Delta G)=0$. The narrow topology on $\operatorname{Pr}_{\P}(\Omega\times M)$ can then be metrized by the following metric (see Theorem 4.16 in \cite{Hans02}):
\begin{equation}\label{metric_PrOM}
  d_{rp}(\mu,\nu)=\sum_{m\in\mathbb{N}}\frac{1}{2^m}\sup_{g\in BL(M),0\leq g\leq 1,\ [g]_L\leq 1}\left|\int_{G_m}(\mu_\omega(g)-\nu_\omega(g))\,d\mathbb P(\omega)\right|,
\end{equation}
where
$BL(M)$ is the set of Lipschitz functions on $M$, and $[g]_L$ is the Lipschitz constant of $g$. The absolute value is outside the integral. It does not change the metric, since the test class is closed under $g\mapsto1-g$.

 We denote by $BL_1(\Omega\times M)$ the space of continuous functions $f:\Omega\times M\to [-1,1]$ satisfying $[f]_L\leq 1$, where the metric on the product space $\Omega\times M$ is $d_{\Omega\times M}:=\max\{d_\Omega,d_M\}$.
   The following standard characterization of narrow topology $\tau$ on $\operatorname{Pr}(\Omega\times M)$ can be found in Section 8.3 of \cite{Bogachev}.
  \begin{lemma}\label{lemma ws top}
    The narrow topology $\tau$ on $\operatorname{Pr}(\Omega\times M)$ is generated by $\mu\mapsto\mu(g)$ for $g\in BL_1(\Omega\times M)$. Moreover, since $\Omega\times M$ is separable, the topology $\tau$ is induced by the following metric; see Theorem 8.3.2 in \cite{Bogachev}:
    \begin{equation*}
      d_0(\mu,\nu):=\sup\left\{\int_{\Omega\times M} fd\mu-\int_{\Omega\times M} fd\nu:\ f\in BL_1(\Omega\times M)\right\}.
    \end{equation*}
  \end{lemma}

 \begin{lemma}[An empirical Portmanteau lemma]\label{lemma empirical port lemma}
  Let $\mathcal{G}\subset L^1(\Omega,C(M))$ be a countable subset. There exists a $\theta$-invariant $\P$-full-measure set $\Omega_\mathcal{G}$ such that the following holds. Suppose that
  \begin{itemize}
   \item $\mu_i$ converges to $\mu$ with respect to $(\operatorname{Pr}(\Omega\times M),\tau)$ as $i\to +\infty$,
     \item the marginal measure of $\mu_i$ on $\Omega$ is $\frac{1}{n_i}\sum_{j=0}^{n_i-1}\delta_{\theta^j\omega}$ for some $\omega\in \Omega_\mathcal{G}$, and $n_i\to+\infty$ as $i\to+\infty$,
     \item $\mu\in \operatorname{Pr}_\P(\Omega\times M)$.
   \end{itemize}
Then
   \begin{equation*}
     \int_{\Omega\times M} g d\mu_i\to \int_{\Omega\times M}g d\mu,\ \mbox{ for every $g\in\mathcal{G}$ as $i\to +\infty$.}
   \end{equation*}
\end{lemma}
\begin{proof}[Proof of Lemma \ref{lemma empirical port lemma}]
To obtain the $\theta$-invariant set $\Omega_G$ of full $\mathbb{P}$-measure for a countable set $G$, it suffices to find, for each $g\in G$, a $\theta$-invariant $\P$-full-measure set on which the conclusion holds and then take the intersection over all $g\in G$.

We first suppose that $|g|\leq R$ for some $R>0$. The map $\omega\mapsto g(\omega,\cdot)$ is a Borel map into the separable Banach space $C(M)$: measurability follows by testing on a countable dense subset of $M$. By Lusin's theorem for Borel maps on a Polish probability space (see \cite{Bogachev}), for every $m\in\mathbb N$ there is a compact set $\Omega_{g,m}\subset\Omega$ such that
  \begin{equation*}
    \P(\Omega\setminus\Omega_{g,m})<\frac1m,
  \end{equation*}
and the map $\omega\mapsto g(\omega,\cdot)$ is continuous on $\Omega_{g,m}$ in the uniform norm. The restriction $g|_{\Omega_{g,m}\times M}$ is therefore jointly continuous. By the Tietze extension theorem, we can extend $g|_{\Omega_{g,m}\times M}$ to a bounded jointly continuous function $g_m$ on $\Omega\times M$ with $\|g_m\|_\infty\leq R$. By Birkhoff's ergodic theorem, there exists a $\theta$-invariant $\P$-full-measure set such that
  \begin{equation}\label{eq choosing omega}
    \lim_{n\to+\infty}\frac{1}{n}\sum_{i=0}^{n-1}1_{\Omega\backslash \Omega_{g,m}}(\theta^i\omega)=\P(\Omega\backslash \Omega_{g,m})<\frac{1}{m},\mbox{ for all $m\geq 2$ and $\P$-a.s. $\omega\in\Omega$}.
  \end{equation}
Now, for any measure $\mu_i\in \operatorname{Pr}(\Omega\times M)$ with marginal measure  $\frac{1}{n_i}\sum_{j=0}^{n_i-1}\delta_{\theta^j\omega}$ on $\Omega$ for some $\omega$ satisfying \eqref{eq choosing omega}, we have
  \begin{equation*}
    \limsup_{i\to +\infty} |\int g d\mu_i-\int g_m d\mu_i|\leq 2R\P(\Omega\backslash \Omega_{g,m})<\frac{2R}{m},
  \end{equation*}
and
  \begin{equation*}
    |\int g_m d\mu-\int gd\mu|\leq 2R\cdot \mu (\Omega\times M\backslash(\Omega_{g,m}\times M))<\frac{2R}{m}
  \end{equation*}
  Therefore,
  \begin{equation*}
  \begin{split}
  &\limsup_{i\to +\infty}|\int g d\mu_i-\int g d\mu|   \\
  \leq &  \limsup_{i\to +\infty}( |\int g d\mu_i-\int g_m d\mu_i|+|\int g_m d\mu_i-\int g_m d\mu|+|\int g_m d\mu-\int gd\mu|)\\
       < & \frac{4R}{m}.
  \end{split}
  \end{equation*}
Letting $m\to+\infty$ and taking the countable intersection of the $\theta$-invariant $\P$-full-measure sets completes the proof of the bounded case.

  We now remove the boundedness assumption on \(g\). We write $
    b_g(\omega):=\sup_{x\in M}|g(\omega,x)|.
$
By the definition of $L^1(\Omega, C(M))$, we have
\(b_g\in L^1(\mathbb P)\). For \(R\in\mathbb N\),  we define the
truncation function
\begin{equation*}
  g^{(R)}(\omega,x)
    :=\max\{-R,\min\{g(\omega,x),R\}\}\in L^1(\Omega,C(M)).
\end{equation*}
The bounded case proved above applies to \(g^{(R)}\). For every \((\omega,x)\in\Omega\times M\), we have
\begin{equation*}
  |g(\omega,x)-g^{(R)}(\omega,x)|
    =\bigl(|g(\omega,x)|-R\bigr)^+
    \le \bigl(b_g(\omega)-R\bigr)^+.
\end{equation*}
By Birkhoff's ergodic theorem, we may further choose $\omega$ from a $\theta$-invariant $\P$-full-measure set such that
\begin{equation*}
  \lim_{n\to +\infty}\frac{1}{n}\sum_{j=0}^{n-1}(b_g(\theta^j\omega)-R)^+=\int_\Omega (b_g(\omega)-R)^+d\P(\omega).
\end{equation*}
It follows that
\begin{equation*}
 \limsup_{i\to +\infty} |\int gd\mu_i-\int g^{(R)}d\mu_i|\leq  \limsup_{i\to +\infty} \frac{1}{n_i}\sum_{j=0}^{n_i-1}(b_g(\theta^j\omega)-R)^+=\int_\Omega(b_g(\omega)-R)^+d\P(\omega).
\end{equation*}
Using the triangle inequality and the bounded case, we obtain
\begin{equation*}
\begin{split}
    & \limsup_{i\to +\infty} |\int gd\mu_i-\int g d\mu| \\
    \leq  & \limsup_{i\to +\infty} (|\int gd\mu_i-\int g^{(R)} d\mu_i|+|\int g^{(R)}d\mu_i-\int g^{(R)} d\mu|+|\int gd\mu-\int g^{(R)} d\mu|)\\
    \leq & 2\int_\Omega(b_g(\omega)-R)^+d\P(\omega).
\end{split}
\end{equation*}
Finally,
\[
    0\le (b_g-R)^+\le b_g,
    \qquad
    (b_g-R)^+\longrightarrow0
\]
pointwise as \(R\to+\infty\). Since \(b_g\in L^1(\mathbb P)\), the
dominated convergence theorem implies
$
    \int_\Omega(b_g-R)^+\,d\mathbb P\to 0
$ as $ R\to +\infty$.
Letting \(R\to\infty\) and taking the countable intersection of the $\theta$-invariant $\P$-full-measure sets, we conclude that
$
    \int_{\Omega\times M}gd\mu_i\to
    \int_{\Omega\times M}gd\mu
$ as $i\to +\infty$. This completes the proof of Lemma \ref{lemma empirical port lemma}.
\end{proof}

\begin{lemma}\label{lemma lemma diffrent converges}
  If $\{\mu,\mu_n\}_{n\in \mathbb{N}}\subset \operatorname{Pr}_\P(\Omega\times M)$, and $\mu_n\to \mu$ with respect to $(\operatorname{Pr}(\Omega\times M),\tau)$, then $\mu_n$ also converges to $\mu$ with respect to $(\operatorname{Pr}_\P(\Omega\times M),\tau_{rn})$. As a consequence, $\tau$ and
$\tau_{rn}$ induce the same topology on $\operatorname{Pr}_\P(\Omega\times M)$.
\end{lemma}
\begin{proof}[Proof of Lemma \ref{lemma lemma diffrent converges}]
For the first part, by lemma \ref{lemma basis of narrow topology}, it suffices to show that $\int g\,d\mu_n\to \int g\,d\mu$ for any fixed $g\in L^1(\Omega,C(M))$ with $0\leq g\leq 1$. The proof is much simpler than that of Lemma \ref{lemma empirical port lemma}, since each $\mu_n$ has marginal $\P$ on $\Omega$.

As in the proof of Lemma~\ref{lemma empirical port lemma}, Lusin's theorem applied to the $C(M)$-valued map $\omega\mapsto g(\omega,\cdot)$ gives, for each $\epsilon\in(0,1/2)$, a compact set $\Omega_\epsilon\subset\Omega$ such that
  \begin{equation*}
    \P(\Omega\backslash \Omega_{\epsilon})<\epsilon\mbox{ and $g|_{\Omega_{\epsilon}\times M}$ is jointly continuous.}
  \end{equation*}
By the Tietze extension theorem, we can extend $g|_{\Omega_{\epsilon}\times M}$ to a bounded jointly continuous function $g_\epsilon$ on $\Omega\times M$ with $\|g_\epsilon\|_\infty\leq 1$. We then have
  \begin{equation*}
    \limsup_{i\to +\infty} |\int g d\mu_i-\int g_\epsilon d\mu_i|\leq 2\P(\Omega\backslash \Omega_{\epsilon})<2\epsilon,
  \end{equation*}
and
  \begin{equation*}
    |\int g_\epsilon d\mu-\int gd\mu|\leq 2\cdot \mu (\Omega\times M\backslash(\Omega_{\epsilon}\times M))<2\epsilon.
  \end{equation*}
  Therefore,
  \begin{equation*}
  \begin{split}
  &\limsup_{i\to +\infty}|\int g d\mu_i-\int g d\mu|   \\
  \leq &  \limsup_{i\to +\infty}( |\int g d\mu_i-\int g_\epsilon d\mu_i|+|\int g_\epsilon d\mu_i-\int g_\epsilon d\mu|+|\int g_\epsilon d\mu-\int gd\mu|)\\
       < & 4\epsilon.
  \end{split}
  \end{equation*}
Letting $\epsilon\to 0$ completes the proof.

For the second part, we note that $C_b(\Omega\times M)\subset L^1(\Omega,C(M))$.
Therefore, convergence with respect to $\tau_{rn}$ implies convergence
with respect to the restriction of $\tau$ to $\operatorname{Pr}_\P(\Omega\times M)$.
The converse follows from
the first part.  Consequently, $\tau$ and
$\tau_{rn}$ induce the same topology on $\operatorname{Pr}_\P(\Omega\times M)$.
This completes the proof of Lemma \ref{lemma lemma diffrent converges}.
\end{proof}
The following is the classical Prohorov Theorem.
\begin{lemma}[Classical Prohorov Theorem]\label{lemma classical prohorov thm}
  A set of probability measures $G\subset \operatorname{Pr}(X)$ on a Polish space $X$ is precompact in the narrow topology if and only if it is tight. Here, tightness means that for every $\epsilon>0$, there exists a compact set $K_\epsilon\subset X$ with $\rho(K_\epsilon)\geq 1-\epsilon$ for all $\rho\in G$.
  \end{lemma}
\begin{remark}\label{remark compactness}
  The classical Prohorov theorem implies that $\operatorname{Pr}_\P(\Omega\times M)$ is compact in $\operatorname{Pr}(\Omega\times M)$ with respect to $\tau.$ Indeed, the singleton $\{\P\}$ is tight in $\operatorname{Pr}(\Omega)$, therefore there exists a compact set $\Omega_\epsilon\subset \Omega$ such that $\P(\Omega_\epsilon)\geq 1-\epsilon$. Then for any measure $\mu\in \operatorname{Pr}_\P(\Omega\times M)$, $\mu(\Omega_\epsilon\times M)=\P(\Omega_\epsilon)\geq 1-\epsilon$. The set $\operatorname{Pr}_\P(\Omega\times M)$ is also closed: it is the inverse image of $\{\P\}$ under the continuous push-forward map induced by $\pi_\Omega$. Tightness and Prohorov's theorem therefore imply compactness in $(\operatorname{Pr}(\Omega\times M),\tau)$. The set $I_\P(\Omega\times M)$ is closed in
$(\operatorname{Pr}_\P(\Omega\times M),\tau)$ and hence compact.
Indeed, suppose that $\rho_n\in I_\P(\Omega\times M)$
and $\rho_n\to\rho$ in $\tau$.
By Lemma~\ref{lemma lemma diffrent converges},
the convergence also holds in $\tau_{rn}$.
For every $g\in L^1(\Omega,C(M))$,
we have $g\circ\Theta\in L^1(\Omega,C(M))$. Indeed,
\[
 \sup_{x\in M}|g(\Theta(\omega,x))|
 \leq\|g(\theta\omega,\cdot)\|_\infty,
\]
and $\theta$ preserves $\P$. Consequently,
\begin{equation*}
  \int g\circ\Theta\,d\rho
=\lim_{n\to \infty}\int g\circ\Theta\,d\rho_n
=\lim_{n\to \infty}\int g\,d\rho_n
=\int g\,d\rho.
\end{equation*}
Since these functions separate the measures in $\operatorname{Pr}_\P(\Omega\times M)$, we have $\Theta_*\rho=\rho$.
\end{remark}

\subsection{Energy functions}\label{subsection energy function}
We now give a natural sufficient condition for an upper normal integrand.
\begin{lemma}\label{lemma check integrand}
  If an $\mathcal{F}\otimes\mathcal{B}(M)$-measurable function $T:\Omega\times M\to\mathbb{R}$ satisfies
  \begin{enumerate}
    \item for $\P$-a.s. $\omega\in\Omega$, $x\mapsto T(\omega,x)$ is upper semicontinuous;
    \item $b(\omega):=\sup_{x\in M}T(\omega,x)\in L^1(\P)$.
  \end{enumerate}
Then $T$ is an upper normal integrand.
\end{lemma}
\begin{proof}[Proof of Lemma \ref{lemma check integrand}]

Let $\Omega_0\subset\Omega$ be a measurable set with $\mathbb P(\Omega_0)=1$ such that
$x\mapsto T(\omega,x)$ is upper semicontinuous on $M$ and
$\sup_{x\in M}T(\omega,x)<\infty$ for every $\omega\in\Omega_0$.
By redefining $T$ on $\Omega\setminus\Omega_0$ by setting $T(\omega,x)=0$, we may assume that
$x\mapsto T(\omega,x)$ is upper semicontinuous and bounded above on $M$ for every
$\omega\in\Omega$. For each $k\in\mathbb N$, define
\[
\phi_k(\omega,x):=\sup_{y\in M}\{T(\omega,y)-k d_M(x,y)\}.
\]
Choosing $y=x$ in the supremum gives $\phi_k(\omega,x)\geq T(\omega,x)$.

We first show that $\phi_k$ is measurable. For fixed $k$, the function $
(\omega,x,y)\mapsto T(\omega,y)-k d_M(x,y)$
is measurable with respect to $\mathcal F\otimes\mathcal B(M)\otimes\mathcal B(M)$, and for
each fixed $(\omega,x)$, the map
$
y\mapsto T(\omega,y)-k d_M(x,y)
$
is upper semicontinuous on the compact metric space $M$. For any $a\in\mathbb{R}$, we have the following equality
\begin{equation*}
  \{(\omega,x)\in\Omega\times M:\ \phi_k(\omega,x)>a\}=   proj_{\Omega\times M}\{(\omega,x,y)\in\Omega\times M\times M:\ T(\omega,y)-kd_M(x,y)>a\}.
\end{equation*}
Let $A=\{(\omega,x,y)\in\Omega\times M\times M:\ T(\omega,y)-kd_M(x,y)>a\}$. For any $(\omega,x)\in\Omega\times M$, the $(\omega,x)$-section of $A$
\begin{equation*}
  A_{(\omega,x)}:=\{y\in M:\ T(\omega,y)-kd_M(x,y)>a\}=\cup_{m=1}^\infty\{y\in M:\ T(\omega,y)-kd_M(x,y)\geq a+\frac{1}{m}\}
\end{equation*}
is $\sigma$-compact by upper semicontinuity. By the projection theorem for measurable sets with $\sigma$-compact sections (see, e.g., \cite[Theorem 5.12.1]{GTM180}), the set $\{(\omega,x)\in\Omega\times M:\ \phi_k(\omega,x)>a\}$
is $\mathcal F\otimes\mathcal B(M)$-measurable. Therefore, $\phi_k$ is measurable.

Next, for each fixed $\omega$, the map $x\mapsto \phi_k(\omega,x)$ is $k$-Lipschitz. Indeed,
for any $x,x'\in M$,
\begin{equation*}
  |\phi_k(\omega,x)-\phi_k(\omega,x') |\leq \sup_{y\in M} k |d_M(x',y)-d_M(x,y)|\leq  k d_M(x,x').
\end{equation*}

Let $b(\omega)=\sup_{y\in M}T(\omega,y)$. Then
\[
 b(\omega)-k\operatorname{diam}M\le \phi_k(\omega,x)\le b(\omega).
\]
Since $b\in L^1(\mathbb P)$, it follows that
$
\phi_k\in L^1(\Omega,C(M)).$

It remains to prove that $\phi_k(\omega,x)\searrow T(\omega,x)$. The monotonicity is immediate
from the definition: if $k_1<k_2$, then
$
T(\omega,y)-k_2d_M(x,y)\le T(\omega,y)-k_1d_M(x,y),
$
hence
$
\phi_{k_2}(\omega,x)\le \phi_{k_1}(\omega,x).
$

Fix $\omega\in\Omega_0$ and $x\in M$. Let
\[
L(\omega,x):=\lim_{k\to\infty}\phi_k(\omega,x).
\]
We already know that $L(\omega,x)\ge T(\omega,x)$. Suppose, for contradiction, that
$
L(\omega,x)>T(\omega,x).
$
Then there exists $\eta>0$ such that for infinitely many $k$,
$
\phi_k(\omega,x)>T(\omega,x)+2\eta.
$
For each such $k$, by definition of $\phi_k$, we can choose $y_k\in M$ satisfying
\begin{equation}\label{eq geq T+eta}
  T(\omega,y_k)-k d_M(x,y_k)>T(\omega,x)+\eta.
\end{equation}

Since $T(\omega,\cdot)\le b(\omega)$, we get
$
k d_M(x,y_k)< b(\omega)-T(\omega,x)-\eta,
$
and therefore $
d_M(x,y_k)\to0 $ as $k \to +\infty.$
Thus $y_k$ converges to $x$. Since $T(\omega,\cdot)$ is upper semicontinuous,
\[
\limsup_{k\to\infty}T(\omega,y_k)\le T(\omega,x).
\]
This contradicts \eqref{eq geq T+eta}.
Therefore
$
L(\omega,x)\le T(\omega,x).
$
Hence
$
\phi_k(\omega,x)\searrow T(\omega,x)
$
for every $\omega\in\Omega_0$ and every $x\in M$. This completes the proof.

\end{proof}

\begin{lemma}\label{lemma upper semi-cotinuous}
  If $T:\Omega\times M\to\mathbb{R}_{\geq 0}$ is an upper normal integrand, then $\rho\mapsto \int_{\Omega\times M}Td\rho$ is upper semicontinuous on $\operatorname{Pr}_\P(\Omega\times M)$ with respect to $(\operatorname{Pr}_\P(\Omega\times M),\tau_{rn})$. As a consequence, if $T$ is both an lower and an upper energy function for some measure $\nu\in \operatorname{Pr}_\P(\Omega\times M)$, then $\rho\mapsto \int_{\Omega\times M}Td\rho$ is continuous on $\operatorname{Pr}_\P(\Omega\times M)$ with respect to $(\operatorname{Pr}_\P(\Omega\times M),\tau_{rn})$.
\end{lemma}
\begin{proof}[Proof of Lemma \ref{lemma upper semi-cotinuous}]
Let $\rho_n\in \operatorname{Pr}_\P(\Omega\times M)$ be a sequence of random probability measures converging to $\rho \in \operatorname{Pr}_\P(\Omega\times M)$ with respect to the narrow topology $\tau_{rn}$. Let $\phi_k\in L^1(\Omega,C(M))$ satisfy $\phi_k\searrow T$ for $\P$-a.s. $\omega\in\Omega$ and all $x\in M$. Since $\phi_k\in L^1(\Omega,C(M))$, we have $\int\phi_k d\rho_n\to \int \phi_k d\rho$ as $n\to +\infty$. Since $T\leq \phi_k$ for all $k$, $\int Td\rho_n\leq \int\phi_k d\rho_n$. Taking the limit superior gives, for every $k$,
\begin{equation*}
  \limsup_{n\to +\infty} \int Td\rho_n\leq \limsup_{n\to +\infty} \int \phi_kd\rho_n=\int \phi_k d\rho.
\end{equation*}
Since the nonnegative functions $\phi_k\geq T$ decrease to $T$ and $\phi_1\in L^1(\Omega\times M,\rho)$, the dominated convergence theorem gives
\begin{equation*}
  \limsup_{n\to +\infty} \int Td\rho_n\leq \inf_k\int \phi_k d\rho=\int Td\rho.
  \end{equation*}
The upper semicontinuity is proved.

 If $T$ is both a lower and an upper energy function, we can choose a constant $B$ with $0\leq T\leq B$ on almost every fiber. Both $T$ and $B-T$ are nonnegative upper normal integrands.
Applying the first part to these two functions proves both upper and lower semicontinuity of $\rho\mapsto\int T\,d\rho$.
\end{proof}
We now prove Proposition \ref{proposition lower energy function}.
\begin{proof}[Proof of Proposition \ref{proposition lower energy function}]
 First, assume that $T$ is a lower energy function for $\nu$. Then $T$ is an upper normal integrand, and therefore there exists a sequence $\phi_k\in L^1(\Omega,C(M))$ such that $\phi_k(\omega,x)\searrow T(\omega,x)$ for $\P$-a.s. $\omega\in\Omega$ and all $x\in M$. Set $b(\omega):=\sup_{x\in M}|\phi_1(\omega,x)|$. Then $b\in L^1(\P)$ by the definition of $L^1(\Omega,C(M))$. After restricting to a $\theta$-invariant $\P$-full-measure set, we may assume that
 \begin{equation*}
   0\leq T(\Theta^l(\omega,x))\leq b(\theta^l\omega),\ \forall l\in\mathbb{Z} \mbox{ and $\P$-a.s. $\omega\in\Omega$}.
 \end{equation*}
Since $\theta$ is $\P$-ergodic, we have
 \begin{equation*}
   \lim_{n\to +\infty}\frac{1}{n}\sum_{j=0}^{n-1}b(\theta^j\omega)=\int_\Omega b(\omega)d\P(\omega)<+\infty,\mbox{ for $\P$-a.s. $\omega\in\Omega$.}
 \end{equation*}
The following inequality then holds for $\P$-a.s. $\omega\in\Omega$:
 \begin{equation*}
   \inf_{x\in M}\frac{1}{n}\log \nu_\omega (B_\omega(x,n,\epsilon))\geq \inf_{x\in M}\{\frac{1}{n}\log \nu_\omega (B_\omega(x,n,\epsilon))+\int Td(\E_n(\omega,x))\}-\frac{1}{n}\sum_{j=0}^{n-1}b(\theta^j\omega).
 \end{equation*}
Hence, for $\P$-a.s. $\omega\in\Omega$,
 \begin{equation*}
 \begin{split}
    &\lim_{\epsilon\searrow0}\liminf_{n\to+\infty}\inf_{x\in M}\{\frac{1}{n}\log \nu_\omega (B_\omega(x,n,\epsilon))\}   \\
\geq      & \lim_{\epsilon\searrow0}\liminf_{n\to+\infty}\inf_{x\in M}\{\frac{1}{n}\log \nu_\omega (B_\omega(x,n,\epsilon))+\int Td(\E_n(\omega,x))\}-\int_\Omega b(\omega)d\P(\omega)\\
\geq &-\int b(\omega)d\P(\omega), \mbox{ by \eqref{eq def lower energy}.}
 \end{split}
 \end{equation*}
Consequently, \eqref{eq proposition exists energy} holds.

 Conversely, assume that \eqref{eq proposition exists energy} holds. Then there exists a constant $C>0$ such that
 \begin{equation*}
 \lim_{\epsilon\searrow0}\liminf_{n\to+\infty}\inf_{x\in M}\{\frac{1}{n}\log \nu_\omega (B_\omega(x,n,\epsilon))\}\geq -C
 \end{equation*}
for $\P$-a.s. $\omega\in\Omega$. Let $T_C\equiv C$. Then $T_C$ is a nonnegative random continuous function and $\int T_Cd(\E_n(\omega,x))=C$. Moreover,
 \begin{equation*}
   \lim_{\epsilon\searrow0}\liminf_{n\to+\infty}\inf_{x\in M}\{\frac{1}{n}\log \nu_\omega (B_\omega(x,n,\epsilon))+\int T_Cd(\E_n(\omega,x))\}\geq 0
 \end{equation*}
for $\P$-a.s. $\omega\in\Omega.$ Hence $T_C$ is a lower energy function for $\nu$.

 Next, assume that $U:\Omega\times M\to\mathbb{R}_{\geq 0}$ is a nondegenerate upper energy function for $\nu$. By definition, $-U$ is an upper normal integrand, and therefore $x\mapsto U(\omega,x)$ is lower semicontinuous on a $\P$-full-measure set $\Omega_U$. Define
 \begin{equation*}
   m(\omega):=\begin{cases}
               \min_{x\in M} U(\omega,x) , & \mbox{for $\omega\in\Omega_U$}  \\
                0, & \mbox{otherwise}.
              \end{cases}
 \end{equation*}
We claim that $\omega\mapsto m(\omega)$ is measurable. Indeed, for any $a>0$,
 \begin{equation*}
   \{\omega\in\Omega:\ m(\omega)\leq a\}=(Proj_\Omega\{(\omega,x)\in\Omega_U\times M:\ U(\omega,x)\leq a\})\cup (\Omega\backslash\Omega_U),
 \end{equation*}
which is measurable by the projection theorem for measurable sets with $\sigma$-compact sections; see, e.g., \cite[Theorem 5.12.1]{GTM180}. Since $$\operatorname*{ess\,sup}_{\omega}\sup_{x\in M}U(\omega,x)<+\infty,$$ we have $m\in L^1(\P)$. For $\P$-a.s. $\omega\in\Omega$, we then have
 \begin{equation*}
   \sup_{x\in M}\frac{1}{n}\log \nu_\omega (B_\omega(x,n,\epsilon))\leq \sup_{x\in M}\{\frac{1}{n}\log \nu_\omega (B_\omega(x,n,\epsilon))+\int Ud(\E_n(\omega,x))\}-\frac{1}{n}\sum_{j=0}^{n-1}m(\theta^j\omega).
 \end{equation*}
Thus, for $\P$-a.s. $\omega\in\Omega$,
 \begin{equation*}
   \begin{split}
 &\lim_{\epsilon\searrow 0}\limsup_{n\to+\infty}\sup_{x\in M}\frac{1}{n}\log \nu_\omega (B_\omega(x,n,\epsilon))\\
 \leq &\lim_{\epsilon\searrow 0}\limsup_{n\to+\infty}\sup_{x\in M}\{\frac{1}{n}\log \nu_\omega (B_\omega(x,n,\epsilon))+\int Ud(\E_n(\omega,x))\}-\int_\Omega m(\omega)d\P(\omega)\\
 \leq & -\int_\Omega m(\omega)d\P(\omega)<0,\mbox{ by \eqref{eq def upper energy}}.
   \end{split}
 \end{equation*}
Consequently, \eqref{eq proposition exists upper energy} holds. The converse is proved as in the preceding case. This completes the proof.
\end{proof}

\begin{lemma}\label{lemma lower energy finite entropy}
Let $\nu\in \operatorname{Pr}_\P(\Omega\times M)$ and suppose that
$T:\Omega\times M\to\mathbb R_{\geq0}$ is a lower energy
function for $\nu$. Then there exists a constant $C_T<\infty$
such that
\begin{equation}\label{eq uniform finite entropy}
 h_{\mathrm{top}}(F)=\sup_{\rho\in I_\P(\Omega\times M)}h_\rho(F)\leq C_T<\infty
\end{equation}
\end{lemma}
\begin{proof}[Proof of Lemma \ref{lemma lower energy finite entropy}]
By Proposition~\ref{proposition lower energy function}, there exists
a constant $C_T<\infty$ such that, for $\P$-a.s.\ $\omega\in\Omega$,
\[
 \lim_{\epsilon\searrow0}\liminf_{n\to\infty}\inf_{x\in M}
 \frac1n\log\nu_\omega(B_\omega(x,n,\epsilon))
 \geq -C_T.
\]
Fix such an $\omega$ for which \eqref{eq fiber top pressure} also
holds with $f=0$, and let $\epsilon,\eta>0$.
By monotonicity with respect to the radius, the same lower bound
holds for every fixed $\epsilon>0$. Hence there exists
$N=N(\omega,\epsilon,\eta)$ such that, for all $n\geq N$ and $x\in M$,
\begin{equation}\label{eq uniform lower Bowen mass}
 \nu_\omega(B_\omega(x,n,\epsilon))
 \geq \exp\bigl(-n(C_T+\eta)\bigr).
\end{equation}

Let $E_n(\omega,2\epsilon)\subset M$ be an
$(\omega,n,2\epsilon)$-separated set of maximal cardinality.
The Bowen balls $B_\omega(x,n,\epsilon)$, $x\in E_n(\omega,2\epsilon)$,
are pairwise disjoint. Therefore, for every $n\geq N$,
\[
 1\geq
 \sum_{x\in E_n(\omega,2\epsilon)}
 \nu_\omega(B_\omega(x,n,\epsilon))
 \geq
 \#E_n(\omega,2\epsilon)\exp\bigl(-n(C_T+\eta)\bigr),
\]
and consequently
\[
 \limsup_{n\to\infty}\frac1n\log\#E_n(\omega,2\epsilon)
 \leq C_T+\eta.
\]
Letting $\eta\searrow0$ and then $\epsilon\searrow0$, and using
\eqref{eq fiber top pressure} with $f=0$, we obtain
\[
 h_{\mathrm{top}}(F)
 =\lim_{\epsilon\searrow0}\limsup_{n\to\infty}
   \frac1n\log\#E_n(\omega,2\epsilon)
 \leq C_T<\infty.
\]
The random variational principle
\cite[Theorem~1.2.13]{KL06} now gives
\eqref{eq uniform finite entropy}.
\end{proof}

\subsection{Measure approximation}\label{section 3.4}
The relative entropy $h_\mu^{(r)}(\Theta)$ of $\Theta$ which is also called the fiber measure-theoretic entropy $h_\mu(F)$ of RDS $F$ with respect to $\mu\in I_\P(\Omega\times M)$ is defined by
\begin{equation*}
  h_\mu^{(r)}(\Theta)=h_\mu(F)=\sup_{\mathcal{Q}}h_\mu(F,\mathcal{Q}),
\end{equation*}
where the supremum is taken over all finite measurable partitions $\mathcal{Q}$ of $\Omega\times M$,
\begin{equation*}
  h_\mu(F,\mathcal{Q})=\lim_{n\to\infty}\frac{1}{n}H_\mu(\vee_{i=0}^{n-1}(\Theta^i)^{-1}\mathcal{Q}|\pi_\Omega^{-1}\mathcal{F}),
\end{equation*}
and $\pi_\Omega:\Omega\times M\to\Omega$ is the projection onto the first coordinate.

 The following is a measure-approximation lemma for random dynamical systems; its proof is similar to that of \cite[p.~535]{Young1990}.
\begin{lemma}[Measure approximation]\label{lemma measure approxi}
Assume that $(\Omega,\mathcal{F},\P,\theta)$ is ergodic and that $M$ is compact. Given a finite sequence of random continuous functions $\{f_i\}_{i=1}^n\subset L^1(\Omega,C(M))$, for any $\mu\in I_\P(\Omega\times M)$ with $h_\mu(F)<\infty$ and every $\delta,\gamma>0$, there exist $\nu\in I_\P(\Omega\times M)$ and $\{\nu_j\}_{j=1}^k\subset I_\P^e(\Omega\times M)$, all of finite fiber entropy, with the following properties:
  \begin{enumerate}
    \item $\nu=\sum_{j=1}^{k}\lambda_j\nu_j$, where $\lambda_j>0$, $\sum_{j=1}^{k}\lambda_j=1$;
    \item $h_\nu(F)\geq h_\mu(F)-\gamma$,
    \item $|\int_{\Omega\times M} f_i d\mu-\int_{\Omega\times M}f_i d\nu|<\delta$ for all $i=1,\ldots,n$.
  \end{enumerate}
\end{lemma}
\begin{proof}[Proof of Lemma \ref{lemma measure approxi}]
By Remark~\ref{remark compactness} and Lemma~\ref{lemma lemma diffrent converges},
$I_\P(\Omega\times M)$ is compact in the random narrow topology. Let $\pi$
be the ergodic decomposition of $\mu$. Since the base measure $\P$ is ergodic,
$\pi$-almost every ergodic component has marginal $\P$; compare
\cite[Lemma~6.19]{Hans02}. For $f\in L^1(\Omega,C(M))$,
\[
 \mu(f)=\int_{I_\P^e(\Omega\times M)}\tau(f)\,d\pi(\tau),\qquad
 h_\mu(F)=\int_{I_\P^e(\Omega\times M)}h_\tau(F)\,d\pi(\tau).
\]
The second identity is the ergodic decomposition formula for relative
entropy, not merely its affinity for finite convex combinations. In
particular, $h_\tau(F)<\infty$ for $\pi$-almost every $\tau$.

The finitely many maps $\tau\mapsto\tau(f_j)$ are uniformly continuous on
this compact space. Choose $\epsilon_0>0$ so that
$d_{rp}(\tau_1,\tau_2)<\epsilon_0$ implies
$|\tau_1(f_j)-\tau_2(f_j)|<\delta$ for all $j$.
Take a finite Borel partition into sets of diameter less than $\epsilon_0$,
discard its $\pi$-null atoms, and denote the remaining atoms by $Q_1,\ldots,Q_k$.
Put $\lambda_i=\pi(Q_i)>0$. In the full-measure subset of finite-entropy
ergodic components, choose $\nu_i\in Q_i$ such that
\[
 h_{\nu_i}(F)\geq\frac1{\lambda_i}\int_{Q_i}h_\tau(F)\,d\pi(\tau)-\gamma.
\]
This choice is possible because the average is finite. Set
$\nu=\sum_i\lambda_i\nu_i$. Affinity for this finite combination gives
\[
 h_\nu(F)=\sum_i\lambda_i h_{\nu_i}(F)
 \geq h_\mu(F)-\gamma.
\]
For every $j$,
\[
 |\mu(f_j)-\nu(f_j)|
 \leq\sum_i\int_{Q_i}|\tau(f_j)-\nu_i(f_j)|\,d\pi(\tau)<\delta.
\]
This proves all the assertions.
\end{proof}

\subsection{Random Katok entropy theorem}

For $\epsilon>0$, $\delta\in(0,1)$, $\omega\in\Omega$, and $\mu\in I_\P(\Omega\times M)$, define
\begin{equation*}
  S_\mu(\omega,n,\epsilon,\delta)=\min\{S_\mu(\omega,n,\epsilon,K)|\ K\subset M,\ \mu_\omega(K)\geq 1-\delta\},
\end{equation*}
where $S_\mu(\omega,n,\epsilon,K)$ denotes the smallest cardinality of any $(\omega,n,\epsilon)$-spanning set of $K$. A subset $G\subset K$ is called an $(\omega,n,\epsilon)$-spanning set of $K$ if, for any $y\in K$, there exists an $x\in G$ such that $d_{\omega}^n(x,y)\leq \epsilon$. The following is the random version of Katok's entropy theorem.
The upper-limit formula is given in
\cite[Theorem~3.1]{zhutwonotes}.
Both entropy identities also follow from the weighted Katok entropy
formula in \cite[Theorem~8.2]{YangChenLinZhou2025}.
\begin{lemma}\label{lemma katok entropy}
   For any $\mu\in I^e_\P(\Omega\times M)$ with $h_{\mu}(F)<\infty$, $\epsilon>0$, and $0<\delta<1$. Then, for $\P$-a.s. $\omega\in\Omega$,
  \begin{align}\label{eq katok delta}
    h_\mu(F)&=\lim_{\epsilon\rightarrow 0}\limsup_{n\rightarrow\infty}\frac{1}{n}\log S_\mu(\omega,n,\epsilon,\delta)=\lim_{\epsilon\rightarrow 0}\liminf_{n\rightarrow\infty}\frac{1}{n}\log S_\mu(\omega,n,\epsilon,\delta).
  \end{align}
\end{lemma}

\section{Proof of the large deviation upper bound}\label{sec:upper-bound}
In this section, we prove Theorem \ref{thm 2}. We write
$\P_{n,\omega}=n^{-1}\sum_{j=0}^{n-1}\delta_{\theta^j\omega}$
for the empirical base marginal.
For a nonempty subset $D\subset \operatorname{Pr}(\Omega\times M)$ and $\epsilon>0$, let $\Gamma_{D}(\omega,n,\epsilon)$ be an $(\omega,n,\epsilon)$-separated subset of $\{x\in M:\E_{n}(\omega,x)\in D\}$ with the maximal cardinality. Let $CH(D)$ denote the convex hull of $D$, and let $\overline{CH(D)}$ be its closure in $(\operatorname{Pr}(\Omega\times M),\tau)$.

\begin{lemma}\label{lemma of Thm 2}
There exists a $\P$-full-measure set $\Omega_{sep}$ such that for every $\omega\in\Omega_{sep},$ any $\epsilon>0$ and every  subset $D\subset \operatorname{Pr}(\Omega\times M)$, either $\Gamma_D(\omega,n,\epsilon)$ is empty for all sufficiently large $n$, or
there exists a measure $$\rho \in \overline{CH(D)}\cap I_\P(\Omega\times M)$$ satisfying
  \begin{equation*}
    \limsup_{n\to+\infty}\frac{1}{n}\log \#\Gamma_D(\omega,n,\epsilon)\leq  h_\rho(F).
  \end{equation*}
  If, furthermore, the fiber entropy map $\mu\mapsto h_\mu(F)$ is upper semicontinuous on $I_\P(\Omega\times M)$, then we can choose $\rho \in \overline{D}\cap I_\P(\Omega\times M)$ in the second alternative.
\end{lemma}
\begin{proof}[Proof of Lemma \ref{lemma of Thm 2}]
To construct $\rho\in I_\P(\Omega\times M)$ by using fiberwise information, we need to use the empirical Portmanteau lemma \ref{lemma empirical port lemma}.
We first choose an admissible countable family for the empirical Portmanteau lemma.

   Let $\{G_m\}_{m\in\mathbb N}$ be the countable algebra introduced
in the definition of $d_{rp}$, which generates $\mathcal F$ modulo
$\P$. Since $M$ is a compact metric space, there exists a countable
subset $\{\varphi_r\}_{r\in\mathbb N}\subset BL_1(M)$ which determines
Borel probability measures on $M$. For $m,r\in\mathbb N$, we define
\begin{equation*}
  g_{m,r}(\omega,x):=1_{G_m}(\omega)\varphi_r(x),
\end{equation*}
and let
$
\mathcal G_0:=\{g_{m,r}:m,r\in\mathbb N\}.
$
Then $\mathcal G_0\subset L^1(\Omega,C(M))$ is countable, and
separates the measures in $\operatorname{Pr}_\P(\Omega\times M)$ by the metric defined in \eqref{metric_PrOM}. Moreover, since $\Theta(\omega,x)=(\theta\omega,F_\omega x)$ and
$x\mapsto F_\omega x$ is continuous, we have
$
g\circ\Theta\in L^1(\Omega,C(M))$
for every $g\in\mathcal G_0$.
For $q\geq 1,$ we define $R_q:\Omega\times M\to \Omega\times M^q$ by
\begin{equation*}
  R_q(\omega,x):=(\omega,x,F_\omega x,...,F_\omega^{q-1}x).
\end{equation*}
We choose a countable convergence-determining family $\mathcal{C}_q\subset C_b(\Omega\times M^q)$. Since $x\mapsto F_\omega^j x$ is continuous for any $j\geq 0$ and fixed $\omega\in\Omega$, we have  $\psi\circ R_q\in L^1(\Omega, C(M))$ for any $\psi\in \mathcal{C}_q$.
Now we define
\begin{equation}\label{eq def mathcalG}
\mathcal G:=
\mathcal G_0\cup\{g\circ\Theta:g\in\mathcal G_0\}  \cup\bigcup_{q\geq 1}\{\psi\circ R_q:\ \psi\in \mathcal{C}_q\}
\end{equation}
to be the admissible countable family of $L^1(\Omega,C(M))$ in Lemma \ref{lemma empirical port lemma}. By Lemma \ref{lemma empirical port lemma}, there exists a
$\theta$-invariant $\P$-full-measure set
$\Omega_{\mathcal G}$ corresponding to
$\mathcal G$ such that the conclusion of Lemma \ref{lemma empirical port lemma} holds. Intersect $\Omega_\mathcal{G}$ with the set of generic points of $\P$, and we denote the resulting set by $\Omega_{sep}$. From now on, we fix
$\omega\in\Omega_{sep}$, $\epsilon>0$, and prove Lemma \ref{lemma of Thm 2}.

   If $\Gamma_D(\omega,n,\epsilon)=\emptyset$ for all sufficiently large $n$, then the logarithmic limsup is $-\infty$, and the result is trivial. Now, we assume $\Gamma_D(\omega,n,\epsilon)\not=\emptyset$ for infinitely many $n\in\mathbb{N}$. Choose a sequence $n_i\to\infty$ such that
$\Gamma_D(\omega,n_i,\epsilon)\ne\varnothing$ and
  \begin{equation*}
     \limsup_{n\to+\infty}\frac{1}{n}\log \#\Gamma_D(\omega,n,\epsilon)=\lim_{i\to+\infty}\frac{1}{n_i}\log \# \Gamma_D(\omega,n_i,\epsilon).
  \end{equation*}
For every $n$ such that
$\Gamma_D(\omega,n,\epsilon)\ne\varnothing$, define the following three measures
  \begin{equation}
  \begin{split}
   \sigma_{n,\omega}:=   & \frac{1}{\#\Gamma_D(\omega,n,\epsilon)}\sum_{x\in \Gamma_D(\omega,n,\epsilon)}\delta_x\in \operatorname{Pr}(M), \\
     \sigma_{n}:=  & \delta_\omega\times \sigma_{n,\omega}=\frac{1}{\#\Gamma_D(\omega,n,\epsilon)}\sum_{x\in \Gamma_D(\omega,n,\epsilon)}\delta_{(\omega,x)}\in \operatorname{Pr}(\Omega\times M),\label{eq sigma n}
  \end{split}
  \end{equation}
and
  \begin{equation}\label{eq approximating measure}
    \mu_n=\frac{1}{\#\Gamma_D(\omega,n,\epsilon)}\sum_{x\in \Gamma_D(\omega,n,\epsilon)}\E_{n}(\omega,x)=\frac{1}{n}\sum_{j=0}^{n-1}(\Theta_*)^j(\sigma_{n})\in \operatorname{Pr}(\Omega\times M).
  \end{equation}
  \begin{claim}\label{claim in thm 2}
  After passing to a further subsequence, still denoted
by $\{n_i\}$, the measures $\mu_{n_i}$ converge to a measure $\rho$ with respect to $(\operatorname{Pr}(\Omega\times M),\tau)$.
  \end{claim}
  \begin{proof}[Proof of Claim \ref{claim in thm 2}]
Since $\omega$ is generic for $\P$, the sequence $\P_{n,\omega}$ converges
narrowly to $\P$. Its closure is therefore compact and hence uniformly
tight. For every $\epsilon>0$, choose a compact $\Omega_\epsilon\subset\Omega$
with $\P_{n,\omega}(\Omega_\epsilon)>1-\epsilon$ for every $n$.
Whenever $\mu_n$ is defined, its marginal is $\P_{n,\omega}$, so
$\mu_n(\Omega_\epsilon\times M)>1-\epsilon$. Since $M$ is compact,
Prohorov's theorem makes the sequence $\{\mu_{n_i}\}$ relatively compact.
Passing to a subsequence proves the claim.
\end{proof}

\begin{claim}\label{claim 5.2 in thm 2}
  The resulting measure $\rho$ belongs to $\overline{CH(D)}\cap I_\P(\Omega\times M)$.
\end{claim}
\begin{proof}[Proof of Claim \ref{claim 5.2 in thm 2}] It is clear that $\rho \in \overline{CH(D)}$.  Moreover, for any bounded continuous function $g:\Omega\to\mathbb{R}$, we have
  \begin{equation*}
  \begin{split}
  \int_{\Omega\times M} g(\omega^\prime)d\rho(\omega^\prime,x^\prime)=& \lim_{i\to\infty}\int_{\Omega\times M}g(\omega^\prime)d(\frac{1}{\#\Gamma_D(\omega,n_i,\epsilon)}\sum_{x\in \Gamma_D(\omega,n_i,\epsilon)}\E_{n_i}(\omega,x))(\omega^\prime,x^\prime)\\
  =& \lim_{i\to\infty}\int_{\Omega}g(\omega^\prime)d(\P_{n_i,\omega})(\omega^\prime)=\int_\Omega g(\omega)d\P(\omega),
  \end{split}
  \end{equation*}
since $\omega$ is a generic point of $\P$. Hence $\rho\in \operatorname{Pr}_\P(\Omega\times M)$. Next, we prove that $\rho$ is $\Theta$-invariant. We note that $\Theta$ on $\Omega\times M$ may not be continuous in the current setting.

For every
$g\in\mathcal G_0$, it follows from \eqref{eq approximating measure}
that
\[
\begin{split}
\int_{\Omega\times M}g\circ\Theta\,d\mu_n
-\int_{\Omega\times M}g\,d\mu_n
&=
\frac1n
\sum_{j=0}^{n-1}
\left(
\int g\,d(\Theta_*)^{j+1}\sigma_n
-\int g\,d(\Theta_*)^j\sigma_n
\right)\\
&=
\frac1n
\left(
\int g\,d(\Theta_*)^n\sigma_n
-\int g\,d\sigma_n
\right).
\end{split}
\]
Consequently,
\begin{equation}\label{eq asymptotic invariance mun}
\left|
\int g\circ\Theta\,d\mu_n-\int g\,d\mu_n
\right|
\leq\frac{2\|g\|_\infty}{n}.
\end{equation}

Recall that $\mu_{n_i}\to\rho$ with respect to
$(\operatorname{Pr}(\Omega\times M),\tau)$, that the marginal measure of
$\mu_{n_i}$ is
\[
\frac1{n_i}\sum_{j=0}^{n_i-1}\delta_{\theta^j\omega},
\]
and that $\rho\in \operatorname{Pr}_\P(\Omega\times M)$. Since
$\omega\in\Omega_{\mathcal G}$, Lemma
\ref{lemma empirical port lemma} gives
\begin{equation*}
  \lim_{i\to+\infty}\int g\,d\mu_{n_i}
=
\int g\,d\rho,\mbox{ and }\lim_{i\to+\infty}\int g\circ\Theta\,d\mu_{n_i}
=
\int g\circ\Theta\,d\rho,\ \forall g\in\mathcal G_0.
\end{equation*}
 Letting $i\to+\infty$ in
\eqref{eq asymptotic invariance mun}, we obtain
\begin{equation*}
  \int_{\Omega\times M}g\circ\Theta\,d\rho
=
\int_{\Omega\times M}g\,d\rho,
\qquad
\forall g\in\mathcal G_0.
\end{equation*}

Furthermore,
\[
(\pi_\Omega)_*(\Theta_*\rho)
=
\theta_*((\pi_\Omega)_*\rho)
=
\theta_*\P
=
\P,
\]
so both $\Theta_*\rho$ and $\rho$ belong to
$\operatorname{Pr}_\P(\Omega\times M)$. Since $\mathcal G_0$ separates the measures
in $\operatorname{Pr}_\P(\Omega\times M)$, it follows that
$
\Theta_*\rho=\rho.
$
Thus, we obtain $
\rho\in I_\P(\Omega\times M).$
Combining this with $\rho\in\overline{CH(D)}$, we conclude that
$
\rho\in\overline{CH(D)}\cap I_\P(\Omega\times M).
$ This completes the proof of Claim \ref{claim 5.2 in thm 2}.
\end{proof}
By Lemma \ref{lemma empirical port lemma} and \eqref{eq def mathcalG}, we also have that for $q\geq 1$,
\begin{equation*}
  \int_{\Omega\times M}\psi\circ R_q d\mu_{n_i}\to \int_{\Omega\times M}\psi\circ R_q d\rho,\ \forall \psi\in \mathcal{C}_q\mbox{ as }i\to+\infty.
\end{equation*}
As a consequence,
\begin{equation*}
  (R_q)_*\mu_{n_i}\to (R_q)_*\rho
\end{equation*}
in $\operatorname{Pr}(\Omega\times M^q)$ with respect to the narrow topology generated by $C_b(\Omega\times M^q)$ for all $q\geq 1$.

  Next, we prove $\lim\limits_{i\to+\infty}\frac{1}{n_i}\log \#\Gamma_D(\omega,n_i,\epsilon)\leq h_\rho(F)$ by following a classical argument of Misiurewicz; see, e.g., \cite[Proposition 2.2]{Kifervariationalprinciple}.
  We first choose a finite measurable partition $\mathcal{Q}=\{\Omega\times P_i\}_{i=1}^{l}$ of $\Omega\times M$ with $\max_{1\leq i\leq l}$diam$(P_i)< \epsilon$, such that
\begin{equation}\label{boundary 0}
  \rho(\partial\mathcal{Q})=\sum_{i=1}^{l}\int_\Omega\rho_\omega(\partial P_i)d\mathbb{P}(\omega)=0.
\end{equation}
Let $\mathcal P=\{P_i\}_{i=1}^l$ and
$\mathcal Q^{(q)}=\bigvee_{r=0}^{q-1}(\Theta^r)^{-1}\mathcal Q$.
Since $\sigma_n$ is supported on a single base fiber, and no atom of the
$n$-step fiber partition contains two points of $\Gamma_D(\omega,n,\epsilon)$,
\begin{equation}\label{eq h=log gamma}
 H_{\sigma_n}(\mathcal Q^{(n)}\mid\mathcal B_\Omega)
 =H_{\sigma_{n,\omega}}\left(\bigvee_{r=0}^{n-1}(F_\omega^r)^{-1}\mathcal P\right)
 =\log\#\Gamma_D(\omega,n,\epsilon).
\end{equation}
Fix $1\leq q<n$. For each $m\in\{0,\ldots,q-1\}$, write
$n=m+a(m)q+r_m$, where $0\leq r_m<q$. Subadditivity of conditional
entropy gives
\[
 H_{\sigma_n}(\mathcal Q^{(n)}\mid\mathcal B_\Omega)
 \leq\sum_{j=0}^{a(m)-1}
 H_{(\Theta^{m+jq})_*\sigma_n}(\mathcal Q^{(q)}\mid\mathcal B_\Omega)
 +2q\log l.
\]
Here invertibility of $\theta$ ensures that the base conditioning
$\sigma$-algebra is unchanged by pull-back. Summing over $m$ and using
concavity of conditional entropy in the underlying probability measure,
we obtain
\[
 \begin{aligned}
 q\log\#\Gamma_D(\omega,n,\epsilon)
 &\leq\sum_{j=0}^{n-1}
 H_{(\Theta^j)_*\sigma_n}(\mathcal Q^{(q)}\mid\mathcal B_\Omega)+2q^2\log l\\
 &\leq n H_{\mu_n}(\mathcal Q^{(q)}\mid\mathcal B_\Omega)+2q^2\log l.
 \end{aligned}
\]
The use of concavity does not require the measures being averaged to
have the same base marginal; see also \cite[Theorem~4.3]{Waltersergodic}
for the elementary entropy properties. Dividing by $n$ and passing to
the chosen subsequence gives
\begin{equation}\label{eq q leq H rho}
 q\lim_{i\to\infty}\frac1{n_i}\log\#\Gamma_D(\omega,n_i,\epsilon)
 \leq\limsup_{i\to\infty}
 H_{\mu_{n_i}}(\mathcal Q^{(q)}\mid\mathcal B_\Omega).
\end{equation}
We claim that
\begin{equation}\label{eq sep conditional entropy usc}
 \limsup_{i\to\infty}H_{\mu_{n_i}}(\mathcal Q^{(q)}\mid\mathcal B_\Omega)
 \leq H_\rho(\mathcal Q^{(q)}\mid\mathcal B_\Omega).
\end{equation}
Since $\Omega$ is a Polish space, we can choose an increasing sequence of finite Borel partitions $\mathcal A_m$
of $\Omega$ that generate $\mathcal F$ and whose atoms
have $\P$-null boundaries. Fix $q\ge1$ and $m\ge1$. Since $\rho$ is $\Theta$-invariant,
the first-coordinate marginal of $(R_q)_*\rho$ is $\P$,
and each of its fiber-coordinate marginals is
$(\pi_M)_*\rho$. Consequently, by \eqref{boundary 0},
every rectangle
\[
A\times P_{a_0}\times\cdots\times P_{a_{q-1}},
\qquad A\in\mathcal A_m,
\]
has boundary of zero $(R_q)_*\rho$-measure.

The narrow convergence
$(R_q)_*\mu_{n_i}\to(R_q)_*\rho$
and the ordinary Portmanteau theorem therefore imply
convergence of every atom probability of
$\mathcal Q^{(q)}\vee\pi_\Omega^{-1}\mathcal A_m$.
By continuity of finite-partition entropy,
\begin{equation*}
  \lim_{i\to\infty}H_{\mu_{n_i}}
  (\mathcal Q^{(q)}|\pi_\Omega^{-1}\mathcal A_m)= H_\rho
  (\mathcal Q^{(q)}|\pi_\Omega^{-1}\mathcal A_m)
\end{equation*}
It follows that
\begin{equation*}
\limsup_{i\rightarrow+\infty}H_{\mu_{n_i}}(\mathcal Q^{(q)}|\mathcal B_\Omega)
  \leq   \limsup_{i\rightarrow+\infty}H_{\mu_{n_i}}(\mathcal Q^{(q)}|\pi_\Omega^{-1}\mathcal{A}_m)
    =  H_\rho
  (\mathcal Q^{(q)}|\pi_\Omega^{-1}\mathcal A_m).
\end{equation*}
Letting $m\to\infty$ proves
\eqref{eq sep conditional entropy usc},
by convergence of conditional entropy along increasing
generating sigma-algebras.

Combining \eqref{eq q leq H rho} and \eqref{eq sep conditional entropy usc}, dividing by $q$, and then letting $q\rightarrow\infty$, we obtain
\begin{equation*}
 \limsup_{n\to+\infty}\frac{1}{n}\log \#\Gamma_D(\omega,n,\epsilon)=\lim_{i\to+\infty}\frac{1}{n_i}\log\#\Gamma_D(\omega,n_i,\epsilon)\leq h_\rho(F,\mathcal Q)\leq h_\rho(F).
\end{equation*}

Next, assume further that $\mu\mapsto h_\mu(F)$ is upper semicontinuous on $I_\P(\Omega\times M)$. We note that $\bar{D}\subset \operatorname{Pr}(\Omega\times M)$ may not be compact. By the proof of Claim \ref{claim in thm 2}, the set $\mathbb{A}:=\overline{\{\E_n(\omega,x):\ x\in M,\ n\in\mathbb{N}\}}$ is compact in $\operatorname{Pr}(\Omega\times M)$. Moreover,
\begin{equation*}
   \limsup_{n\to+\infty}\frac{1}{n}\log \#\Gamma_D(\omega,n,\epsilon)= \limsup_{n\to+\infty}\frac{1}{n}\log \#\Gamma_{D\cap \A}(\omega,n,\epsilon)
\end{equation*}
For each $\delta>0$, there exists a finite cover of $D\cap \A$ by open balls $\{B_{d_0}(\nu_i,\delta)\}$ of radius $\delta$ in $\operatorname{Pr}(\Omega\times M)$ with respect to the metric in Lemma \ref{lemma ws top}. The finite index set and the choice of centers depend on $\delta$. Each open ball $B_{d_0}(\nu_i,\delta)$ is convex in $\operatorname{Pr}(\Omega\times M)$ with respect to this metric. Restriction of a
separated set to the members of the finite cover gives
\[
\#\Gamma_D(\omega,n,\epsilon)\leq\sum_{i}
 \#\Gamma_{D\cap\mathbb A\cap B_{d_0}(\nu_i,\delta)}(\omega,n,\epsilon).
\]
Since $i$ belonging to a finite index set which is independent of $n$, there exists an $i^*$ such that
\begin{equation*}
  \limsup_{n\to +\infty}\frac{1}{n}\log \#\Gamma_{D\cap \A\cap B_{d_0}(\nu_{i^*},\delta)}(\omega,n,\epsilon)\geq \limsup_{n\to+\infty}\frac{1}{n}\log \#\Gamma_{D\cap \A}(\omega,n,\epsilon).
\end{equation*}
The preceding construction shows that there exists $\rho_\delta \in I_\P(\Omega\times M)\cap \overline{CH(D\cap\A\cap B_{d_0}(\nu_{i^*},\delta))}\subset I_\P(\Omega\times M)\cap \overline{B_{d_0}(\nu_{i^*},\delta)} $ satisfying
\begin{equation}\label{eq log gamma leq hrho}
   \limsup_{n\to+\infty}\frac{1}{n}\log \#\Gamma_{D}(\omega,n,\epsilon)\leq \limsup_{n\to +\infty}\frac{1}{n}\log \#\Gamma_{D\cap \A\cap B_{d_0}(\nu_{i^*},\delta)}(\omega,n,\epsilon)\leq h_{\rho_\delta}(F).
\end{equation} Choose any sequence $\delta_j\downarrow0$. By compactness of $I_\P(\Omega\times M)$ (Remark~\ref{remark compactness}), pass to a subsequence, relabeled by $j$, such that $\rho_{\delta_j}\to\rho^*\in I_\P(\Omega\times M)$ in $\tau$. Let $\nu_j$ be the center of the ball selected for $\delta_j$, and choose  $\eta_j\in D\cap\mathbb A\cap B_{d_0}(\nu_j,\delta_j).$ Since $\rho_{\delta_j}\in \overline{B_{d_0}(\nu_j,\delta_j)}$, we have $ d_0(\rho_{\delta_j},\eta_j)\le2\delta_j.$ Thus $\eta_j\to\rho^*$ and $\rho^*\in\overline D$. Hence $\rho^*\in I_\P(\Omega\times M)\cap \overline{D}$.

Lemma \ref{lemma lemma diffrent converges} converts the convergence $\rho_{\delta_j}\to \rho^*$ into narrow convergence with respect to $\operatorname{Pr}_\P(\Omega\times M,\tau_{rn})$. The upper semicontinuity of $\mu\mapsto h_\mu(F)$ on $I_\P(\Omega\times M)$ and \eqref{eq log gamma leq hrho} imply
\begin{equation*}
  \limsup_{n\to+\infty}\frac{1}{n}\log \#\Gamma_{D}(\omega,n,\epsilon)\leq h_{\rho^*}(F),
\end{equation*}
This completes the proof of Lemma \ref{lemma of Thm 2}.
\end{proof}
  Since $T$ is an upper energy function, the following holds on a $\P$-full-measure set:
  \begin{equation}\label{eq upper energy}
    \lim_{\epsilon \searrow 0}\limsup_{N\to +\infty}\sup_{x\in M}\left(\frac{1}{N}\log\nu_\omega(B_\omega(x,N,\epsilon))+\int Td(\E_N(\omega,x))\right)\leq 0,
  \end{equation}
  and there exists a sequence $\{\phi_l\}_{l\in\mathbb{N}}\subset L^1(\Omega,C(M))$ such that $\phi_l(\omega,x)\searrow (-T)(\omega,x)$ for $\P$-a.s. $\omega\in\Omega$ and all $x\in M$.

Choose a deterministic $B>0$ larger than
$\operatorname*{ess\,sup}_\omega\sup_xT(\omega,x)$. Take the invariant
full-measure core of the set where the bound $0\leq T\leq B$ and the
representatives $\phi_l\downarrow-T$ hold on the whole fiber. Intersect it
with the sets supplied by \eqref{eq upper energy}, Lemma~\ref{lemma empirical port lemma}
for $\mathcal G=\{\phi_l:l\geq1\}$, and Lemma~\ref{lemma of Thm 2}.
Call the resulting invariant full-measure set $\Omega_{\rm up}$.
For every $\omega\in\Omega_{\rm up}$, these properties hold at all
$\theta^j\omega$, $j\in\mathbb Z$. This set is independent of the closed
set $K$, the energy bands and the separation radius. Fix such an $\omega$.

\subsection{The case that $K\subset Pr(\Omega\times M)$ is a closed convex subset}
We first prove \eqref{eq thm 2} in the convex case.
  Given $\delta>0$, choose $0=a_0<a_1<\cdots<a_k=B$ with $a_j-a_{j-1}<\delta$. The choice $B>0$ also covers $T\equiv0$. For $j=1,\ldots,k$, we define
  $$K_j:=\{\mu\in K:\ a_{j-1}\leq \int_{\Omega\times M}Td\mu\le a_j\}.$$
   We choose a sufficiently small $\epsilon=\epsilon(\omega)>0$ and $n_{\epsilon}=n_{\epsilon}(\omega)\in\mathbb{N}$ so that for all $x\in M$
  \begin{equation}\label{eq leq delta}
\frac{1}{n}\log\nu_\omega(B_\omega(x,n,\epsilon))+\int Td(\E_n(\omega,x))\leq \delta,\ \forall n\geq n_{\epsilon}.
  \end{equation}
Let $\Gamma_{K_j}(\omega,n,\epsilon)$ be an $(\omega,n,\epsilon)$-separated set of $\{x\in M:\E_{n}(\omega,x)\in K_j\}$ with maximal cardinality. If $\Gamma_{K_j}(\omega,n,\epsilon)$ is empty for all
sufficiently large $n$, then
\[
\limsup_{n\to\infty}\frac1n
\log\nu_\omega\{x:\mathcal E_n(\omega,x)\in K_j\}
=-\infty,
\]
and no further argument is needed for this index. For every remaining index, we have
\begin{equation*}
  \{x\in M:\E_{n}(\omega,x)\in K_j\}\subset \bigcup_{x\in \Gamma_{K_j}(\omega,n,\epsilon)}B_\omega(x,n,\epsilon).
\end{equation*}
For $n\geq n_\epsilon(\omega)$ with $\Gamma_{K_j}(\omega,n,\epsilon)\ne\varnothing$, we have
\begin{equation}\label{eq with convex}
  \begin{split}
     &\frac{1}{n}\log\nu_\omega\{x\in M:\ \E_n(\omega,x)\in K_j\} \\
     \leq & \frac{1}{n}\left(\log\# \Gamma_{K_j}(\omega,n,\epsilon)+ \sup_{x\in \Gamma_{K_j}(\omega,n,\epsilon)}\log \nu_\omega(B_\omega(x,n,\epsilon))\right)\\
      \leq  & \frac{1}{n}\log\# \Gamma_{K_j}(\omega,n,\epsilon)+ \delta-\inf_{x\in \Gamma_{K_j}(\omega,n,\epsilon)}\int Td(\E_n(\omega,x)),\ \mbox{ by \eqref{eq leq delta}}\\
      \leq &\frac{1}{n}\log\# \Gamma_{K_j}(\omega,n,\epsilon)+ \delta-\inf_{\mu\in K_j}\int Td\mu\\
      \leq &\frac{1}{n}\log\# \Gamma_{K_j}(\omega,n,\epsilon)+ \delta-a_{j-1}.
  \end{split}
\end{equation}
When the separated set is empty, the deviation probability is zero, so the same upper bound is understood with logarithm $-\infty$ without using an infimum over an empty set. Since $K_j$ is convex, Lemma \ref{lemma of Thm 2} implies that there exists $\rho_j\in  \overline{K_j}\cap I_\P(\Omega\times M)$ such that
\begin{equation}\label{eq with convex 2}
  \begin{split}
      \limsup_{n\to+\infty}  \frac{1}{n}\log\nu_\omega\{x\in M:\ \E_n(\omega,x)\in K_j\}
      \leq  &h_{\rho_j}(F)+\delta-a_{j-1}\\
      \leq & (h_{\rho_j}(F)-a_{j})+2\delta.
  \end{split}
\end{equation}
Recall that $\rho_j$ is the limit of a subsequence of the measures in \eqref{eq approximating measure} with respect to $(\operatorname{Pr}(\Omega\times M),\tau)$. This subsequence and $\rho_j$ satisfy the hypotheses of Lemma \ref{lemma empirical port lemma} for $\mathcal{G}=\{\phi_l\}_{l\in\mathbb{N}}$.
\begin{claim}\label{cliam 2 in thm2}
  \begin{equation}\label{eq int T aj}
    \int Td\rho_j\leq a_j
  \end{equation}
\end{claim}
\begin{proof}[Proof of Claim \ref{cliam 2 in thm2}]
We retain the notation from the proof of Lemma \ref{lemma of Thm 2}. Since $j$ is fixed, we omit $j$ in this proof.
  By Lemma \ref{lemma empirical port lemma}, we obtain
  \begin{equation}\label{eq nontrivial}
    \lim_{i\to+\infty}\int \phi_l d\mu_{n_i}=\int \phi_l d\rho.
\end{equation}
Recall that $\phi_l(\omega,x)\searrow (-T)(\omega,x)$ for $\P$-a.s. $\omega\in\Omega$ and all $x\in M$. Then, for each $l$,
   \begin{equation*}
     \limsup_{i\to +\infty}\int(-T)d\mu_{n_i}\leq \lim_{i\to +\infty}\int\phi_l d\mu_{n_i}\overset{\eqref{eq nontrivial}}=\int \phi_l d\rho.
   \end{equation*}
   Note that $-a_k\leq -T\leq \phi_l\leq \phi_1$ for $\P$-a.s. $\omega$ and all $x\in M$. Then $|\phi_l|\leq a_k+|\phi_1|$.
   By the dominated convergence theorem, we have
$
    \lim_{l\to+\infty}\int \phi_l d\rho=\int (-T)d\rho.
$ Therefore,
$$
  \int Td\rho\leq \liminf\int Td\mu_{n_i}.$$ Recall that $\mu_n$ is constructed from points in $\{x\in M:\ \E_n(\omega,x)\in K_j\}$. We therefore obtain
  \begin{equation*}
     \int Td\rho\leq \liminf\int Td\mu_{n_i}\leq a_j.
  \end{equation*}
This completes the proof of the claim.
\end{proof}
Note that $K$ is closed, so $\rho_j\in I_\P(\Omega\times M)\cap \overline{K_j}\subset K\cap I_\P(\Omega\times M)$.  Hence, by \eqref{eq with convex 2} and \eqref{eq int T aj},
\begin{equation*}
\begin{split}
   \limsup_{n\to+\infty}  \frac{1}{n}\log\nu_\omega\{x\in M:\ \E_n(\omega,x)\in K_j\}\leq    & (h_{\rho_j}(F)-\int Td\rho_j)+2\delta\\
     \leq & \sup_{\rho\in K\cap I_\P(\Omega\times M)} (h_\rho(F)-\int Td\rho)+2\delta.
\end{split}
\end{equation*}
To complete the proof, first note that
  \begin{equation*}
    \limsup_{n\to+\infty}\frac{1}{n}\log\nu_\omega\{x\in M:\ \E_n(\omega,x)\in K\}=\max_{1\leq j\leq k}\limsup_{n\to+\infty}\frac{1}{n}\log\nu_\omega\{x\in M:\ \E_n(\omega,x)\in K_j\},
  \end{equation*}
Letting $\delta\downarrow0$ completes the proof of the convex case.

\subsection{The general closed-set case}
Now assume further that $\mu\mapsto h_\mu(F)$ is upper semicontinuous on $I_\P(\Omega\times M)$ and that $K\subset \operatorname{Pr}(\Omega\times M)$ is closed but not necessarily convex. The proof up to \eqref{eq with convex} requires no change; we replace the remainder by the following argument.

By the second part of Lemma \ref{lemma of Thm 2}, there exists $\rho_j^*\in  \overline{K_j}\cap I_\P(\Omega\times M)$ such that
\begin{equation}\label{eq with convex 2-2}
  \begin{split}
      \limsup_{n\to+\infty}  \frac{1}{n}\log\nu_\omega\{x\in M:\ \E_n(\omega,x)\in K_j\}
      \leq  &h_{\rho_j^*}(F)+\delta-a_{j-1}\\
      \leq & (h_{\rho_j^*}(F)-a_{j})+2\delta.
  \end{split}
\end{equation}
  Recall that, in the second part of the proof of Lemma \ref{lemma of Thm 2}, $\rho_j^*$ is the limit of measures $\rho_{\delta_n}\in I_\P(\Omega\times M)\cap \overline{CH(K_j\cap \A\cap B(\mu_{i^*},\delta_n))}$ as $\delta_n \to 0$ with respect to $(\operatorname{Pr}(\Omega\times M),\tau)$.
\begin{claim}\label{cliam 22 in thm2}
  $\int Td\rho_j^*\leq a_j$.
\end{claim}
\begin{proof}[Proof of Claim \ref{cliam 22 in thm2}]
We retain the notation from the proof of Lemma \ref{lemma of Thm 2}.
 By Lemma \ref{lemma lemma diffrent converges}, we obtain
  \begin{equation}\label{eq nontrivial 2}
    \lim_{n\to+\infty}\int \phi_l d\rho_{\delta_n}=\int \phi_l d\rho_j^*.
  \end{equation}
Recall that $\phi_l(\omega,x)\searrow (-T)(\omega,x)$ for $\P$-a.s. $\omega\in\Omega$ and all $x\in M$. Then, for each $l$,
   \begin{equation*}
     \limsup_{n\to +\infty}\int(-T)d\rho_{\delta_n}\leq \lim_{n\to +\infty}\int\phi_l d\rho_{\delta_n}\overset{\eqref{eq nontrivial 2}}=\int \phi_l d\rho_j^*.
   \end{equation*} By the dominated convergence theorem, we have
$
    \lim_{l\to+\infty}\int \phi_l d\rho_j^*=\int (-T)d\rho_j^*.
$ Therefore,
$$
  \int Td\rho_j^*\leq \liminf\int Td\rho_{\delta_n}.$$ By the proof of Claim \ref{cliam 2 in thm2}, we have
  \begin{equation*}
    \int Td\rho_{\delta_n}\leq a_j.
  \end{equation*}  We therefore obtain
  \begin{equation*}
     \int Td\rho_j^*\leq \liminf\int Td\rho_{\delta_n}\leq a_j.
  \end{equation*}
\end{proof}
Taking the maximum over the finitely many energy bands and letting $\delta\downarrow0$ completes the proof, as in the convex case.

\section{Proof of the large deviation lower bound}\label{sec:lower-bound}
In this section, we prove Theorem \ref{thm main 1}. To keep the main argument transparent and avoid interrupting its flow, we defer the proofs of all auxiliary lemmas used in this section to Section \ref{sec proof of lemma}. We first fix any open subset $G\subset \operatorname{Pr}(\Omega\times M)$.  Assume the reference measure $\nu \in \operatorname{Pr}_\P(\Omega\times M)$ admits a lower energy function, and let $T:\Omega\times M\to \mathbb{R}_{\geq 0}$ be any lower energy function for $\nu$. Lemma~\ref{lemma lower energy finite entropy} implies that $h_{\mathrm{top}}(F)=\sup\{h_\rho(F):\ \rho\in I_\P(\Omega\times M)\}<+\infty$. Therefore, we can apply Lemma \ref{lemma measure approxi} and Lemma \ref{lemma katok entropy} in the following proof.

\subsection{Invariant measures with positive fiber entropy}
 We show that, for $\P$-a.s. $\omega\in\Omega$,
  \begin{equation}\label{eq liminf ineq}
    \liminf_{n\to +\infty}\frac{1}{n}\log \nu_\omega\{x\in M:\ \mathcal{E}_n(\omega,x)\in G\}\geq h_\rho(F)-\int_{\Omega\times M} Td\rho,
  \end{equation}
where $T$ is a lower energy function for $\nu\in \operatorname{Pr}_\P(\Omega\times M).$

   Fix any $\rho \in G\cap I_\P(\Omega\times M)$ and $\eta>0$. Since $G$ is open with respect to $(\operatorname{Pr}(\Omega\times M),\tau)$, Lemma \ref{lemma ws top} gives $\{f_1,\ldots,f_{l-1}\}\subset BL_1(\Omega\times M)$ and $\delta>0$ such that
   \begin{equation*}
     \{\mu\in \operatorname{Pr}(\Omega\times M):|\int f_jd\mu-\int f_j d\rho|<6\delta,\ 1\leq j\leq l-1\}\subset G.
   \end{equation*}
   We remark that although $G\cap \operatorname{Pr}_\P(\Omega\times M)$ is open relative to $\operatorname{Pr}_\P(\Omega\times M)$, it may not be open in $\operatorname{Pr}(\Omega\times M)$. Moreover, this relative openness cannot by itself be applied to $\E_n(\omega,x)$, since the marginal measure $(\pi_\Omega)_*\E_n(\omega,x)=\frac{1}{n}\sum_{j=0}^{n-1}\delta_{\theta^j\omega}$ is generally not equal to $\P$. The next lemma supplies the additional energy control needed in the proof.
   \begin{lemma}\label{lemma energy control}
     Recall that $T:\Omega\times M\to \mathbb{R}_{\geq 0}$ is an upper normal integrand and $\rho\in \operatorname{Pr}_\P(\Omega\times M)$. For any $\eta>0$, there exist a constant $C\geq 1$, a function $q\in L^1(\Omega,C(M))$ with $q:\Omega\times M\to [0,1]$ and $[q(\omega,\cdot)]_L\leq 1$ for all $\omega\in\Omega$, and a nonnegative function $r\in L^1(\P)$ such that
\begin{equation}\label{eq T leq CQ+r}
  T(\omega,x)\leq Cq(\omega,x)+r(\omega),\ \mbox{ for $\P$-a.s. $\omega\in\Omega$ and every $x\in M$},
\end{equation}
and
\begin{equation}\label{eq C int q+r leq T}
  C\int_{\Omega\times M} q(\omega,x) d\rho(\omega,x)+\int_\Omega r(\omega) d\P(\omega)\leq \int_{\Omega\times M}T(\omega,x)d\rho(\omega,x)+\frac{\eta}{3}.
\end{equation}
   \end{lemma}
We rename the function $q$ in Lemma \ref{lemma energy control} as $f_l$ and reduce $\delta$, if necessary, so that
\begin{equation}\label{eq delta eta}
  6C\delta<\frac{\eta}{3}.
\end{equation}
Set
\begin{equation*}
  O_\rho:= \{\mu\in \operatorname{Pr}(\Omega\times M):|\int f_jd\mu-\int f_j d\rho|<6\delta,\ 1\leq j\leq l\}.
\end{equation*}
By the preceding construction, $O_\rho\subset G$, although $O_\rho$ need not be open.

In the next lemma, we construct points whose fiber empirical measure is close to $\rho$. We emphasize that $\{f_1,..,f_{l-1}\}\subset BL_1(\Omega\times M)$ and $f_l\in L^1(\Omega,C(M))$ satisfies $0\leq f_l\leq1$ and $[f_l(\omega,\cdot)]_L\leq 1$ for all $\omega\in\Omega$. Fix $\{f_1,...,f_l\}$ from now on.
  \begin{lemma}\label{Lemma lower estimates}
Assume that the skew-product system has the $\P$-a.s. fiber specification property, and let $\rho\in I_\P(\Omega\times M)$ with $0<h_\rho(F)<\infty$. For any $\delta>0$, $0<\gamma<\min\{1,h_\rho(F)/4\}$, and any sufficiently small $\xi>0$ such that
\begin{equation}\label{xi estimation}
			\xi\leq \min\left\{ \frac{\delta}{16}, \ \frac{1}{4}\left(1-\frac{h_{\rho}(F)-4\gamma}{h_{\rho}(F)-3\gamma}\right)\right\},
		\end{equation}
there exist $\epsilon_0>0$, a $\mathbb{P}$-full-measure set $\Omega^\prime$, and a measurable function $\bar{n}:\Omega^\prime\to\mathbb{N}$ such that the following holds. For any $\omega\in\Omega^\prime$ and any $N>\bar{n}(\omega)$, there exists $\hat{n}(\omega)=\hat{n}_N(\omega)\in\mathbb{N}$, depending on $N$, such that
\begin{enumerate}
  \item the following deviation set is not empty:
      \begin{equation*}
  P(\rho,4\delta,\hat{n}(\omega),\omega):=\left\{x\in M:\left|\frac{1}{\hat{n}(\omega)}\sum_{i=0}^{\hat{n}(\omega)-1}f_j(\Theta^i(\omega,x))-\int_{\Omega\times M} f_jd\rho\right|<4\delta,\ j=1,\ldots,l\right\}\not=\emptyset;
\end{equation*}
  \item let $M(\omega,4\delta,\frac{\epsilon_0}{2},\hat{n}(\omega))$ be an $(\omega,\hat{n}(\omega),\frac{\epsilon_0}{2})$-separated subset of $P(\rho,4\delta,\hat{n}(\omega),\omega)$ of maximal cardinality. Then we have the estimate
\begin{equation}\label{eq estimate 1/n log M}
 \frac{1}{\hat{n}(\omega)}\log \# M(\omega,4\delta,\frac{\epsilon_0}{2},\hat{n}(\omega))\geq h_\rho(F)-4\gamma;
\end{equation}
  \item $N\cdot (1-2\xi)\leq \hat{n}(\omega)\leq N$.
\end{enumerate}
    \end{lemma}
   By Birkhoff's ergodic theorem, we may assume that, for $\omega\in\Omega^\prime$ and $N>\bar{n}(\omega)$, one has
   \begin{equation}\label{eq r leq r+eta}
     \frac{1}{N}\sum_{j=0}^{N-1}r(\theta^j\omega)\leq \int_\Omega r(\omega')d\P(\omega')+\frac{\eta}{3},
   \end{equation}
where the nonnegative function $r\in L^1(\P)$ was constructed in Lemma \ref{lemma energy control}. This follows by intersecting $\Omega^\prime$ with another $\P$-full-measure set and increasing $\bar{n}(\omega)$ if necessary. We further intersect $\Omega'$ with the full-measure set for the
lower energy condition and with the invariant core of the set where
\eqref{eq T leq CQ+r} holds on the whole fiber. Hence
\[
 T(\theta^j\omega,x)\leq C f_l(\theta^j\omega,x)+r(\theta^j\omega)
 \qquad(j\geq0,\ x\in M)
\]
for every base point under consideration.

   From now on, we fix any $\omega\in\Omega^\prime$, pick any $N>\bar{n}(\omega)$, and let $\hat{n}(\omega)$ be chosen as in Lemma \ref{Lemma lower estimates}.
For any $x\in M(\omega,4\delta,\frac{\epsilon_0}{2},\hat{n}(\omega))\subset P(\rho,4\delta,\hat{n}(\omega),\omega)$, since $N(1-2\xi)\leq \hat{n}(\omega)\leq N$, we have
\begin{equation*}
 \begin{split}
   & |\frac{1}{N}\sum_{i=0}^{N-1}f_j(\Theta^i(\omega,x))- \frac{1}{\hat{n}(\omega)}\sum_{i=0}^{\hat{n}(\omega)-1}f_j(\Theta^i(\omega,x))|  \\
      \leq  & |\frac{1}{N}-\frac{1}{\hat{n}(\omega)}|\cdot \hat{n}(\omega)\cdot \|f_j\|_{\infty}+\frac{1}{N}|\sum_{i=\hat{n}(\omega)}^{N-1}f_j(\Theta^i(\omega,x))|\\
      \leq &2\frac{N-\hat{n}(\omega)}{N}\|f_j\|_{\infty}\leq 4\xi \|f_j\|_{\infty}\leq \frac{\delta}{2}, \mbox{ by \eqref{xi estimation}}.
  \end{split}
\end{equation*}
By the triangle inequality,
\begin{equation*}
  M(\omega,4\delta,\frac{\epsilon_0}{2},\hat{n}(\omega))\subset \left\{x\in M:\left|\frac{1}{N}\sum_{i=0}^{N-1}f_j(\Theta^i(\omega,x))-\int_{\Omega\times M} f_jd\rho\right|<\frac{9}{2}\delta,\ j=1,\ldots,l\right\}
\end{equation*}
Since $[f_j(\omega,\cdot)]_L\leq 1$, we pick $\varrho=\frac{\delta}{2}$ such that
\begin{equation*}
  d_M(x,y)\leq \varrho \mbox{ implies that }|f_j(\omega,x)-f_j(\omega,y)|\leq \frac{\delta}{2}.
\end{equation*}
For any $0<\epsilon<\min\{\frac{\epsilon_0}{4},\varrho\}$, the set $M(\omega,4\delta,\frac{\epsilon_0}{2},\hat{n}(\omega))$ is $(\omega,\hat{n}(\omega),2\epsilon)$-separated. Moreover, since $N\geq\hat n(\omega)$ by Lemma~\ref{Lemma lower estimates}, this set is also $(\omega,N,2\epsilon)$-separated. Hence, for any two distinct points $x,y\in M(\omega,4\delta,\frac{\epsilon_0}{2},\hat{n}(\omega))$, we have
\begin{equation*}
  B_\omega(x,N,\epsilon)\cap B_\omega(y,N,\epsilon)=\emptyset,
\end{equation*}
and
\begin{equation}\label{eq balls sub 5delta}
  \bigcup_{x\in M(\omega,4\delta,\frac{\epsilon_0}{2},\hat{n}(\omega))}B_\omega(x,N,\epsilon)\subset \left\{x\in M:\left|\frac{1}{N}\sum_{i=0}^{N-1}f_j(\Theta^i(\omega,x))-\int_{\Omega\times M} f_jd\rho\right|<5\delta,\ j=1,\ldots,l\right\}.
\end{equation}
   The above relation implies that $\E_N(\omega,y) \in O_\rho$ for any $y\in \bigcup_{x\in M(\omega,4\delta,\frac{\epsilon_0}{2},\hat{n}(\omega))}B_\omega(x,N,\epsilon)$, and moreover
   \begin{equation}\label{eq int T E leq int T rho}
   \begin{split}
    \int_{\Omega\times M} T d\E_N(\omega,y)\leq   &C \int_{\Omega\times M} f_l d\E_N(\omega,y)+\frac{1}{N}\sum_{j=0}^{N-1}r(\theta^j\omega),\ \mbox{ by \eqref{eq T leq CQ+r}}\\
        \leq & C\int_{\Omega\times M} f_l d\rho+6C\delta+\int_\Omega r(\omega)d\P(\omega)+\frac{\eta}{3},\mbox{ by \eqref{eq r leq r+eta}}\\
        \leq &\int_{\Omega\times M} T d\rho+\eta, \mbox{ by \eqref{eq C int q+r leq T} and \eqref{eq delta eta}}.
   \end{split}
   \end{equation}

 Since $O_\rho \subset G$, we obtain
\begin{equation*}
  \begin{split}
     &\log\nu_\omega\{x:\ \E_{N}(\omega,x)\in G\} \geq  \log \nu_\omega\{x:\ \E_{N}(\omega,x)\in O_\rho\}
       \geq  \log \nu_\omega\left( \bigcup_{x\in M(\omega,4\delta,\frac{\epsilon_0}{2},\hat{n}(\omega))}B_\omega(x,N,\epsilon)\right)\\
       \geq &\hat{n}(\omega)(h_\rho(F)-4\gamma)+\inf \{\log \nu_\omega(B_\omega(x,N,\epsilon)):\ x\in M(\omega,4\delta,\frac{\epsilon_0}{2},\hat{n}(\omega))\},\mbox{ by \eqref{eq estimate 1/n log M}}
  \end{split}
  \end{equation*}
  Hence, we have
  \begin{equation*}
  \small
    \begin{split}
     &\frac{1}{N}\log\nu_\omega\{x:\ \E_{N}(\omega,x)\in G\} \\
      \geq &(1-2\xi)(h_\rho(F)-4\gamma)+\inf_{ x\in M(\omega,4\delta,\frac{\epsilon_0}{2},\hat{n}(\omega))}\left(\frac{1}{N}\log \nu_\omega(B_\omega(x,N,\epsilon))\right)\\
       \geq &
       (1-2\xi)(h_\rho(F)-4\gamma)+\inf_{ x\in M(\omega,4\delta,\frac{\epsilon_0}{2},\hat{n}(\omega))}\left(\frac{1}{N}\log \nu_\omega(B_\omega(x,N,\epsilon))+\int T d(\E_N(\omega,x))-\int Td\rho-\eta\right),
  \end{split}
\end{equation*}
where the last inequality follows from \eqref{eq int T E leq int T rho}.
By the definition of a lower energy function (Definition \ref{definition lower energy}), we have
\begin{equation*}
  \liminf_{N\to +\infty}\frac{1}{N}\log\nu_\omega\{x:\ \E_{N}(\omega,x)\in G\}\geq (1-2\xi)(h_\rho(F)-4\gamma)-\int Td\rho-\eta.
\end{equation*}
We complete the proof by first letting $\xi\to0$ along a discrete sequence, then letting $\delta$, $\eta$, and $\gamma$ tend to $0$ along discrete sequences, and taking the intersection of the resulting $\P$-full-measure sets.

\subsection{Invariant measures with zero fiber entropy}
 The preceding argument up to Lemma \ref{Lemma lower estimates} remains valid. Analogously to Lemma \ref{Lemma lower estimates}, we obtain the following lemma, whose proof is given in Section \ref{sec proof of lemma}.
\begin{lemma}\label{Lemma lower estimates 2}
Assume that the skew-product system has the $\P$-a.s. fiber specification property, and let $\rho\in I_\P(\Omega\times M)$ with $h_\rho(F)=0$. Fix the same bounded random Lipschitz test functions as in Lemma~\ref{Lemma lower estimates}. For any $\delta>0$ and any $0<\xi\leq\min\{\frac{1}{4},\frac{\delta}{16}\}$,
there exist $\epsilon_0>0$, a $\mathbb{P}$-full-measure set $\Omega^\prime\subset\Omega$, and a measurable function $\bar{n}:\Omega^\prime\to\mathbb{N}$ such that the following holds. For any $\omega\in\Omega^\prime$ and any $N>\bar{n}(\omega)$, there exists $\hat{n}(\omega)\in\mathbb{N}$ such that
\begin{enumerate}
  \item the deviation set $P(\rho,4\delta,\hat{n}(\omega),\omega)\not=\emptyset$;
  \item $N\cdot (1-2\xi)\leq \hat{n}(\omega)\leq N$.
\end{enumerate}Consequently, if $M(\omega,4\delta,\frac{\epsilon_0}{2},\hat{n}(\omega))$ denotes a $(\omega,\hat{n}(\omega),\frac{\epsilon_0}{2})$-separated set of maximal cardinality in $P(\rho,4\delta,\hat{n}(\omega),\omega)$, then $\# M(\omega,4\delta,\frac{\epsilon_0}{2},\hat{n}(\omega))\geq 1$ and hence
\begin{equation*}
 \frac{1}{\hat{n}(\omega)}\log \# M(\omega,4\delta,\frac{\epsilon_0}{2},\hat{n}(\omega))\geq 0.
\end{equation*}
    \end{lemma}
It remains only to replace $h_\rho(F)-4\gamma$ by $0$ in the proof of the preceding case.

We now complete the lower bound for the fixed open set $G$. If
$G\cap I_\P(\Omega\times M)=\varnothing$, the bound is automatic.
Otherwise set
\[
 S=\sup_{\rho\in G\cap I_\P(\Omega\times M)}
       \left\{h_\rho(F)-\int T\,d\rho\right\}.
\]
The entropy bound and the integrable upper envelope of $T$ imply that $S$
is a finite real number. Choose $\rho_m\in G\cap I_\P(\Omega\times M)$
with $h_{\rho_m}(F)-\int T\,d\rho_m>S-1/m$.
The positive-entropy and zero-entropy arguments above both apply to these measures.
Taking the countable intersection of their full-measure sets and letting
$m\to\infty$ gives the required lower bound for $G$.

\subsection{A common full-measure set for all open sets}
Since $\Omega\times M$ is Polish, $\operatorname{Pr}(\Omega\times M)$ with its narrow topology has a countable base $\{U_j\}_{j\geq1}$. For every $j$, the established lower bound gives a full-measure set $\Omega_j$ on which
\begin{equation*}
    \liminf_{n\to +\infty}\frac{1}{n}\log \nu_\omega\{x\in M:\ \mathcal{E}_n(\omega,x)\in U_j\}\geq -\inf_{\rho\in U_j} I_T(\rho).
  \end{equation*}
Set $\Omega_0=\bigcap_{j\geq1}\Omega_j$. Then $\P(\Omega_0)=1$. Fix any $\omega\in\Omega_0$ and an arbitrary open set $G\subset \operatorname{Pr}(\Omega\times M)$, and $\rho\in G$ with $I_T(\rho)<+\infty$. We choose $j$ such that $\rho\in U_j\subset G$. Then by monotonicity, we have
  \begin{equation*}
  \begin{split}
    \liminf_{n\to +\infty}\frac{1}{n}\log \nu_\omega\{x\in M:\ \mathcal{E}_n(\omega,x)\in G\}\geq  & \liminf_{n\to +\infty}\frac{1}{n}\log \nu_\omega\{x\in M:\ \mathcal{E}_n(\omega,x)\in U_j\} \\
      \geq  & -\inf_{\eta\in U_j}I_T(\eta)\geq -I_T(\rho).
  \end{split}
  \end{equation*}
Taking the supremum over all $\rho\in G$ proves
  \begin{equation*}
     \liminf_{n\to +\infty}\frac{1}{n}\log \nu_\omega\{x\in M:\ \mathcal{E}_n(\omega,x)\in G\}\geq-\inf_{\rho\in G}I_T(\rho).
  \end{equation*}
If the infimum is $+\infty$, the inequality is automatic. Thus the lower bound holds for every open $G\subset \operatorname{Pr}(\Omega\times M)$ on the common set $\Omega_0$.

  \section{Proofs of the auxiliary lemmas for the lower bound}\label{sec proof of lemma}
  \begin{proof}[Proof of Lemma \ref{lemma energy control}]
     Recall that $T$ is an upper normal integrand and that $\rho\in \operatorname{Pr}_\P(\Omega\times M)$. There exists a sequence of nonnegative functions $\phi_k\in L^1(\Omega,C(M))$ with $\phi_k(\omega,x)\searrow T(\omega,x)$ for $\P$-a.s. $\omega\in\Omega$ and all $x\in M$. For $\omega$ in a $\P$-null set and $x\in M$, we set $\phi_k(\omega,x)=0$. By the dominated
convergence theorem, we choose $k$ sufficiently large that
\begin{equation}\label{eq phik leq T}
   \int_{\Omega\times M}\phi_k(\omega,x)d\rho\leq\int_{\Omega\times M} T(\omega,x)d\rho+\frac{\eta}{9}.
\end{equation}
Define $b(\omega)=\sup\limits_{x\in M}\phi_k(\omega,x)$. Then $b\in L^1(\P)$ by the definition of $L^1(\Omega,C(M))$.
For $m\geq1$, we define
$$
 \phi_k^{(m)}(\omega,x)
 :=
 \max_{y\in M}\{\phi_k(\omega,y)-m d_M(x,y)\}.
$$
 The joint measurability of $\phi_k^{(m)}$ follows by taking the supremum over a countable dense subset of $M$. Moreover, a direct calculation shows that $[\phi_k^{(m)}(\omega,\cdot)]_L\leq m$, and
$$
 0\leq \phi_k\leq\phi_k^{(m)}\leq b,
\
 \phi_k^{(m)}\searrow\phi_k, \mbox{ as }m\to +\infty,\mbox{ for $\P$-a.s. $\omega\in\Omega$ and all $x\in M$}
$$
Hence, for some $m$ sufficiently large,
\begin{equation}\label{eq phikm leq phik}
   \int_{\Omega\times M}\phi_k^{(m)}(\omega,x)d\rho(\omega,x)
 \leq\int_{\Omega\times M}\phi_kd\rho(\omega,x)+\frac{\eta}{9}.
\end{equation}
We can choose a constant $R>0$ such that
$$
 \int_\Omega(b(\omega)-R)^+d\P(\omega)<\frac{\eta}{9},
$$
and put
\[
 C=\max\{1,m,R\},\qquad
 q(\omega,x)=\frac{\min\{\phi_k^{(m)}(\omega,x),R\}}{C},
 \qquad
 r(\omega)=(b(\omega)-R)^+.
\]
Then $q\in L^1(\Omega,C(M))$, takes values in $[0,1]$, and satisfies $[q(\omega,\cdot)]_L\leq \frac{m}{C}\leq 1$ for all $\omega\in\Omega$ (after setting $q(\omega,\cdot)\equiv 0$ on a $\P$-null set if necessary). Moreover,
\[
 T\leq\phi_k\leq\phi_k^{(m)}
 \leq \min\{\phi_k^{(m)},R\}+(b-R)^+
 =Cq+r,\mbox{ for $\P$-a.s. $\omega\in\Omega$ and every $x\in M$}
\]
We also have
\[
\begin{split}
 C\int_{\Omega\times M} q(\omega,x)d\rho(\omega,x)+\int_\Omega r(\omega)d\P(\omega)
 &\leq \int_{\Omega\times M}\phi_k^{(m)}(\omega,x)d\rho(\omega,x)+\int r(\omega)d\P(\omega)\\
 &\leq \int_{\Omega\times M} T(\omega,x)d\rho(\omega,x)+\frac{\eta}{3},\mbox{ by \eqref{eq phik leq T} and \eqref{eq phikm leq phik}}.
\end{split}
\]
  \end{proof}

\begin{proof}[Proof of Lemma \ref{Lemma lower estimates}]
We recall that the functions $f_1,\ldots,f_l$ satisfy
$\|f_j\|_\infty\leq1$ and $[f_j(\omega,\cdot)]_L\leq1$
on every relevant fiber. Fix $\delta$, $\gamma$, and $\xi$ as in
the statement.

\emph{Step 1. We first construct the full-measure set $\Omega'$.}

By Lemma \ref{lemma measure approxi}, we can approximate $\rho$ by
a $\Theta$-invariant measure $\nu=\sum_{i=1}^k\lambda_i\nu_i$,
where $\lambda_i>0$, $\sum_{i=1}^k\lambda_i=1$, and
$\nu_i\in I_\P^e(\Omega\times M)$ has finite fiber entropy, such that
\[
 h_\nu(F)\geq h_\rho(F)-\gamma,\qquad
 \left|\int f_j\,d\rho-\int f_j\,d\nu\right|<\delta
 \quad(j=1,\ldots,l).
\]
Since each $\nu_i$ is ergodic, Birkhoff's ergodic theorem and
measure disintegration imply that the measurable functions
\[
 g_N^i(\omega):=(\nu_i)_\omega\left(
 \left\{x\in M:
 \left|\frac1r\sum_{p=0}^{r-1}f_j(\Theta^p(\omega,x))
              -\int f_j\,d\nu_i\right|<\delta,\
 j=1,\ldots,l,\ \text{for every }r\geq N\right\}\right)
\]
increase to $1$ as $N\to\infty$, for $\P$-a.s. $\omega$.
By Egorov's theorem, there exist a measurable set $\Omega_1$
with $\P(\Omega_1)>1-\xi/4$ and integers $N_{1,i}$ such that
\[
 g_N^i(\omega)\geq1-\gamma
 \quad\text{for every }\omega\in\Omega_1,\ N\geq N_{1,i},
 \text{ and }i=1,\ldots,k.
\]

By Lemma \ref{lemma katok entropy}, for $\P$-a.s. $\omega$,
\[
 \lim_{\epsilon\downarrow0}\liminf_{n\to\infty}
 \frac1n\log S_{\nu_i}(\omega,n,\epsilon,\gamma)
 =h_{\nu_i}(F),\qquad i=1,\ldots,k.
\]
Here $S_{\nu_i}(\omega,n,\epsilon,\gamma)$ is the smallest
cardinality of an $(\omega,n,\epsilon)$-spanning set covering
a set of $(\nu_i)_\omega$-measure at least $1-\gamma$.
Decreasing $\epsilon$ along a sequence, and using monotonicity,
we can choose $\epsilon_0>0$ and a measurable set
$\Omega_2\subset\Omega_1$ with $\P(\Omega_2)>1-\xi/2$ such that
\[
 \liminf_{n\to\infty}\frac1n
 \log S_{\nu_i}(\omega,n,\epsilon_0,\gamma)
 \geq h_{\nu_i}(F)-\gamma/2
 \quad(\omega\in\Omega_2,\ i=1,\ldots,k).
\]
For each $i$, the measurable sets
\[
 A_N^i:=\bigcap_{n\geq N}
 \left\{\omega:
 S_{\nu_i}(\omega,n,\epsilon_0,\gamma)
 \geq e^{n(h_{\nu_i}(F)-\gamma)}\right\}
\]
increase with $N$, and their union contains $\Omega_2$.
Thus we can choose integers $N_{2,i}$ such that
$\Omega_3:=\Omega_2\cap\bigcap_{i=1}^k A_{N_{2,i}}^i$
satisfies $\P(\Omega_3)>1-3\xi/4$. We denote
\begin{equation}\label{eq def Ni}
 N_i:=\max\{N_{1,i},N_{2,i}\},\qquad i=1,\ldots,k.
\end{equation}
For $\omega\in\Omega_3$, let $K_i(\omega)$ denote the set
appearing in the definition of $g_{N_i}^i(\omega)$. Then
$(\nu_i)_\omega(K_i(\omega))\geq1-\gamma$.
For $n\geq N_i$, let $C_i(\omega,n,\epsilon_0)$ be an
$(\omega,n,\epsilon_0)$-separated subset of $K_i(\omega)$
of maximal cardinality. Its Bowen balls of radius $\epsilon_0$
cover $K_i(\omega)$, so
\begin{equation}\label{eq def C omega n epsilon0}
 \#C_i(\omega,n,\epsilon_0)
 \geq S_{\nu_i}(\omega,n,\epsilon_0,\gamma)
 \geq e^{n(h_{\nu_i}(F)-\gamma)}.
\end{equation}

Choose $0<\epsilon<\min\{\epsilon_0,\delta/2\}$.
The Lipschitz bound gives
\begin{equation}\label{eq pick of epsilon lemmaproof}
 d_M(x,y)<\epsilon
 \ \Longrightarrow\
 |f_j(\omega,x)-f_j(\omega,y)|<\delta/2,
 \qquad j=1,\ldots,l.
\end{equation}
Let $m=m(\epsilon/4)\geq1$ be a gap supplied by the
$\P$-a.s. fiber specification property, increasing it to $1$
if necessary. By Birkhoff's ergodic theorem, there exists a
$\theta$-invariant $\P$-full-measure set on which
\begin{equation}\label{measure of omega3}
 \lim_{n\to\infty}\frac1n
 \sum_{p=0}^{n-1}1_{\Omega_3}(\theta^p\omega)
 =\P(\Omega_3)>1-\frac{3\xi}{4}.
\end{equation}
Intersect this set with the invariant cores of the full-measure
sets for specification at accuracy $\epsilon/4$ and for the stated
bounds on the test functions. We denote the resulting
$\theta$-invariant $\P$-full-measure set by $\Omega'$.

\emph{Step 2. We construct the measurable function $\bar n$.}

For $\omega\in\Omega'$, let $N_{\Omega_3}(\omega)$ be the least
positive integer such that
\begin{equation}\label{def Nomega3}
 \frac1n\sum_{p=0}^{n-1}1_{\Omega_3}(\theta^p\omega)
 \geq1-\xi
 \quad\text{for every }n\geq N_{\Omega_3}(\omega).
\end{equation}
It is finite by \eqref{measure of omega3}, and measurable because
\[
 \{\omega\in\Omega':N_{\Omega_3}(\omega)\leq q\}
 =
 \bigcap_{n\geq q}
 \left\{\omega\in\Omega':
 \frac1n\sum_{p=0}^{n-1}1_{\Omega_3}(\theta^p\omega)\geq1-\xi\right\}.
\]
Choose a sufficiently large integer $N_*$ such that, for every
integer $n\geq N_*$,
\begin{gather}
 \lfloor\lambda_i n\rfloor\geq N_i\quad(i=1,\ldots,k),
 \qquad n\xi\geq k+1,\label{def of ni}\\
 0\leq n-\sum_{i=1}^k\lfloor\lambda_i n\rfloor
 \leq\frac{\delta n}{4},
 \qquad \frac{2(k-1)m}{n}\leq\xi,\label{pick N6}\\
 \sum_{i=1}^k\lfloor\lambda_i n\rfloor
                  (h_{\nu_i}(F)-\gamma)
 \geq n(h_\nu(F)-2\gamma).\label{eq N03}
\end{gather}
Such an $N_*$ exists since $\lambda_i>0$, all component entropies
are finite, and each rounding error is less than $1$.
We define
\begin{equation}\label{eq pick of n}
 \bar n(\omega):=
 \left\lceil
 \frac{\max\{N_{\Omega_3}(\omega),N_*\}+(k-1)m}{1-\xi}
 \right\rceil.
\end{equation}
The measurability of $\bar n$ follows from that of $N_{\Omega_3}$.

\emph{Step 3. Fix $\omega\in\Omega'$ and $N>\bar n(\omega)$.
We construct $\hat n(\omega)$.}

Define
\begin{equation}\label{eq def ns}
 n_s(\omega):=\lfloor N(1-\xi)-(k-1)m\rfloor,\qquad
 n_i(\omega):=\lfloor\lambda_i n_s(\omega)\rfloor
 \quad(i=1,\ldots,k).
\end{equation}
Then $n_s(\omega)\geq\max\{N_{\Omega_3}(\omega),N_*\}$,
and $n_i(\omega)\geq N_i$ by \eqref{def of ni}.
Since $\omega$ and $N$ are fixed, we suppress their dependence
in the notation below.

Let $l'_0=0$. Inductively, for $i=1,\ldots,k$, let $l_i$
be the first nonnegative integer such that
$\theta^{l'_{i-1}+l_i}\omega\in\Omega_3$, and, for $i<k$, put
\[
 l'_i=l'_{i-1}+l_i+n_i+m.
\]
All these waiting times are finite by \eqref{measure of omega3}.
The $i$th orbit segment starts at time $l'_{i-1}+l_i$,
which we also denote by $s_i$.
Finally, define
\begin{equation}\label{def of hat n}
 \hat n=l'_{k-1}+l_k+n_k
       =\sum_{i=1}^k n_i+(k-1)m+\sum_{i=1}^k l_i.
\end{equation}
To apply \eqref{def Nomega3} at time $\hat n$, first observe that
\[
 0\leq n_s-\sum_{i=1}^k n_i\leq k-1\leq(k-1)m.
\]
Consequently,
\begin{equation}\label{eq hat n geq N0}
 \hat n\geq\sum_{i=1}^k n_i+(k-1)m
 \geq n_s\geq\max\{N_{\Omega_3}(\omega),N_*\}.
\end{equation}
The waiting intervals are disjoint, and all their base iterates
lie outside $\Omega_3$. Applying \eqref{def Nomega3} at time
$\hat n$, we therefore obtain
\begin{equation}\label{proportion of l}
 \sum_{i=1}^k l_i\leq\xi\hat n.
\end{equation}
It follows that
\[
 \hat n\leq
 \frac{\sum_{i=1}^k n_i+(k-1)m}{1-\xi}
 \leq\frac{n_s+(k-1)m}{1-\xi}\leq N.
\]
On the other hand, \eqref{eq def ns} and the rounding bound give
\[
 \hat n\geq\sum_{i=1}^k n_i+(k-1)m
 \geq N(1-\xi)-k\geq N(1-2\xi),
\]
where the last inequality follows from \eqref{def of ni}.
Thus $\hat n=\hat n_N(\omega)$ satisfies the required length bounds.

\emph{Last step. We use orbit gluing to construct the required
separated set.}

For each $i=1,\ldots,k$, we have
$\theta^{l'_{i-1}+l_i}\omega\in\Omega_3$ and $n_i\geq N_i$.
Choose
$x_i\in C_i(\theta^{l'_{i-1}+l_i}\omega,n_i,\epsilon_0)$.
By the definition of $K_i$, these points satisfy
\begin{equation}\label{eq number of C}
 \left|\frac1{n_i}\sum_{p=0}^{n_i-1}
 f_j(\Theta^p(\theta^{l'_{i-1}+l_i}\omega,x_i))
       -\int f_j\,d\nu_i\right|<\delta,
 \qquad j=1,\ldots,l.
\end{equation}
The orbit segments starting at times $l'_{i-1}+l_i$
have lengths $n_i$ and intervening gaps $m+l_{i+1}$.
For every such tuple $(x_1,\ldots,x_k)$, the fiber specification
property gives a point $y=y(x_1,\ldots,x_k)$ such that
\begin{equation}\label{eq lemma y shadowing xi}
 d_{\theta^{l'_{i-1}+l_i}\omega}^{n_i}
 \bigl(x_i,F_\omega^{l'_{i-1}+l_i}y\bigr)<\epsilon/4,
 \qquad i=1,\ldots,k.
\end{equation}
If two tuples differ in their $i$th component, the triangle
inequality and the separation of $C_i$ give
\begin{equation}\label{eq y eta/2 separated}
 d_\omega^{\hat n}
 \bigl(y(x_1,\ldots,x_k),y(x'_1,\ldots,x'_k)\bigr)
 >\epsilon_0-\epsilon/2>\epsilon_0/2.
\end{equation}
In particular, different tuples give different shadowing points.

We next prove that each shadowing point belongs to
$P(\rho,4\delta,\hat n,\omega)$.
Fix $j\in\{1,\ldots,l\}$. On the prescribed orbit segments,
we compare the orbit of $y$ with those of the $x_i$ and then
with the integrals against $\nu_i$. There are
$\sum_{i=1}^k l_i+(k-1)m$ remaining times, at each of which
$|f_j-\int f_j\,d\rho|\leq2$. The triangle inequality gives
\begin{equation}\label{esimation 1}
\begin{split}
 &\left|\sum_{p=0}^{\hat n-1}f_j(\Theta^p(\omega,y))
                   -\hat n\int f_j\,d\rho\right|\\
 &\leq\sum_{i=1}^k
 \left|\sum_{p=0}^{n_i-1}
 \left(f_j(\Theta^{l'_{i-1}+l_i+p}(\omega,y))
       -f_j(\Theta^p(\theta^{l'_{i-1}+l_i}\omega,x_i))\right)\right|\\
 &\quad+\sum_{i=1}^k
 \left|\sum_{p=0}^{n_i-1}
 f_j(\Theta^p(\theta^{l'_{i-1}+l_i}\omega,x_i))
                         -n_i\int f_j\,d\nu_i\right|\\
 &\quad+\left|\sum_{i=1}^k n_i\int f_j\,d\nu_i
                     -\left(\sum_{i=1}^k n_i\right)\int f_j\,d\rho\right|
       +2\left(\sum_{i=1}^k l_i+(k-1)m\right).
\end{split}
\end{equation}
By \eqref{eq lemma y shadowing xi} and
\eqref{eq pick of epsilon lemmaproof}, the first sum is at most
$
 \frac{\delta}{2}\sum_{i=1}^k n_i.
$
By \eqref{eq number of C}, the second sum is at most
$
 \delta\sum_{i=1}^k n_i.
$
For the third term, using $n_i=\lfloor\lambda_i n_s\rfloor$,
$\|f_j\|_\infty\leq1$, and \eqref{pick N6}, we have
\[
\begin{split}
 &\left|\sum_{i=1}^k n_i\int f_j\,d\nu_i
              -\left(\sum_{i=1}^k n_i\right)\int f_j\,d\rho\right|\\
 &\leq n_s\left|\int f_j\,d\nu-\int f_j\,d\rho\right|
            +2\left(n_s-\sum_{i=1}^k n_i\right)
 \leq2n_s\delta.
\end{split}
\]
Since $\sum_{i=1}^k n_i\leq n_s\leq\hat n$, combining these
estimates with \eqref{proportion of l} and \eqref{pick N6} gives
\begin{equation}\label{eq third item}
\begin{split}
 \left|\frac1{\hat n}\sum_{p=0}^{\hat n-1}
 f_j(\Theta^p(\omega,y))-\int f_j\,d\rho\right|
 &\leq\frac{7\delta}{2}+2\xi+\frac{2(k-1)m}{\hat n}\\
 &\leq\frac{7\delta}{2}+3\xi
 \leq\frac{59\delta}{16}<4\delta,
\end{split}
\end{equation}
where the last line uses \eqref{xi estimation}.
Therefore,
\begin{equation}\label{eq y subset P}
 \{y(x_1,\ldots,x_k):
 x_i\in C_i(\theta^{l'_{i-1}+l_i}\omega,n_i,\epsilon_0),\
 i=1,\ldots,k\}
 \subset P(\rho,4\delta,\hat n,\omega).
\end{equation}

By \eqref{eq y eta/2 separated}, these points form an
$(\omega,\hat n,\epsilon_0/2)$-separated set. Hence
\[
\begin{split}
 \#M(\omega,4\delta,\epsilon_0/2,\hat n)
 &\geq\prod_{i=1}^k
 \#C_i(\theta^{l'_{i-1}+l_i}\omega,n_i,\epsilon_0)
 \geq\exp\left(\sum_{i=1}^k n_i(h_{\nu_i}(F)-\gamma)\right)\\
 &\geq\exp\bigl(n_s(h_\nu(F)-2\gamma)\bigr)\geq\exp\bigl(n_s(h_\rho(F)-3\gamma)\bigr),
\end{split}
\]
where we used \eqref{eq def C omega n epsilon0} and \eqref{eq N03}.
Moreover, by \eqref{def of hat n}, \eqref{proportion of l},
and \eqref{pick N6},
\[
 \frac{n_s}{\hat n}\geq
 \frac{\sum_{i=1}^k n_i}{\hat n}
 =1-\frac{\sum_{i=1}^k l_i+(k-1)m}{\hat n}
 \geq1-2\xi.
\]
Since $h_\rho(F)-3\gamma>0$, it follows that
\[
 \frac1{\hat n}\log\#M(\omega,4\delta,\epsilon_0/2,\hat n)
 \geq(1-2\xi)(h_\rho(F)-3\gamma)
 \geq h_\rho(F)-4\gamma,
\]
where the last inequality follows from \eqref{xi estimation}.
This proves \eqref{eq estimate 1/n log M} and completes the proof.
\end{proof}
\begin{proof}[Proof of Lemma \ref{Lemma lower estimates 2}]
Apply Lemma~\ref{lemma measure approxi} with any positive
auxiliary entropy error to obtain a finite convex combination
$\nu=\sum_{i=1}^k\lambda_i\nu_i$, where $\lambda_i>0$,
$\sum_{i=1}^k\lambda_i=1$, and
$\nu_i\in I_\P^e(\Omega\times M)$, such that
\[
 |\nu(f_j)-\rho(f_j)|<\delta,
 \qquad 1\leq j\leq l.
\]
No entropy estimate for the components is needed.

Arguing as in Step~1 of the proof of
Lemma~\ref{Lemma lower estimates}, Birkhoff's ergodic theorem,
disintegration, and Egorov's theorem give integers $N_i\geq1$
and a measurable set $\Omega_3$ with
$\P(\Omega_3)>1-3\xi/4$ such that
\[
 (\nu_i)_\omega(K_i(\omega))\geq\tfrac12
 \quad(\omega\in\Omega_3,\ 1\leq i\leq k),
\]
where
\[
 K_i(\omega)=\left\{x\in M:
 \left|\frac1r S_rf_j(\omega,x)-\nu_i(f_j)\right|<\delta
 \text{ for every }r\geq N_i,\ 1\leq j\leq l\right\}.
\]
In particular, each $K_i(\omega)$ is nonempty.
Choose any $\epsilon_0>0$ and then
$0<\epsilon<\min\{\epsilon_0,\delta/2\}$.
Let $m=m(\epsilon/4)\geq1$ be a specification gap.

Repeating the full-measure-set construction in the preceding
proof with the present choices of $\Omega_3$ and $\epsilon$,
we obtain a $\theta$-invariant $\P$-full-measure set $\Omega'$
on which \eqref{measure of omega3} holds, the test-function
bounds hold along every base iterate, and the fiber
specification property holds at accuracy $\epsilon/4$.

Define $N_{\Omega_3}(\omega)$ as in \eqref{def Nomega3}.
Choose a deterministic integer $N_*$ satisfying
\eqref{def of ni} and \eqref{pick N6}; the entropy condition
\eqref{eq N03} is omitted.
Define $\bar n(\omega)$ by \eqref{eq pick of n}.
Its finiteness and measurability follow as in the preceding proof.

For every $\omega\in\Omega'$ and every $N>\bar n(\omega)$,
construct $n_s,n_i,l_i,l_i',s_i$ and
$\widehat n=\widehat n_N(\omega)$ exactly as in Step~3
of the preceding proof.
The rounding and frequency arguments do not use entropy
and apply to the present range of $\xi$. Hence
\[
 N(1-2\xi)\leq\widehat n\leq N,
 \qquad n_s\leq\widehat n,
 \qquad \sum_{i=1}^k l_i\leq\xi\widehat n.
\]

Choose one point $x_i\in K_i(\theta^{s_i}\omega)$ for each
$i=1,\ldots,k$ and apply specification to these orbit segments.
The resulting point $y$ satisfies
\eqref{eq lemma y shadowing xi}.
Since $n_i\geq N_i$, the definition of $K_i$ gives the same
Birkhoff-average estimates as in the preceding proof.
Therefore, the estimates
\eqref{esimation 1}--\eqref{eq third item}, which use only
the test-function bounds, the approximation of $\rho$,
and the length bounds, give
\[
 \left|\frac1{\widehat n}S_{\widehat n}f_j(\omega,y)
                  -\rho(f_j)\right|
 \leq\frac{7\delta}{2}+3\xi
 \leq\frac{59\delta}{16}<4\delta
 \qquad(1\leq j\leq l),
\]
where we used $\xi\leq\delta/16$.
Thus $y\in P(\rho,4\delta,\widehat n,\omega)$, so this set
is nonempty. Consequently,
\[
 \#M(\omega,4\delta,\epsilon_0/2,\widehat n)\geq1,
 \qquad
 \frac1{\widehat n}
 \log\#M(\omega,4\delta,\epsilon_0/2,\widehat n)\geq0.
\]
Together with the length bounds above, this completes the proof.
\end{proof}

\section{Proof of Theorem \ref{proposition normalized good rate function} and corollaries}
\label{subsection consequences of main results}

\begin{proof}[Proof of Theorem \ref{proposition normalized good rate function}]
Let $\Omega_{\mathrm{LD}}=\Omega_{\rm up}\cap\Omega_{\rm low}$ be the
intersection of the full-measure sets furnished by
Theorems~\ref{thm 2} and~\ref{thm main 1}. On this single set both
bounds hold for every open and every closed set. By
Lemma~\ref{lemma lower energy finite entropy}, all invariant measures
with marginal $\P$ have finite fiber entropy.
Lemma~\ref{lemma upper semi-cotinuous} makes the energy integral
continuous, and the fiber entropy map is upper semicontinuous by
hypothesis. Hence $I_T$ is lower semicontinuous on
$I_\P(\Omega\times M)$. By Remark~\ref{remark compactness} and
Lemma~\ref{lemma lemma diffrent converges}, this set is compact in both
narrow topologies and closed in $\operatorname{Pr}(\Omega\times M)$.
Its extension by $+\infty$ is therefore lower semicontinuous on the
whole probability-measure space. For every finite $a$,
\[
 \{\rho\in \operatorname{Pr}(\Omega\times M):I_T(\rho)\leq a\}
 =\{\rho\in I_\P(\Omega\times M):I_T(\rho)\leq a\}
\]
is compact.

Applying the lower bound to the whole space gives
\begin{equation}\label{eq normalization upper half}
 0\geq\sup_{\rho\in I_\P(\Omega\times M)}
 \left\{h_\rho(F)-\int_{\Omega\times M}T\,d\rho\right\}.
\end{equation}
The upper bound, applied to the same closed convex set, gives the
reverse inequality. Consequently,
\begin{equation}\label{eq free energy normalized}
 \sup_{\rho\in I_\P(\Omega\times M)}
 \left\{h_\rho(F)-\int_{\Omega\times M}T\,d\rho\right\}=0.
\end{equation}
Thus $I_T\geq0$ and $\inf I_T=0$. Together with compact sublevel sets,
this proves that $I_T$ is a good rate function and completes the proof.
\end{proof}

\begin{proof}[Proof of Corollary \ref{corollary quenched level one LDP}]
Since $g$ is bounded and continuous,  the mapping $\mathcal{A}_g:\operatorname{Pr}(\Omega\times M)\to \mathbb{R}$ defined by $\mathcal{A}_g(\rho):=\int_{\Omega\times M}g\,d\rho$ is continuous with respect to $(\operatorname{Pr}(\Omega\times M),\tau)$.
We first prove that $J_g$ is a good rate function.  Since $I_T\geq0$,
we have $J_g\geq0$.  For every
$b\in[0,+\infty)$, we claim that
\begin{equation}\label{eq compact level set J g}
  \{a\in\mathbb R:J_g(a)\leq b\}
  =\mathcal A_g\bigl(\{\rho\in \operatorname{Pr}(\Omega\times M):I_T(\rho)\leq b\}\bigr).
\end{equation}
The direction $``\supseteq"$ is clear. We only need to prove the other direction. For any $a$ such that $J_g(a)\leq b$, by definition, there exists a sequence $\rho_m\in I_\P(\Omega\times M)$ with $\int g\,d\rho_m=a$ and $I_T(\rho_m)=\int T\,d\rho_m-h_{\rho_m}(F)\leq b+\frac{1}{m}$. Since $I_\P(\Omega\times M)$ is compact with respect to $(\operatorname{Pr}(\Omega\times M),\tau)$, there exists a subsequence of $\rho_m$, denoted by $\rho_{m_k}$, converging to a measure $\rho^*\in I_\P(\Omega\times M)$ with respect to $(\operatorname{Pr}(\Omega\times M),\tau)$. Since $g\in C_b(\Omega\times M)$, we obtain $\int g\,d\rho^*=\lim_{k\to+\infty}\int g\,d\rho_{m_k}=a$. Since $I_T$ is a good rate function, by lower semicontinuity, we obtain
\begin{equation*}
  I_T(\rho^*)\leq \liminf_{k\to+\infty}I_T(\rho_{m_k})\leq b.
\end{equation*}
Hence $a\in \mathcal A_g\bigl(\{\rho\in \operatorname{Pr}(\Omega\times M):I_T(\rho)\leq b\}\bigr)$. The claim is proved.

Fix $\omega\in\Omega_{\mathrm{LD}}$. Both level--2 bounds hold for every open and every closed set on this same fiber.
If $O\subset\mathbb R$ is open, then the subset
$\mathcal{A}_g^{-1}(O)=\{\rho\in \operatorname{Pr}(\Omega\times M):\ \int g\, d\rho\in O\}$ is open in $(\operatorname{Pr}(\Omega\times M),\tau)$.
By \eqref{eq LDP lower I T},
\begin{equation}\label{eq open O}
\begin{split}
 &\liminf_{n\to+\infty}\frac{1}{n}
  \log\nu_\omega\{x\in M:\ \int_{\Omega\times M} g\,d(\E_n(\omega,x))\in O \}\\
  =&\liminf_{n\to+\infty}\frac{1}{n}
  \log\nu_\omega\{x\in M:\ \mathcal{E}_n(\omega,x)\in \mathcal{A}_g^{-1}(O)\}\\
  \geq&-\inf_{\rho\in\mathcal A_g^{-1}(O)}I_T(\rho)
   =-\inf_{a\in O}J_g(a).
\end{split}
\end{equation}  If $K\subset\mathbb R$
is closed, then $\mathcal A_g^{-1}(K)$ is also closed.  By
\eqref{eq LDP upper I T},
\begin{equation}\label{eq closed K}
\begin{split}
  &\limsup_{n\to+\infty}\frac{1}{n}\log\nu_\omega\{x\in M:\ \int_{\Omega\times M} g\,d(\E_n(\omega,x))\in K \}\\
  =&\limsup_{n\to+\infty}\frac{1}{n}
  \log\nu_\omega\{x\in M:\ \mathcal{E}_n(\omega,x)\in \mathcal{A}_g^{-1}(K)\}\\
  \leq&-\inf_{\rho\in\mathcal A_g^{-1}(K)}I_T(\rho)
   =-\inf_{a\in K}J_g(a),
\end{split}
\end{equation}
which completes the proof.

\end{proof}

\begin{proof}[Proof of Corollary \ref{corollary quenched Laplace principle}]
Fix $\omega\in\Omega_{\mathrm{LD}}$ and $\Phi\in C_b(\operatorname{Pr}(\Omega\times M))$. We first prove the lower bound in \eqref{eq quenched Laplace principle}.
For $n\in\mathbb N$, define $\mathbb{Q}_{n,\omega}:=(\E_n(\omega,\cdot))_*\nu_\omega \in \operatorname{Pr}(\operatorname{Pr}(\Omega\times M))$, that is, for every Borel set $U\subset \operatorname{Pr}(\Omega\times M)$,
\begin{equation*}
  \mathbb{Q}_{n,\omega}(U):=\nu_\omega\{x\in M:\ \E_n(\omega,x)\in U\}.
\end{equation*}
For $n\in\mathbb N$ and $\omega\in\Omega$, it follows that
\begin{equation}\label{eq exponential pushforward identity}
  \int_M e^{n\Phi(\E_n(\omega,x))}d\nu_\omega(x)
  =\int_{\operatorname{Pr}(\Omega\times M)}e^{n\Phi(\rho)}d\mathbb Q_{n,\omega}(\rho).
\end{equation}
Set
\begin{equation*}
  L_\Phi:=
  \sup_{\rho\in \operatorname{Pr}(\Omega\times M)}
  \{\Phi(\rho)-I_T(\rho)\}.
\end{equation*}
Since $\Phi$ is bounded and $\inf_{\rho\in \operatorname{Pr}(\Omega\times M)} I_T(\rho)=0$, the number $L_\Phi$ is finite.
Choose $\rho_k\in \operatorname{Pr}(\Omega\times M)$ such that
\begin{equation}\label{eq approximate Laplace supremum}
  \Phi(\rho_k)-I_T(\rho_k)>L_\Phi-\frac1k,
  \qquad k\in\mathbb N.
\end{equation}
Since $\Phi$ is continuous,
\begin{equation*}
  G_k:=\left\{\rho\in \operatorname{Pr}(\Omega\times M):
  \Phi(\rho)>\Phi(\rho_k)-\frac1k\right\}
\end{equation*}
is an open neighborhood of $\rho_k$.  By
\eqref{eq exponential pushforward identity},
\begin{equation*}
  \int_M e^{n\Phi(\E_n(\omega,x))}d\nu_\omega(x)\geq \int_{G_k}e^{n\Phi(\rho)}d\mathbb Q_{n,\omega}(\rho)
  \geq e^{n(\Phi(\rho_k)-1/k)}\mathbb Q_{n,\omega}(G_k).
\end{equation*}
Applying \eqref{eq LDP lower I T} on the fixed fiber to each $G_k$, we obtain
\begin{equation*}
\begin{split}
  \liminf_{n\to+\infty}\frac1n
  \log\int_M e^{n\Phi(\E_n(\omega,x))}d\nu_\omega(x)
  \geq & \Phi(\rho_k)-\frac1k+\liminf_{n\to+\infty}\frac{1}{n}\log(\mathbb{Q}_{n,\omega}(G_k))\\
  \geq & \Phi(\rho_k)-\frac1k-\inf_{\rho\in G_k}I_T(\rho),\mbox{ by \eqref{eq LDP lower I T}}\\
  \geq &\Phi(\rho_k)-\frac1k-I_T(\rho_k)
  >L_\Phi-\frac2k.
\end{split}
\end{equation*}
Letting $k\to+\infty$ gives
\begin{equation}\label{eq Laplace lower proof}
  \liminf_{n\to+\infty}\frac1n
  \log\int_M e^{n\Phi(\E_n(\omega,x))}d\nu_\omega(x)
  \geq L_\Phi.
\end{equation}

We now prove  $``\leq "$ in \eqref{eq quenched Laplace principle}.  Write
\begin{equation*}
  a:=\inf_{\rho\in \operatorname{Pr}(\Omega\times M)}\Phi(\rho),
  \qquad
  b:=\sup_{\rho\in \operatorname{Pr}(\Omega\times M)}\Phi(\rho).
\end{equation*}
If $a=b$, the assertion \eqref{eq quenched Laplace principle} follows immediately from
\eqref{eq normalization I T}.  Suppose that $a<b$.  For every
$k\in\mathbb N$, choose
\begin{equation*}
  a=t_{k,0}<t_{k,1}<\cdots<t_{k,N_k}=b
\end{equation*}
such that $t_{k,i}-t_{k,i-1}\leq1/k$ for every $i$.  Define
\begin{equation*}
  K_{k,i}:=\{\rho\in \operatorname{Pr}(\Omega\times M):
  t_{k,i-1}\leq\Phi(\rho)\leq t_{k,i}\}.
\end{equation*}
The sets $K_{k,i}$ are closed and, for each $k$, they cover
$\operatorname{Pr}(\Omega\times M)$.  Therefore
\begin{equation*}
\begin{split}
  \int_{\operatorname{Pr}(\Omega\times M)}e^{n\Phi(\rho)}d\mathbb Q_{n,\omega}(\rho)
  &\leq
  \sum_{i=1}^{N_k}e^{nt_{k,i}}
  \mathbb Q_{n,\omega}(K_{k,i}).
\end{split}
\end{equation*}
Apply \eqref{eq LDP upper I T} to the closed sets $K_{k,i}$ on the fixed fiber.
Since $N_k$ does not
depend on $n$, we obtain
\begin{equation*}
\begin{split}
  \limsup_{n\to+\infty}\frac1n
  \log\int_{\operatorname{Pr}(\Omega\times M)}e^{n\Phi(\rho)}d\mathbb Q_{n,\omega}(\rho)\leq &  \max_{1\leq i\leq N_k}
  \left\{t_{k,i}+ \limsup_{n\to+\infty}\frac1n
  \log\mathbb{Q}_{n,\omega}(K_{k,i}) \right\}\\
 \leq &
  \max_{1\leq i\leq N_k}
  \left\{t_{k,i}-\inf_{\rho\in K_{k,i}}I_T(\rho)\right\}.
\end{split}
\end{equation*}
For every $\rho\in K_{k,i}$, we have
$t_{k,i}\leq\Phi(\rho)+1/k$.  Hence
\begin{equation*}
  t_{k,i}-\inf_{\rho\in K_{k,i}}I_T(\rho)
  \leq\frac1k+
  \sup_{\rho\in K_{k,i}}\{\Phi(\rho)-I_T(\rho)\}
  \leq\frac1k+L_\Phi.
\end{equation*}
Letting $k\to\infty$ and combining the resulting upper bound with
\eqref{eq Laplace lower proof} proves \eqref{eq quenched Laplace principle}.
Taking $\Phi(\rho)=\int g\,d\rho$, where $g\in C_b(\Omega\times M)$,
gives the stated pressure-type formula. Since $\omega\in\Omega_{\mathrm{LD}}$
was arbitrary, this same set works for every such $\Phi$ and $g$.
\end{proof}
\section*{}
\subsection*{Acknowledgments}
The paper was partially supported by the National Natural Science Foundation of China (Nos. 12401230,~12171175 and~12371197),
Jiangsu Provincial Scientific Research Center of Applied Mathematics (BK20233002), and the Qinglan Project of Jiangsu Province of China.
\subsection*{Declaration of AI use}
The main mathematical ideas and proof strategies in this work were developed by the authors. During the preparation of the manuscript, the authors used ChatGPT-6 to assist with language editing, literature searches, and the refinement and critical examination of details in the mathematical arguments. Such assistance does not replace the authors’ independent verification. The authors take full responsibility for the mathematical results and proofs, the accuracy of the references, and the entire content of the manuscript.


\begin{thebibliography}{10}

\bibitem{ArnoldBook1998}
{\sc Arnold, L.}
\newblock {\em Random dynamical systems}.
\newblock Springer Monographs in Mathematics. Springer-Verlag, Berlin, 1998.

\bibitem{Bogachev}
{\sc Bogachev, V.~I.}
\newblock {\em Measure theory. {V}ol. {I}, {II}}.
\newblock Springer-Verlag, Berlin, 2007.

\bibitem{ChungNakanoWittsten2023}
{\sc Chung, Y.~M., Nakano, Y., and Wittsten, J.}
\newblock Quenched limit theorems for random {$U(1)$} extensions of expanding
  maps.
\newblock {\em Discrete Contin. Dyn. Syst. 43}, 1 (2023), 338--377.

\bibitem{ClimenhagaThompsonYamamoto2017}
{\sc Climenhaga, V., Thompson, D.~J., and Yamamoto, K.}
\newblock Large deviations for systems with non-uniform structure.
\newblock {\em Trans. Amer. Math. Soc. 369}, 6 (2017), 4167--4192.

\bibitem{CommanRiveraLetelier2011}
{\sc Comman, H., and Rivera-Letelier, J.}
\newblock Large deviation principles for non-uniformly hyperbolic rational
  maps.
\newblock {\em Ergodic Theory Dynam. Systems 31}, 2 (2011), 321--349.

\bibitem{Hans02}
{\sc Crauel, H.}
\newblock {\em Random probability measures on {P}olish spaces}, vol.~11 of {\em
  Stochastics Monographs}.
\newblock Taylor \& Francis, London, 2002.

\bibitem{DragicevicFroylandGonzalezTokmanVaienti2018ASIP}
{\sc Dragi\v{c}evi\'{c}, D., Froyland, G., Gonz\'{a}lez-Tokman, C., and
  Vaienti, S.}
\newblock Almost sure invariance principle for random piecewise expanding maps.
\newblock {\em Nonlinearity 31}, 5 (2018), 2252--2280.

\bibitem{DragicevicFroylandGonzalezTokmanVaienti2018}
{\sc Dragi\v{c}evi\'{c}, D., Froyland, G., Gonz\'{a}lez-Tokman, C., and
  Vaienti, S.}
\newblock A spectral approach for quenched limit theorems for random expanding
  dynamical systems.
\newblock {\em Comm. Math. Phys. 360}, 3 (2018), 1121--1187.

\bibitem{DragicevicFroylandGonzalezTokmanVaienti2020}
{\sc Dragi\v{c}evi\'{c}, D., Froyland, G., Gonz\'{a}lez-Tokman, C., and
  Vaienti, S.}
\newblock A spectral approach for quenched limit theorems for random hyperbolic
  dynamical systems.
\newblock {\em Trans. Amer. Math. Soc. 373}, 1 (2020), 629--664.

\bibitem{DragicevicHafouta2020}
{\sc Dragi\v{c}evi\'{c}, D., and Hafouta, Y.}
\newblock Limit theorems for random expanding or {A}nosov dynamical systems and
  vector-valued observables.
\newblock {\em Ann. Henri Poincar\'{e} 21}, 12 (2020), 3869--3917.

\bibitem{DragicevicHafouta2021ASIP}
{\sc Dragi\v{c}evi\'{c}, D., and Hafouta, Y.}
\newblock Almost sure invariance principle for random dynamical systems via
  {Gou\"ezel}'s approach.
\newblock {\em Nonlinearity 34}, 10 (2021), 6773--6798.

\bibitem{Kifer2000specification}
{\sc Gundlach, V.~M., and Kifer, Y.}
\newblock Expansiveness, specification, and equilibrium states for random
  bundle transformations.
\newblock {\em Discrete Contin. Dyn. Syst. 6}, 1 (2000), 89--120.

\bibitem{Hafouta2022}
{\sc Hafouta, Y.}
\newblock Limit theorems for random non-uniformly expanding or hyperbolic maps
  with exponential tails.
\newblock {\em Ann. Henri Poincar\'{e} 23}, 1 (2022), 293--332.

\bibitem{Hafouta2024MixingBase}
{\sc Hafouta, Y.}
\newblock Large deviations, moment estimates and almost sure invariance
  principles for skew products with mixing base maps and expanding-on-average
  fibers.
\newblock {\em Ergodic Theory Dynam. Systems 44}, 1 (2024), 118--158.

\bibitem{HLLScichina}
{\sc Huang, W., Lian, Z., and Lu, K.}
\newblock Dynamical complexity of {A}nosov systems driven by a quasi-periodic
  force.
\newblock {\em Sci. China Math. 68}, 1 (2025), 89--136.

\bibitem{KellerNowicki1992}
{\sc Keller, G., and Nowicki, T.}
\newblock Spectral theory, zeta functions and the distribution of periodic
  points for {Collet--Eckmann} maps.
\newblock {\em Comm. Math. Phys. 149}, 1 (1992), 31--69.

\bibitem{Kifer1990}
{\sc Kifer, Y.}
\newblock Large deviations in dynamical systems and stochastic processes.
\newblock {\em Trans. Amer. Math. Soc. 321}, 2 (1990), 505--524.

\bibitem{Kifer1991RandomExpanding}
{\sc Kifer, Y.}
\newblock Large deviations for random expanding maps.
\newblock In {\em Lyapunov Exponents}, L.~Arnold, H.~Crauel, and J.-P. Eckmann,
  Eds., vol.~1486 of {\em Lecture Notes in Mathematics}. Springer, Berlin,
  1991, pp.~178--186.

\bibitem{Kifer1996RandomEnvironments}
{\sc Kifer, Y.}
\newblock Perron--{F}robenius theorem, large deviations, and random
  perturbations in random environments.
\newblock {\em Math. Z. 222}, 4 (1996), 677--698.

\bibitem{Kifer1998LimitTheorems}
{\sc Kifer, Y.}
\newblock Limit theorems for random transformations and processes in random
  environments.
\newblock {\em Trans. Amer. Math. Soc. 350}, 4 (1998), 1481--1518.

\bibitem{Kifervariationalprinciple}
{\sc Kifer, Y.}
\newblock On the topological pressure for random bundle transformations.
\newblock In {\em Topology, ergodic theory, real algebraic geometry}, vol.~202
  of {\em Amer. Math. Soc. Transl. Ser. 2}. Amer. Math. Soc., Providence, RI,
  2001, pp.~197--214.

\bibitem{KL06}
{\sc Kifer, Y., and Liu, P.-D.}
\newblock Random dynamics.
\newblock In {\em Handbook of dynamical systems. {V}ol. 1{B}}. Elsevier B. V.,
  Amsterdam, 2006, pp.~379--499.


\bibitem{Liu2025variatonalprinciple}
{\sc Liu, N., and Liu, X.}
\newblock On the variational principle for a class of skew product
  transformations.
\newblock {\em J. Differential Equations 442\/} (2025), Paper No. 113498, 70
  pp.

\bibitem{LiuLiu2026Irregular}
{\sc Liu, N., and Liu, X.}
\newblock Irregular sets for a class of skew product transformations.
\newblock {\em Discrete Contin. Dyn. Syst. 51\/} (2026), 128--160.

\bibitem{LiuZhao2003}
{\sc Liu, P.-D., and Zhao, Y.}
\newblock Large deviations in random perturbations of {A}xiom {A} basic sets.
\newblock {\em J. London Math. Soc. (2) 68}, 1 (2003), 148--164.

\bibitem{LiuMa2026}
{\sc Liu, X., and Ma, X.}
\newblock Bowen's {S}pecification {P}roperty for {I}nduced {D}ynamics on
  {S}paces of {R}andom {P}robability {M}easures.
\newblock {\em Qual. Theory Dyn. Syst. 25}, 2 (2026), Paper No. 69.

\bibitem{PfisterSullivan2005}
{\sc Pfister, C.-E., and Sullivan, W.~G.}
\newblock Large deviations estimates for dynamical systems without the
  specification property. {A}pplications to the {$\beta$}-shifts.
\newblock {\em Nonlinearity 18}, 1 (2005), 237--261.

\bibitem{PfisterSullivan2018}
{\sc Pfister, C.-E., and Sullivan, W.~G.}
\newblock Weak {G}ibbs measures and large deviations.
\newblock {\em Nonlinearity 31}, 1 (2018), 49--53.

\bibitem{ReyBelletYoung2008}
{\sc Rey-Bellet, L., and Young, L.-S.}
\newblock Large deviations in non-uniformly hyperbolic dynamical systems.
\newblock {\em Ergodic Theory Dynam. Systems 28}, 2 (2008), 587--612.

\bibitem{GTM180}
{\sc Srivastava, S.~M.}
\newblock {\em A Course on Borel Sets}, vol.~180 of {\em Graduate Texts in
  Mathematics}.
\newblock Springer, New York, 1998.

\bibitem{Waltersergodic}
{\sc Walters, P.}
\newblock {\em An introduction to ergodic theory}, vol.~79 of {\em Graduate
  Texts in Mathematics}.
\newblock Springer-Verlag, New York-Berlin, 1982.

\bibitem{YangChenLinZhou2025}
{\sc Yang, K., Chen, E., Lin, Z., and Zhou, X.}
\newblock Weighted topological entropy of random dynamical systems.
\newblock {\em Discrete Contin. Dyn. Syst. 45}, 8 (2025), 2688--2734.

\bibitem{Young1990}
{\sc Young, L.-S.}
\newblock Some large deviation results for dynamical systems.
\newblock {\em Trans. Amer. Math. Soc. 318}, 2 (1990), 525--543.

\bibitem{zhutwonotes}
{\sc Zhu, Y.~J.}
\newblock Two notes on measure-theoretic entropy of random dynamical systems.
\newblock {\em Acta Math. Sin. (Engl. Ser.) 25}, 6 (2009), 961--970.

\end{thebibliography}
\end{document}